\documentclass[11pt]{article}
\newcommand{\keywords}[1]{\textit{Keywords:}\quad #1}

\usepackage{amsfonts}
\usepackage{amssymb}
\usepackage{mathrsfs}
\usepackage{latexsym}
\usepackage{pgf,tikz}
\usetikzlibrary{arrows}
\usetikzlibrary[patterns]
\usepackage{url}

\usepackage[all]{nowidow}

\usepackage{setspace}

\usepackage{graphicx}
\usepackage{enumitem}
\usepackage{lipsum} 
\usepackage{subfig}
\usepackage{amsmath,amssymb,amsthm}
\usepackage{varioref} 
\usepackage{bbold}
\usepackage{scalerel}
\usepackage[makeroom]{cancel}
\usepackage{mathtools}
\usepackage{float}

\usepackage[text={6.5in,9in}]{geometry}
\theoremstyle{plain}

\newcommand{\Q}{\mathcal{Q}}

\newcommand{\T}{\mathcal{T}}

\newcommand{\BBQ}{\mbox{$\mathbb{Q}$}}
\newcommand{\BBR}{\mbox{$\mathbb{R}$}}

\newcommand{\bsigma}{\mbox{\boldmath{$\sigma$}}}

\newcommand{\btau}{\mbox{\boldmath{$\tau$}}}

\newcommand{\bvarepsilon}{\mbox{\boldmath{$\varepsilon$}}}

\newcommand{\half}{\mbox{$\frac{1}{2}$}}

\newcommand{\divV}{\nabla \cdot}

\newcommand{\ldisc}{[\hspace{-0.3ex}[}   
\newcommand{\rdisc}{]\hspace{-0.3ex}]}   
\newcommand{\ldiscV}{\lfloor \hspace{-0.6ex} \lfloor}  
\newcommand{\rdiscV}{\rfloor \hspace{-0.6ex} \rfloor}

\newcommand{\lrcb}[1]{ \left\lbrace #1 \right\rbrace}
\newcommand{\lravg}[1]{ \{ \hspace{-4pt} \{  #1 \} \hspace{-4pt} \} } 
\newcommand{\lrsb}[1]{ \left[ #1 \right]}
\newcommand{\lrb}[1]{ \left( #1 \right)}

\newcommand{\lrjmp}[1]{ \ldisc #1 \rdisc}
\newcommand{\lrjmpV}[1]{\ldiscV #1 \rdiscV}

\newcommand{\ubf}{\mathbf{f}}
\newcommand{\ubg}{\mathbf{g}}
\newcommand{\ubh}{\mathbf{h}}

\newcommand{\ubn}{\mathbf{n}}

\newcommand{\ubu}{\mathbf{u}}
\newcommand{\ubv}{\mathbf{v}}

\newcommand{\idten}{\mathbf{1}}

\newcommand{\elBound}[1]{\partial \Omega_ #1}

\newcommand{\DBound}{\Gamma_D}
\newcommand{\NBound}{\Gamma_N}
\newcommand{\intBound}{\Gamma_{int}}

\newcommand{\partition}{\T_h}

\newcommand{\elDom}{\Omega_e}

\newcommand{\Ltd}{L^2 \lrb{\Omega}}

\newcommand{\Hoe}{H^1 \lrb{\elDom}}

\newcommand{\Hop}{H^1 \lrb{\partition}}

\newcommand{\sumalle}{\sum_{\elDom \in \partition}}

\newcommand{\sumDE}{\sum_{E \in \DBound}}
\newcommand{\sumNE}{\sum_{E \in \NBound}}
\newcommand{\sumiDE}{\sum_{E \in \Gamma_{iD}}}

\newcommand{\intEl}{\int_{\elDom}}

\newcommand{\epsu}{\bvarepsilon \lrb{\ubu} }

\newcommand{\epsv}{\bvarepsilon \lrb{\ubv} }

\newcommand{\sigu}{\bsigma \lrb{\ubu}}

\newcommand{\defined}{\colon \hspace{-5pt}=}

\title{Convergence in the incompressible limit of new discontinuous Galerkin methods with general quadrilateral and hexahedral elements
}

\author{Beverley~J. Grieshaber\footnotemark[2]\
\and Andrew~T. McBride\footnotemark[3]\
\and B.~Daya Reddy\footnotemark[4]}
\date{}

\begin{document}
\maketitle
\renewcommand{\thefootnote}{\fnsymbol{footnote}}
\footnotetext[2]{Centre for Research in Computational and Applied Mechanics, and Department of Mathematics and Applied Mathematics, University of Cape Town, 7701 Rondebosch, South Africa; beverley.grieshaber@uct.ac.za; Tel. +27 21 650 3817; Fax. +27 21 685 2281}
\footnotetext[3]{Glasgow Computational Engineering Centre,
The University of Glasgow, Glasgow G12 8QQ}
\footnotetext[4]{Centre for Research in Computational and Applied Mechanics, and Department of Mathematics and Applied Mathematics, University of Cape Town, 7701 Rondebosch, South Africa}

\begin{abstract}
Standard low-order finite elements, which perform well for problems involving compressible elastic materials, are known to under-perform when nearly incompressible materials are involved, commonly exhibiting the locking phenomenon.
Interior penalty (IP) discontinuous Galerkin methods have been shown to circumvent locking when simplicial elements are used. The same IP methods, however, result in locking on meshes of quadrilaterals. The authors have shown in earlier work that  under-integration of specified terms in the IP formulation eliminates the locking problem for rectangular elements. Here it is demonstrated through an extensive numerical investigation that the effect of using under-integration carries over successfully to meshes of more general quadrilateral elements, as would likely be used in practical applications, and results in accurate displacement approximations. Uniform convergence with respect to the compressibility parameter is shown numerically. Additionally, a stress approximation obtained here by postprocessing shows good convergence in the incompressible limit.

\end{abstract}

\section{Introduction}

The finite element method is well established as a method for solving boundary value problems approximately. Numerical implementations are generally supported by rigorous analyses of the method, certainly for linear problems, and for an increasingly wide range of nonlinear problems. Despite this, significant challenges remain.

In the context of solid mechanics, and in particular in problems for elastic materials, the standard Galerkin (SG) finite element method, while performing very well for compressible materials, may exhibit the phenomenon known as ``volumetric locking" for materials that are nearly incompressible, if low-order (linear or bilinear, or trilinear in three dimensions) elements are used. 
Manifesting particularly in the case of bending-dominated problems, this pathological behaviour results from the too-severe constraint placed on the solution by the incompressibility condition.
The adverse effect of the degree of compressibility on the performance of the SG method may nevertheless also be seen in poor displacement approximations that are not specifically of the locking type. 

The problem may be circumvented by the use of high-order elements. Low-order
elements remain an attractive option, though, and for this reason various extensions to the SG method have been studied, and shown to be effective in remedying locking when low-order approximations for the displacement are used. One class of extensions is mixed methods (see, for example, \cite{Brezzi1991}).
Related to the pressure-displacement mixed method is the method of selective reduced integration (SRI), also effective in producing locking-free results. Finally, discontinuous Galerkin (DG) methods, specifically the range of interior penalty (IP) DG methods, have been used effectively with low-order elements, within a limited scope.

Various DG mixed methods have been proven to be robust for near-incompressible isotropic elasticity,  allowing overall for meshes of quadrilateral, hexahedral, as well as simplicial elements (see \cite{Hansbo2002,Cockburn2006,Houston2006,Wihler2012}). 
In contrast, DG primal formulations have been established as having optimal performance independent of material parameters for meshes of triangular elements (\cite{Wihler2004,Wihler2006,Hansbo2002,Hansbo2003,Bridgeman2011}), but not for meshes of quadrilateral elements. 

In a numerical investigation, Liu et al.\ \cite{Liu2009} consider the matter of locking in the context of low-order hexahedral elements. This work, on a specific benchmark problem, shows the under-performance of the SG method and the superiority of three well-known primal IP methods, NIPG, SIPG and IIPG (Nonsymmetric, Symmetric and Incomplete Interior Penalty Galerkin methods, developed in \cite{Riviere1999,Riviere2000,Wheeler1978,Dawson2004}), as well as of the method of Oden, Babu\u{s}ka and Baumann \cite{Oden1998} (known as OBB), a penalty-free version of NIPG. 
However, the conditions of the benchmark problem are not those leading to the severest form of locking, and no accompanying analysis or convergence data is included. Therefore, while the results are positive, questions remain about the scope of the superior performance of the IP methods with non-simplicial elements.

In \cite{Grieshaber2015}, the authors showed using several numerical examples in two and three dimensions that the three IP methods do in fact produce poor approximations, and notably locking-type behaviour, when quadrilateral/hexahedral elements are used. A new method, in which selected edge terms of the IP formulation are under-integrated, circumvents the problem, as shown through an analytical proof that the new method is locking-free for rectangular elements, and through numerical examples that demonstrate the optimal performance of this formulation. 

While the technique of under-integration (or SRI) has long been used in other contexts to eliminate volume locking, it has until recently typically been used on integrals on element domains, while in the formulation of \cite{Grieshaber2015} it is used in integrals on element edges only. Hansbo and Larson \cite{Hansbo2003} use under-integration similarly on an edge-based stabilization term in their formulation for triangular elements, relaxing what would be otherwise be a severely constraining term. In other applications, edge under-integration is used in \cite{Grieshaber20inpress} with simplicial elements to eliminate extensional locking within transversely isotropic elasticity, and in \cite{Bayat2018} this technique is used to circumvent shear locking in beams.

The analysis in \cite{Grieshaber2015} of the new IP formulation accommodates both essential and natural boundary conditions, and the numerical examples incorporate both, as would be expected in a realistic, practical model. However, both the theoretical and numerical analyses presented are concerned with rectangular elements and simple domains. A single example in that paper shows the method locking-free with quadrilaterals as well, but is not conclusive regarding the general case.
In a more recent computational paper, Bayat et al.\ \cite{Bayat2018} have similarly shown several variants of  IIPG to be volumetrically locking-free when selected edge terms are under-integrated. However, these authors consider limited test cases only and their results are therefore inconclusive regarding the general case.

For the method of \cite{Grieshaber2015} to be of practical value in solving a broad range of boundary value problems, it is necessary to establish its robustness when general quadrilaterals or hexahedrals are employed, allowing for its use on a variety of domain shapes and meshes (particularly unstructured).

This paper seeks to address this issue computationally through an extensive numerical investigation using non-rectangular elements with a variety of model problems. We systematically compare the performance of the new method, particularly in the incompressible limit, as meshes with decreasing element shape regularity are used. We demonstrate that the formulation of \cite{Grieshaber2015} is effective in alleviating locking and generating accurate approximations in the general case.

In practice the stress field generated from the displacement approximation by postprocessing is often of interest. As a second component of this work we therefore study the accuracy of the stress field approximation obtained from the new IP methods, considering both error convergence rates and approximation quality at individual refinement levels.

Following this introduction, the boundary value problem and DG framework and formulation are presented in \S \ref{preliminaries}. Section \ref{numerics} gives the numerical results for both the displacement and the stress approximations, before the conclusion in \S \ref{conc}.

\section{Preliminaries}\label{preliminaries}

\subsection{The boundary value problem of linear elasticity}

Let a homogeneous, isotropic, linear elastic body occupy the bounded domain $\Omega \subset \BBR^d$ ($d=2,3$), with the Lipschitz boundary $\partial \Omega$ consisting of a Dirichlet portion, $\DBound$, of positive measure, and a Neumann portion, $\NBound$, such that $\DBound \cap \NBound = \emptyset \text{ and } \DBound \cup \NBound = \partial \Omega$, and with outward unit normal $\ubn$.

A body force $\ubf \in \lrsb{\Ltd}^d$ is applied, with prescribed displacement $\ubg \in \lrsb{L^2 \lrb{\DBound}}^d$ on $\DBound$ and prescribed traction $\ubh \in \lrsb{L^2 \lrb{\NBound}}^d$ on $\NBound$. The resultant displacement is $\ubu$, and the strain $\bvarepsilon$ is expressed as a tensor defined in index notation as 
\begin{align}
&\bvarepsilon_{ij} (\ubu) \defined \half \lrb{u_{i,j} + u_{j,i}}, &1 \leq i,j \leq d. \notag 
\end{align}
The stress $\sigu$ is related to the strain via the constitutive law
\begin{align}\label{const_rel}
\sigu \defined 2 \mu \epsu + \lambda \text{tr}\, \epsu \idten
= 2 \mu \epsu + \lambda \lrb{\divV \ubu} \idten, 
\end{align}
where $\lambda$ and the shear modulus $\mu$ are known as the Lam\'e parameters, and $\idten$ is the second-order identity tensor in $\BBR^d$.

The governing equation of the system is the equilibrium equation
\begin{subequations} \label{strong_eqns}
\begin{align} \label{gov_eqn}
- \text{div } \sigu &= \ubf		&\text{in } \Omega, 
\end{align}
and the boundary conditions are
\begin{align} 
\ubu &= \ubg 					&\text{on } \DBound, \label{bc_eqnD}\\
\sigu \ubn &= \ubh 				&\text{on } \NBound. \label{bc_eqnN}
\end{align}
\end{subequations}

The Lam\'e parameters $\lambda$ and $\mu$ are assumed to be positive, and can be expressed in terms of the Young's modulus, $E$, and Poisson's ratio, $\nu$, by
\begin{align}
\lambda = 	\frac{E \nu}{\lrb{1+ \nu} \lrb{1 - 2 \nu}},	
\ \ \ \mu = \frac{E}{2 \lrb{1 + \nu}}. \notag
\end{align}
As $\nu \rightarrow \half$, which corresponds to the incompressible limit, so $\lambda \rightarrow \infty$. 

\subsection{The discontinuous Galerkin framework}

The domain $\Omega \subset \BBR^d\ \lrb{d = 2,\,3},$ is partitioned into a mesh of regular elements $\elDom$, where $\partition = \lrcb{\elDom}$. 
The outward unit normal of $\elDom$ is denoted by $\ubn_e$. 

In the following, all definitions and notation are given for $d=2$, but are equally applicable to $d=3$ if ``edge" is replaced with ``face" in each instance.

Each element has a boundary $\elBound{e}$, consisting of edges $E$.  Define $h_E \defined \text{diam} \lrb{E}$. 

The union of all edges lying in the interior of the domain, rather than on the boundary, will be denoted by $\intBound$. 
Define $\Gamma_{iD} \defined \intBound \cup \DBound$.
By abuse of notation, any symbol denoting a union of edges will also denote the corresponding set of edges.

Use will also be made of the discrete Sobolev space 
\begin{align}
\Hop \defined \lrcb{v \in \Ltd : \ \ v|_{\elDom} \in \Hoe \ \ \forall \ \elDom \in \partition}. \notag
\end{align}
Some kind of weak continuity is nevertheless required between neighbouring elements. For a vector $\ubv$ and a tensor $\btau$, with components in $H^1(\elDom)$ and $H^1(\Omega_f)$ for adjacent elements $\elDom$ and $\Omega_f$ with common edge $\Gamma_{ef}$, the jumps are defined by 
\begin{align}
\lrjmpV{\ubv} &\defined \ubv_e \otimes \ubn_e + \ubv_f \otimes \ubn_f, \ \
\lrjmpV{\btau} \defined \btau_e \ubn_e + \btau_f \ubn_f   \notag \\
\text{and }\hspace{10mm} \lrjmp{\ubv} &\defined \ubv_e \cdot \ubn_e + \ubv_f \cdot \ubn_f  \notag 
\end{align}
and the averages by
\begin{align}
\lravg{\ubv} &\defined \half{\lrb{\ubv_e + \ubv_f}},
\ \
 \lravg{\btau} \defined \half{\lrb{\btau_e + \btau_f}}. \notag 
\end{align}
On edges $E$ such that $E \cap \partial \Omega \neq \emptyset$, the jumps and averages are defined by
\begin{align}
&\lrjmpV{\ubv} = \ubv \otimes \ubn, \ \ \lrjmpV{\btau} = \btau \ubn,  \notag \\
&\lrjmp{\ubv} = \ubv \cdot \ubn, \notag \\
&\lravg{\ubv} = \ubv, \ \ \lravg{\btau} = \btau. \notag 
\end{align}

With $\BBQ_1 \lrb{\Omega}$ the space of polynomials on $\Omega$ with maximum degree one in each variable, define the DG solution space
\begin{align}\label{def_Vh}
V_h = \lrsb{\ubv \in \lrsb{\Ltd}^d: \ \ubv|_{\elDom} \in \lrsb{\BBQ_1 \lrb{\Omega}}^d \ \ \forall \ \elDom \in \partition},
\end{align}
where $V_h \subset \lrsb{\Hop}^d$.

\subsection{Modified interior penalty (IP) formulation}

In the general formulation presented in \cite{Grieshaber2015}, which defines three IP methods, selected edge terms are under-integrated (indicated below by a summation over the Gauss points, $\rm{p}_i$, for $i = 1, ..., ngp$, with $w_i$ referring to the corresponding weighting value). With the context here being the use of bilinear or trilinear elements, the appropriate under-integration employs a single Gauss point (that is, $ngp = 1$). 

With non-negative parameters $k_{\mu}$ and $k_{\lambda}$, and a switch $\theta$ to distinguish between methods, the formulation is defined by the bilinear form 
\begin{align} \label{defUI3}
&a^{\rm{UI}}_h \lrb{\ubu, \ubv} 
= \sumalle \intEl \sigu \colon \epsv dx \notag\\
    &\hspace{10mm}  + \theta \ 2 \mu \sumiDE \int_E  \lrjmpV{\ubu} \colon \lravg{\epsv} ds 
  + \theta \ \lambda \sumiDE \sum_{i=1}^{ngp} \lrsb{\lrjmpV{\ubu} \colon \lravg{\divV \ubv \idten} }|_{{\rm{p}_i}} w_i \notag \\
    & \hspace{10mm}	- 2 \mu \sumiDE \int_E  \lravg{\epsu} \colon \lrjmpV{\ubv} ds 
    - \lambda \sumiDE \sum_{i=1}^{ngp} \lrsb{ \lravg{\divV \ubu \idten} \colon \lrjmpV{\ubv}} |_{{\rm{p}_i}} w_i \notag \\
    & \hspace{10mm} + k_{\mu} \,\mu \sumiDE \frac{1}{h_E} \int_E \lrjmpV{\ubu} \colon \lrjmpV{\ubv} \,ds 
    + k_{\lambda} \,\lambda \sumiDE \frac{1}{h_E} \sum_{i=1}^{ngp} \lrsb{ \lrjmp{\ubu}  \lrjmp{\ubv}} |_{{\rm{p}_i}} w_i,
\end{align}
and linear functional
\begin{align} \label{def_l_UI3}
&l^{\rm{UI}}_h \lrb{\ubv} 
= \sumalle \intEl \ubf \cdot \ubv \,dx    + \sumNE \int_E \ubh \cdot \ubv  \,ds \notag\\
    &\hspace{3mm} 
    + \theta \ 2 \mu \sumDE \int_E  \lrb{\ubg \otimes \ubn}
\colon \epsv  \,ds     
+ \theta \ \lambda \sumDE \sum_{i=1}^{ngp} \lrsb{ \lrb{\ubg \otimes \ubn} \colon \lrb{\divV \ubv \,\idten}}|_{{\rm{p}_i}} w_i  \notag  \\
    & \hspace{3mm}
    		+ k_{\mu} \,\mu \sumDE \frac{1}{h_E} \int_E \lrb{\ubg \otimes \ubn} 
    		\colon \lrb{\ubv \otimes \ubn}  \,ds 
    + k_{\lambda} \,\lambda \sumDE \frac{1}{h_E} \sum_{i=1}^{ngp} \lrsb{\lrb{\ubg \cdot \ubn} 
    		\lrb{\ubv \cdot \ubn}}|_{{\rm{p}_i}} w_i.
\end{align}

The standard IP formulation is recovered by integrating fully all terms (i.e.\ reverting to $ngp = 2$), in which case $\theta = 1$ gives the NIPG method, $\theta = -1$ gives SIPG, and $\theta = 0$ gives IIPG. (The original methods do not in all cases contain the stabilization terms included here.)

We seek $\ubu_h \in V_h$ such that, for all $\ubv \in V_h$, 
\begin{equation}\label{IPeqn}
a^{\rm{UI}}_h \lrb{\ubu_h, \ubv} = l^{\rm{UI}}_h \lrb{\ubv}. 
\end{equation}

As in the original IP methods, both Dirichlet and Neumann boundary conditions, (\ref{bc_eqnD}) and (\ref{bc_eqnN}), are imposed weakly through this formulation.

The formulation of (\ref{IPeqn}) is shown in \cite{Grieshaber2015} to be stable and optimally convergent, and specifically uniformly convergent in the incompressible limit, provided that the domain, applied forces and boundary conditions satisfy the necessary smoothness requirements (as detailed in \cite{Grieshaber2015}).

\section{Computational examples}\label{numerics}
Four model problems are studied, all with analytical solutions available for comparison of results: the first three in two dimensions under plane strain conditions, and the fourth in three dimensions. (Implementation makes use of the deal.ii finite element library \cite{Arndt2017}.)

In each case, two values of Poisson's ratio are considered: $\nu = 0.3$ and $\nu = 0.49995$ representing respectively a compressible and a nearly incompressible material.

\begin{figure}[!ht]
\centering
\subfloat[\textit{df} = 0.0]{\includegraphics[width=.30\columnwidth]{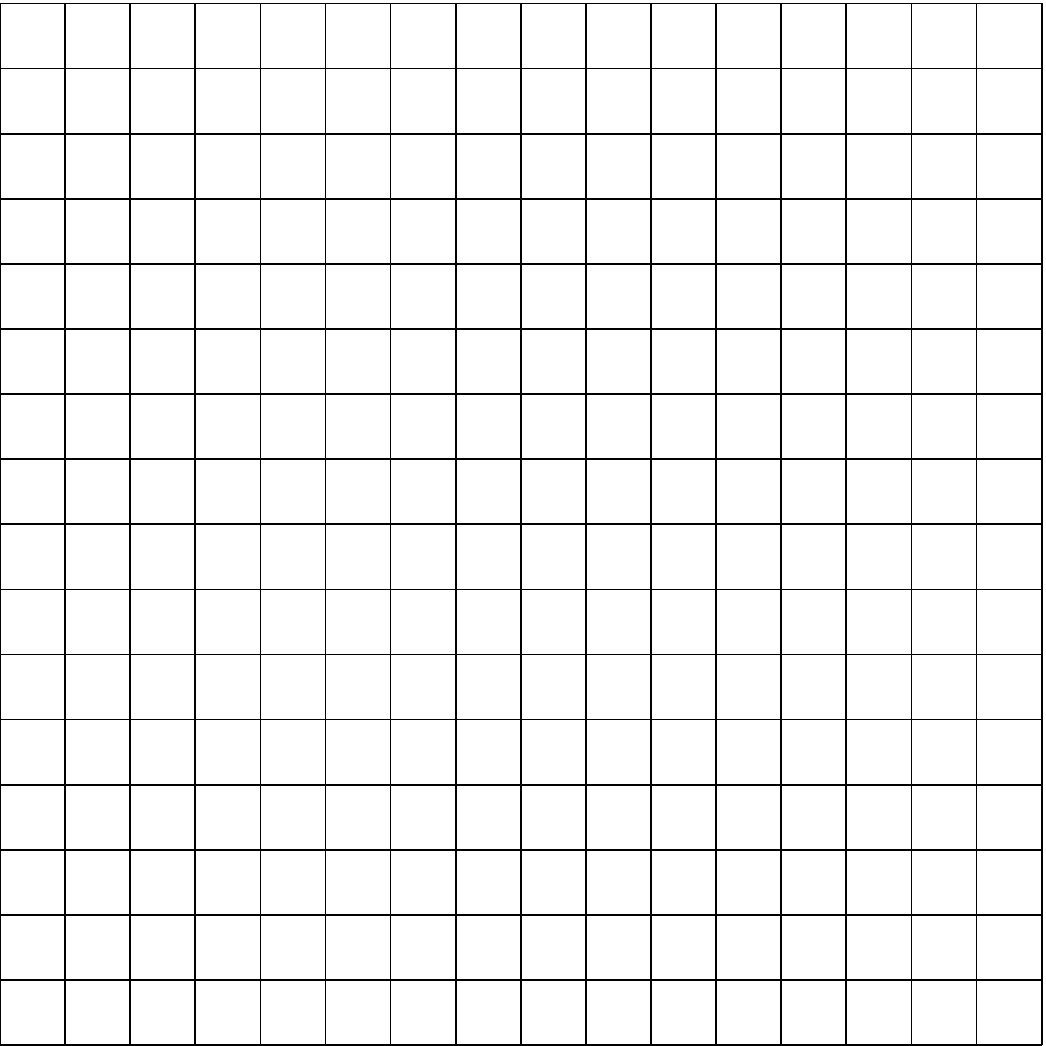}} \hspace{15mm}
\subfloat[\textit{df} = 0.1]{\includegraphics[width=.30\columnwidth]{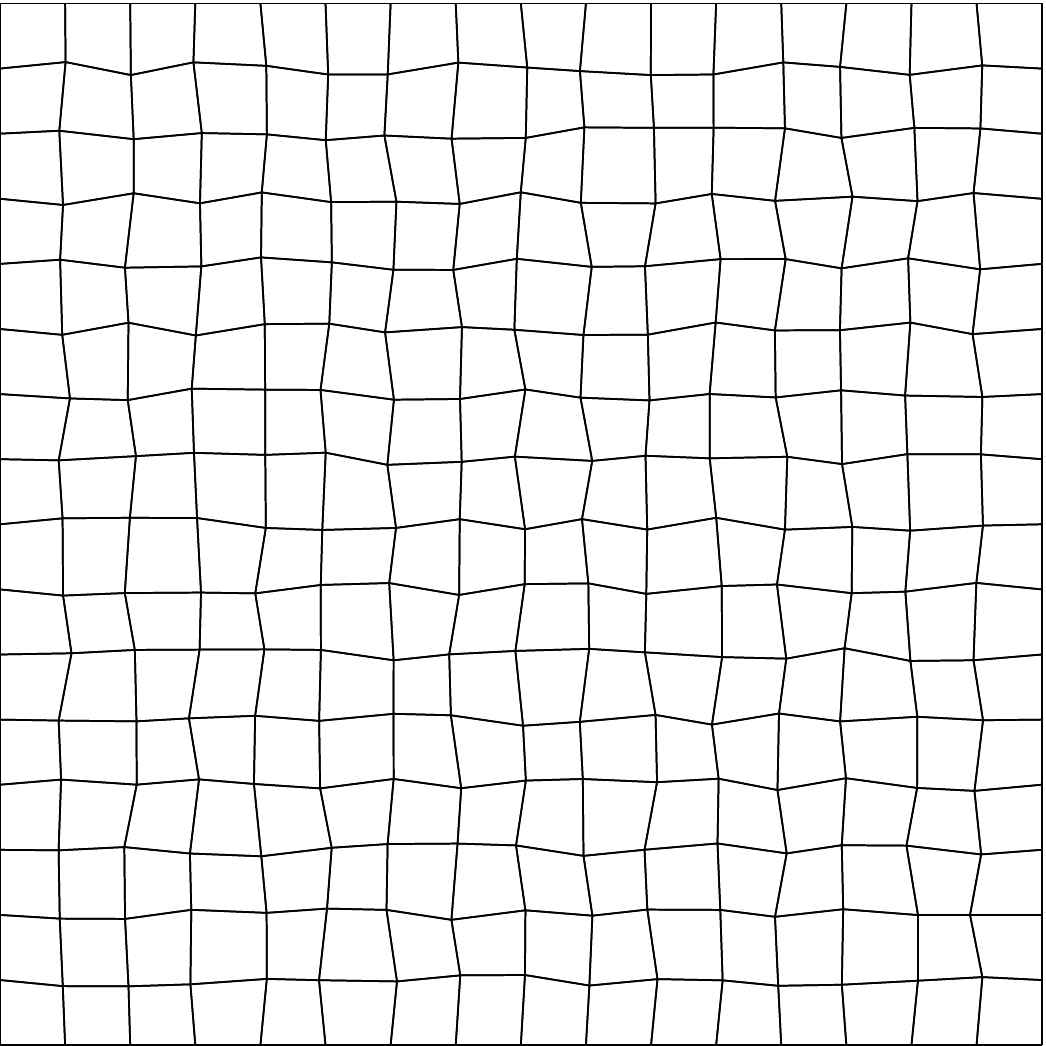}}\\
\subfloat[\textit{df} = 0.3]{\includegraphics[width=.30\columnwidth]{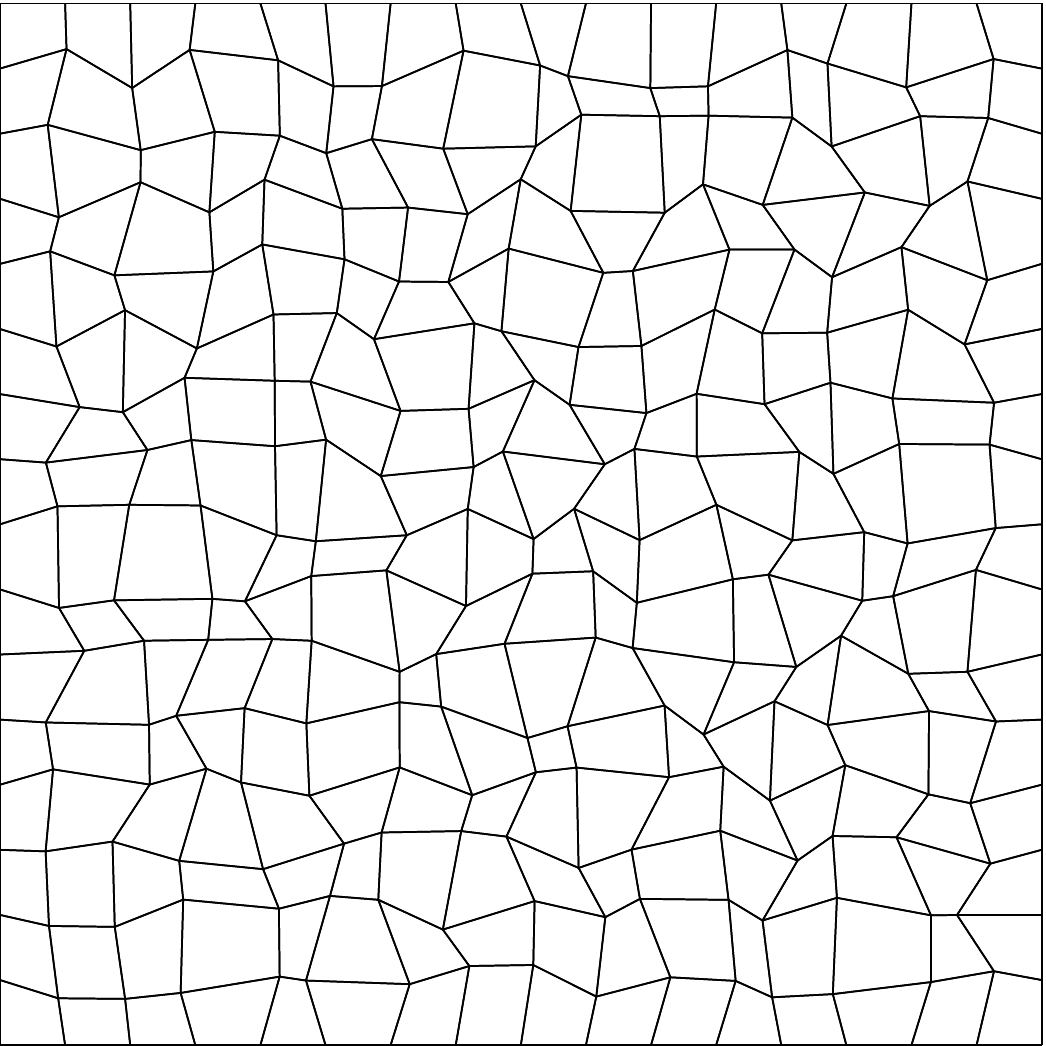}}
\caption{Refinement level 4 for the cantilever beam and square plate problems for different mesh distortion factors}
\label{df_egs}
\end{figure}

In order to ascertain the effects of deviating from the use of rectangular elements with the IP methods, in the first, second and fourth examples, comparative performance of the methods is investigated on meshes with systematically increasing degrees of mesh distortion. At each refinement level, an initial mesh is generated by isotropically refining from a single element. The resulting elements are squares in two dimensions or cubes in three dimensions. Meshes are then modified according to an algorithm (applied within the deal.ii library \cite{Arndt2017}) that distorts each element up to a prescribed distortion factor (\textit{df}), higher values of \textit{df} producing meshes with lower element shape regularity. Meshes of \textit{df} = 0.0 (initial, isotropically refined meshes), \textit{df} = 0.1 and \textit{df} = 0.3 are considered here. (See Figure \ref{df_egs} for an example at a given refinement level.)

In the third example, the L-shaped domain, an unstructured, graded initial mesh is used, with a wide range of element regularity included. This mesh is then refined to investigate convergence behaviour. It should be noted that in general the shape regularity of the elements will be increased by the refinement process in this example.

Unless otherwise indicated, the IP stabilization parameters are set at $k_{\mu} = k_{\lambda} = 10$ for the IIPG and SIPG methods, and $k_{\mu} = 10$, $k_{\lambda} = 0$ for the NIPG method.

Results are given for both the original (with full edge integration) and the new IP methods. It may be assumed henceforth that references to the IP method are to the new method of \cite{Grieshaber2015}, unless explicit mention is made of the original IP methods.

Results obtained using the SG method with $Q_1$ and $Q_2$ elements, and $Q_1$ with SRI, have been included for perspective on the effects of deviating from rectangular elements in both the compressible and near-incompressible cases.

Convergence plots of the $H^1$ error for displacement approximations and the $L^2$ error for post-processed stresses are displayed, where the mesh measure $h$ is the average element diagonal. Contour plots for the displacement approximations use the nodal solution values, while for those of the stress approximations, stress values are calculated at quadrature points and a projection is done over the domain, onto the nodal points, giving continuous stress fields.

Where the exact solution is shown in a contour plot, nodal values are calculated directly from the analytical solution at nodal points, for both displacement and stress fields, except in the case of the L-shaped domain. There, for the stress field, the exact values are calculated from the analytical solution at quadrature points, and a projection, like that for the approximate stresses, is done to obtain nodal values. The reason and implications will be discussed where the results of that example are given.

\begin{table}
\begin{center}
\begin{tabular}{|c|c|c|c|c|} \hline
$n$ & no.\ els & \multicolumn{3}{|c|}{no.\ dofs} \\ \hline 
& & SG $Q_1$ & SG $Q_2$ & IP \\ \hline
1 & 4 & 18 & 50 & 32\\ \hline
2 & 16 & 50 & 162  & 128\\ \hline
3 & 64 & 162 & 578 & 512\\ \hline
4 & 256 & 578 & 2178  & 2048\\ \hline
5 & 1024 & 2178 & 8450 & 8192\\ \hline
6 & 4096 & 8450 & 33282 & 32768\\ \hline
7 & 16384 & 33282 & 132098 & 131072\\ \hline
8 & 65536 & 132098 & 526338 & 524288\\ \hline
9 & 262144 & 526338 & 2101250 & 2097152\\ \hline
\end{tabular}
\end{center}
\caption{Mesh details for standard Galerkin (SG) and interior penalty (IP) methods for the cantilever beam and square plate examples: mesh $n$ has $2^n$ elements per side; number of elements and number of degrees of freedom are shown. (Note: details for SG $Q_2$ and for $n = 9$ apply to the square plate example only.)}\label{mesh_info_CB_SP}\end{table}

\subsection{Cantilever beam}

\begin{figure}[!ht]
\centering
\includegraphics[width=.35\columnwidth]{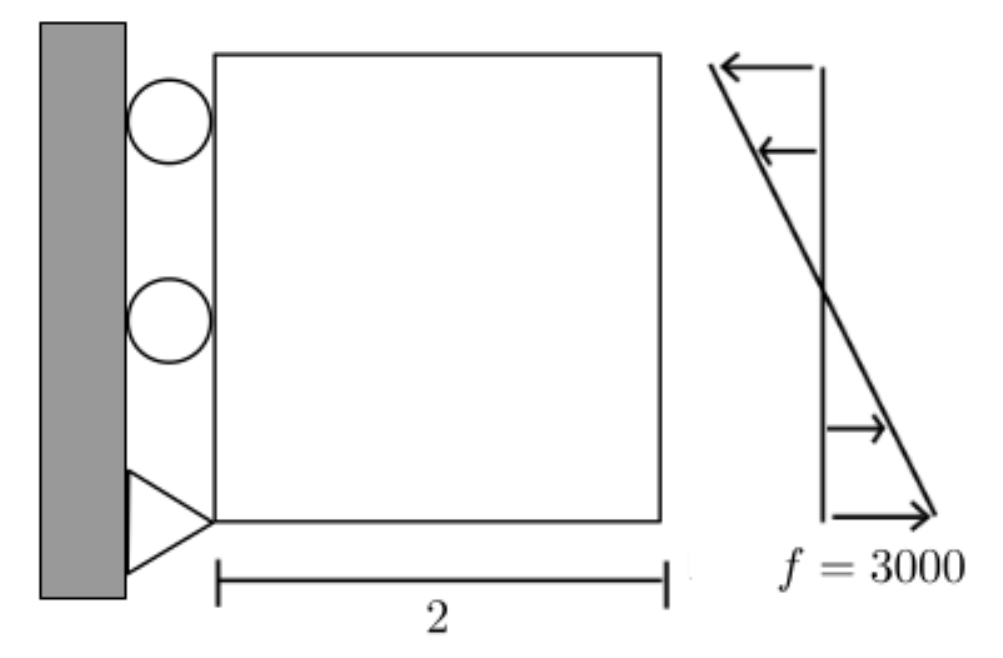}
\caption{Cantilever beam with boundary conditions}
\label{beam_diagram}
\end{figure}

A square beam in two dimensions, with $E = 1 500 000$, is subjected to a linearly varying force on the free end, with the maximum value $f = 3000$, as illustrated in  Figure~\ref{beam_diagram}. 
The analytical solution is given in \cite{Djoko2006}. 

Eight levels of mesh refinement are used, and Table \ref{mesh_info_CB_SP} details the number of elements and degrees of freedom (dofs) for each method at each level.

Displacement results (scaled) of the original and new IP methods with square elements (i.e.\ with \textit{df} = 0.0) have been presented in \cite{Grieshaber2015} and are repeated here for comparison.

\subsubsection{Displacement approximation}

\begin{figure}[!ht]
\centering
\subfloat[SG $Q_1$]{\includegraphics[trim={3.5cm 9cm 4cm 9cm},clip,width=.49\columnwidth]{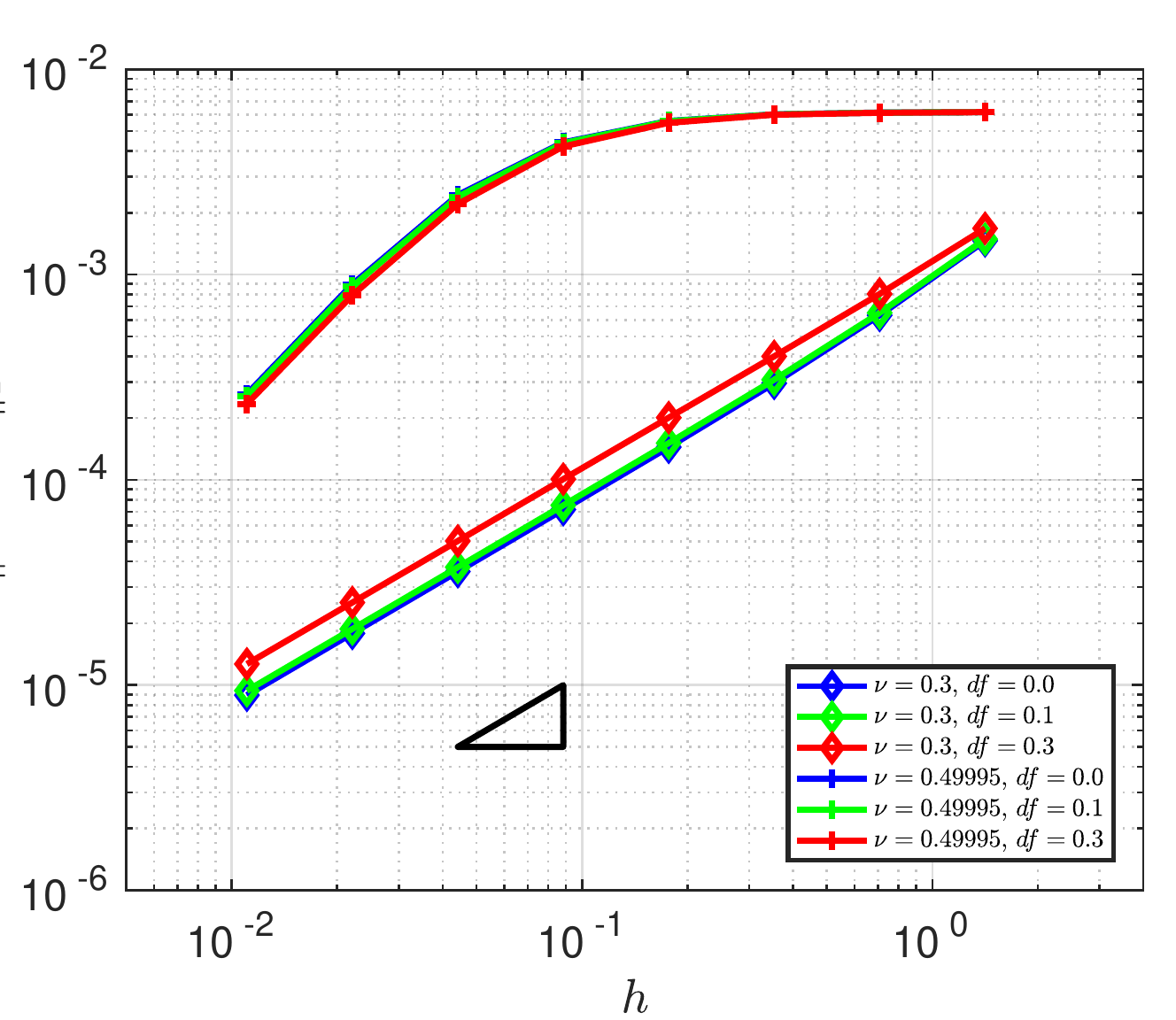}\label{CB_disp_error_SG}}
\subfloat[Original IP, $\nu = 0.49995$]{\includegraphics[trim={3.5cm 9cm 4cm 9cm},clip,width=.49\columnwidth]{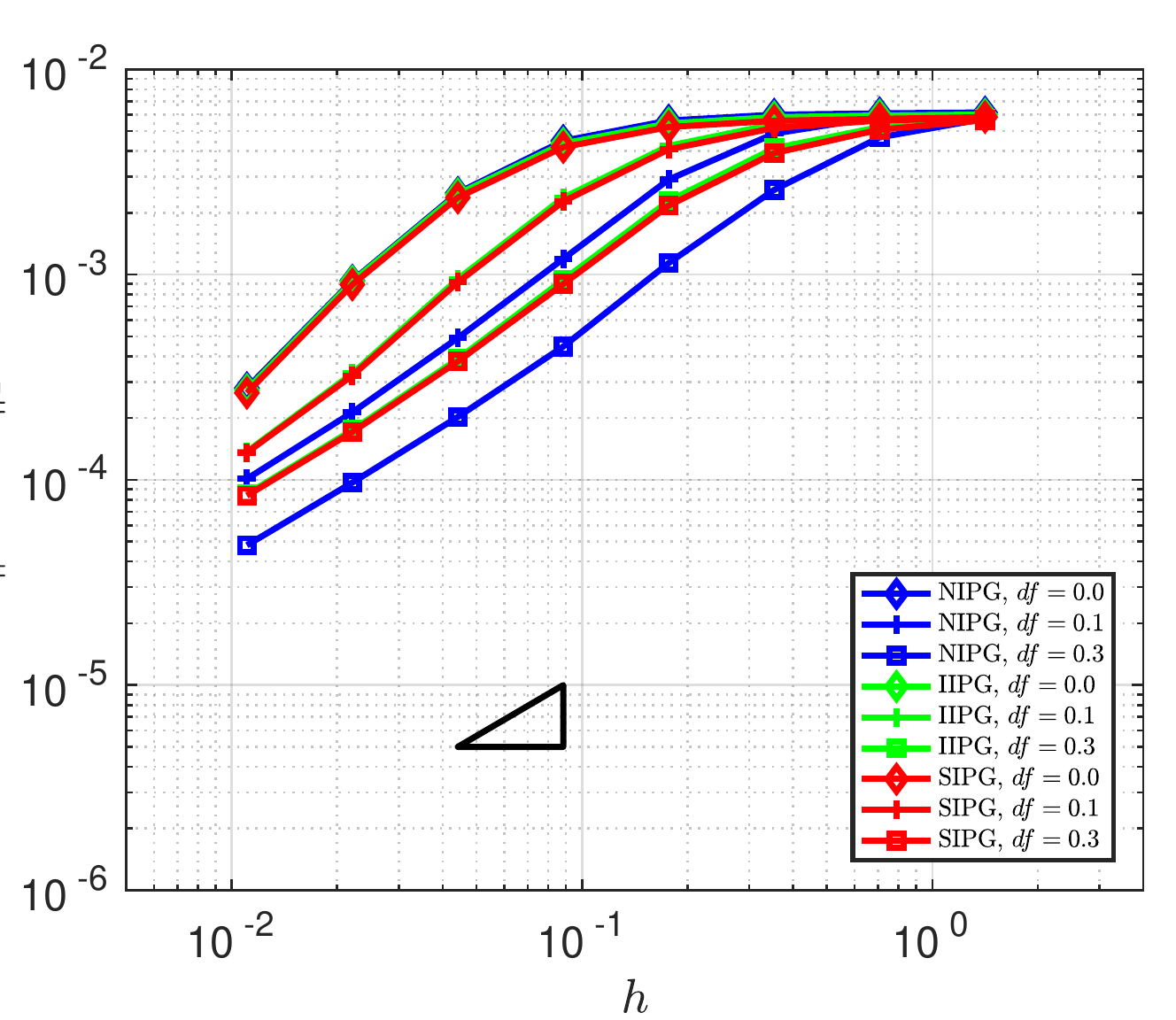}\label{CB_disp_error_oldIPnu2}}\\
\subfloat[New IP, $\nu = 0.49995$]{\includegraphics[trim={3.5cm 9cm 4cm 9cm},clip,width=.49\columnwidth]{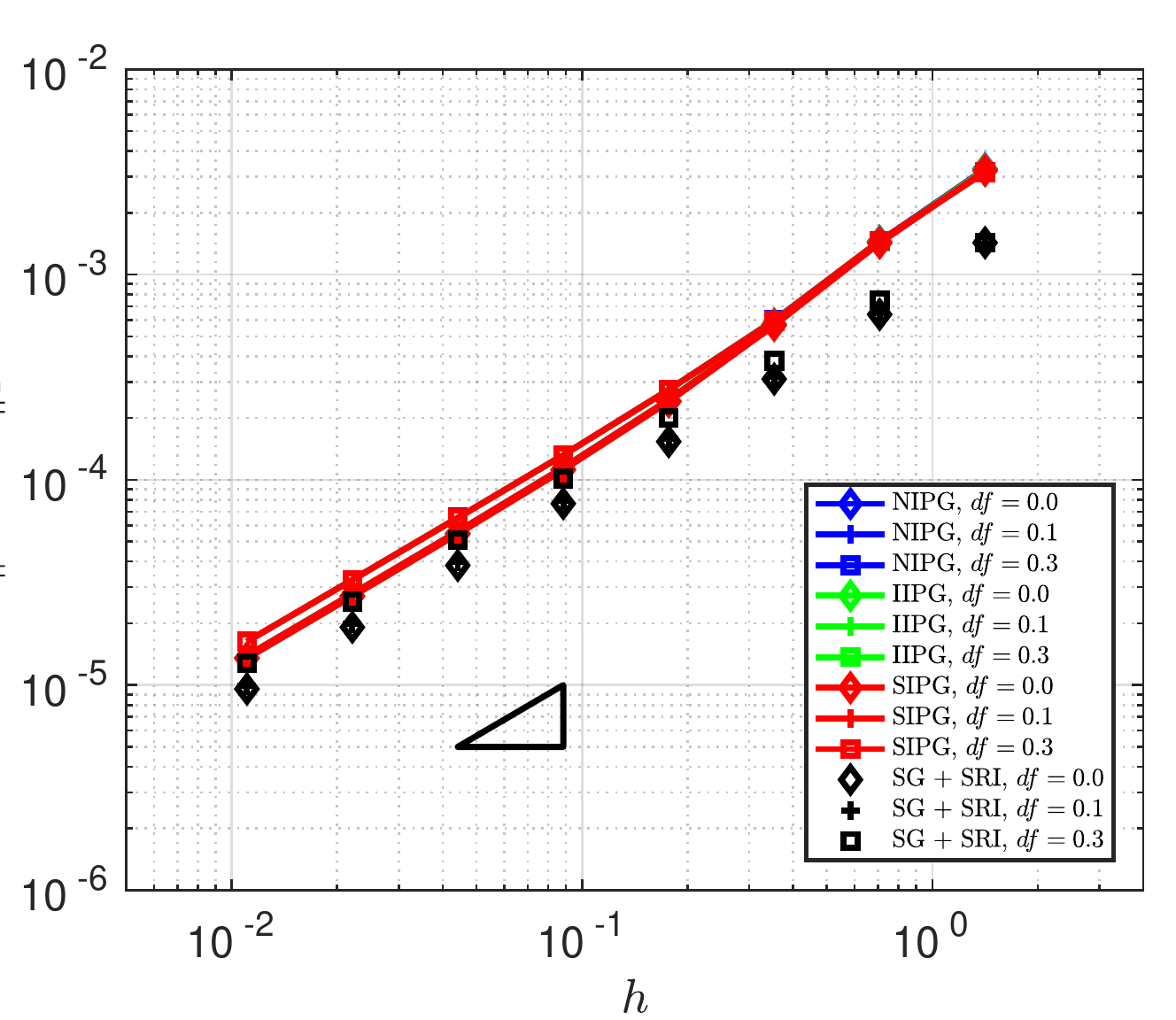}\label{CB_disp_error_newIPnu2}}
\subfloat[New NIPG]{\includegraphics[trim={3.5cm 9cm 4cm 9cm},clip,width=.49\columnwidth]{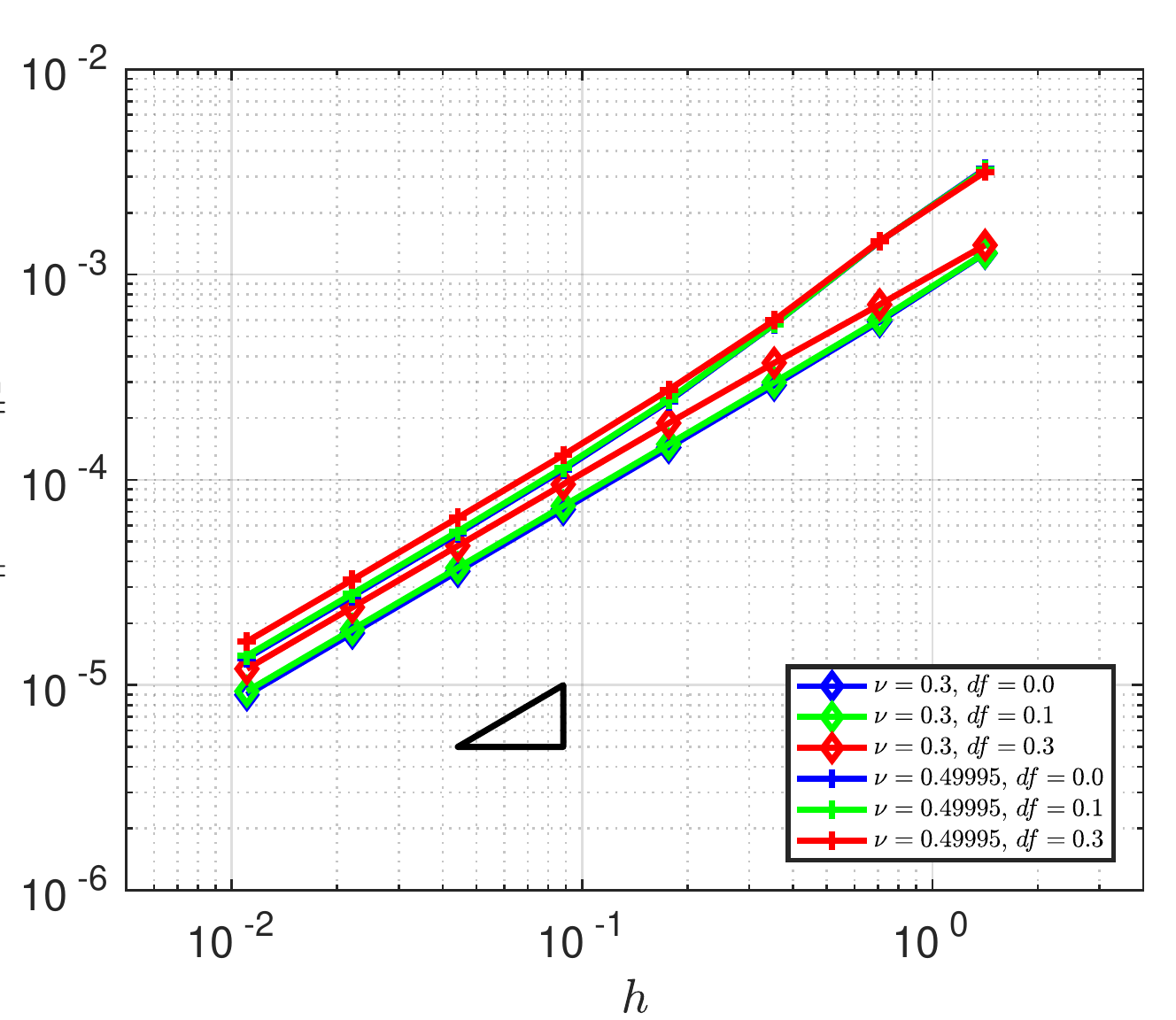}\label{CB_disp_error_NIPG}}
\caption{Cantilever beam: Displacement $H^1$ error convergence. The hypotenuse of the triangle has a slope of 1 in each case.}
\label{CB_disp_error}
\end{figure}

\begin{figure}[!ht]
\centering
\subfloat[Exact solution]{\includegraphics[width=.20\columnwidth]{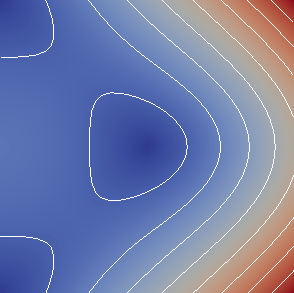}\label{CB_disp_exact}}\hspace{10mm}
\subfloat[NIPG, \textit{df} = 0.0, mesh 5]{\includegraphics[width=.20\columnwidth]{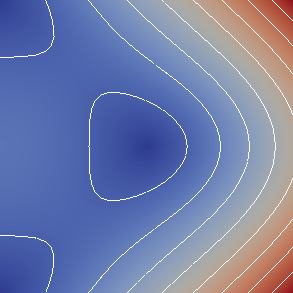}\label{CB_disp_NIPG_00}}\hspace{10mm}
\subfloat[NIPG, \textit{df} = 0.3, mesh 5]{\includegraphics[width=.26\columnwidth]{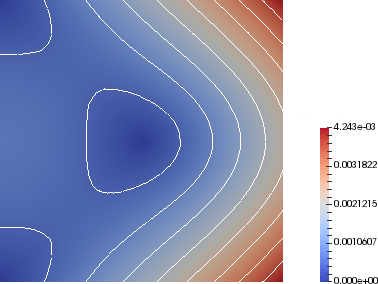}\label{CB_disp_NIPG_03}}
\caption{Cantilever beam: Displacement magnitude, $\nu = 0.49995$}\label{CB_disp_magn}
\end{figure}

Low-order ($Q_1$) elements with the SG method (Figure  \ref{CB_disp_error_SG}) converge optimally for $\nu = 0.3$, as expected, with a slight decrease in accuracy as \textit{df} increases; but show poor convergence, indicating locking, with $\nu = 0.49995$ for all values of \textit{df}, except for high refinement levels. Figure  \ref{CB_disp_error_oldIPnu2} shows the results obtained using the original IP methods for $\nu = 0.49995$, with varying values of \textit{df}. Poor convergence is exhibited for coarser meshes in all cases, although optimal convergence is attained at lower refinements as element regularity decreases. In contrast, the new IP methods (Figure  \ref{CB_disp_error_newIPnu2}) display optimal convergence when $\nu = 0.49995$, irrespective of element regularity, with a small increase in error magnitude as \textit{df} increases. The error of the new IP methods is significantly smaller than that of the original IP methods. The SG method with SRI, included here for comparison, also shows optimal convergence, with slightly better accuracy than the IP methods. Finally, the uniform optimal convergence of the new IP methods with respect to $\lambda$, independent of \textit{df}, is illustrated in Figure  \ref{CB_disp_error_NIPG} (shown for the NIPG method).

Figure \ref{CB_disp_magn} depicts the displacement magnitude obtained using the NIPG method in the near-incompressible case, for refinement level 5, comparing it to the exact solution: the excellent accuracy achieved using square elements (\textit{df} = 0.0) is maintained with the significantly lower element regularity of \textit{df} = 0.3. The performance of the other IP methods is the same.

\subsubsection{Post-processed stress}


\begin{figure}[!ht]
\centering
\subfloat[SG $Q_1$]{\includegraphics[trim={3.5cm 9cm 4cm 9cm},clip,width=.49\columnwidth]{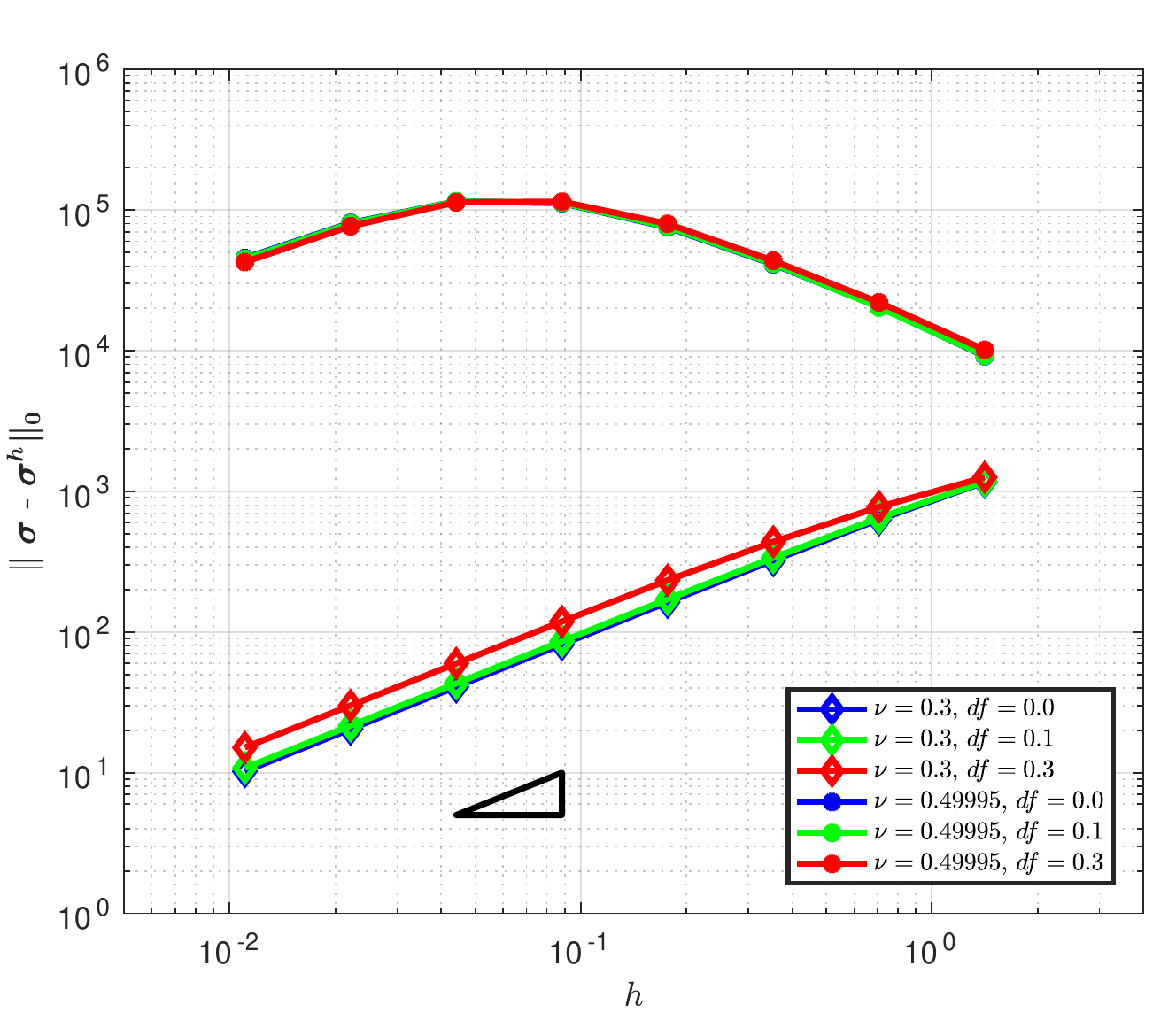}\label{CB_stress_error_SG}}
\subfloat[Original IP, $\nu =0.49995$]{\includegraphics[trim={3.5cm 9cm 4cm 9cm},clip,width=.49\columnwidth]{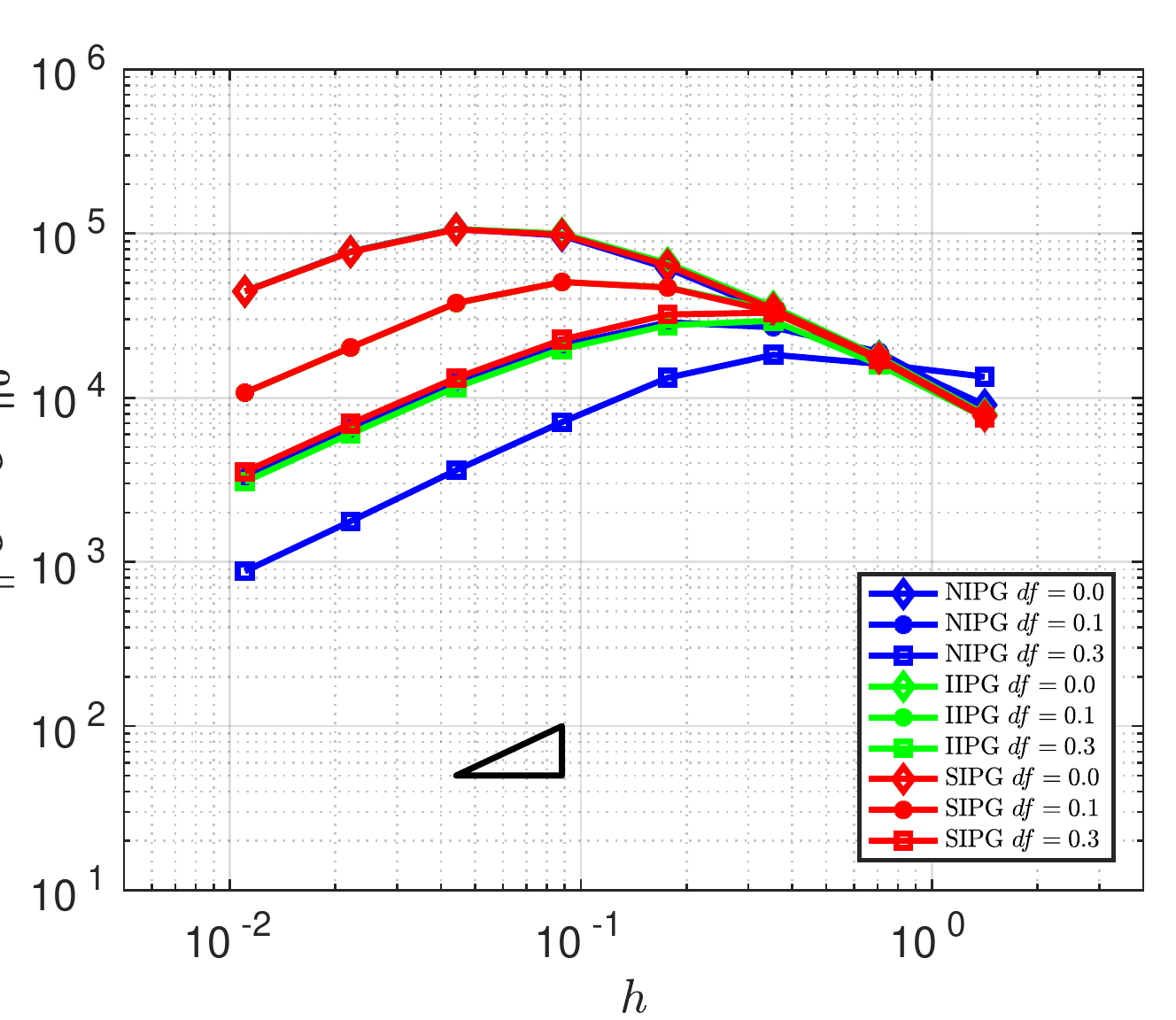}\label{CB_stress_error_oldIPnu2}}\\
\subfloat[New IP, $\nu =0.49995$]{\includegraphics[trim={3.5cm 9cm 4cm 9cm},clip,width=.49\columnwidth]{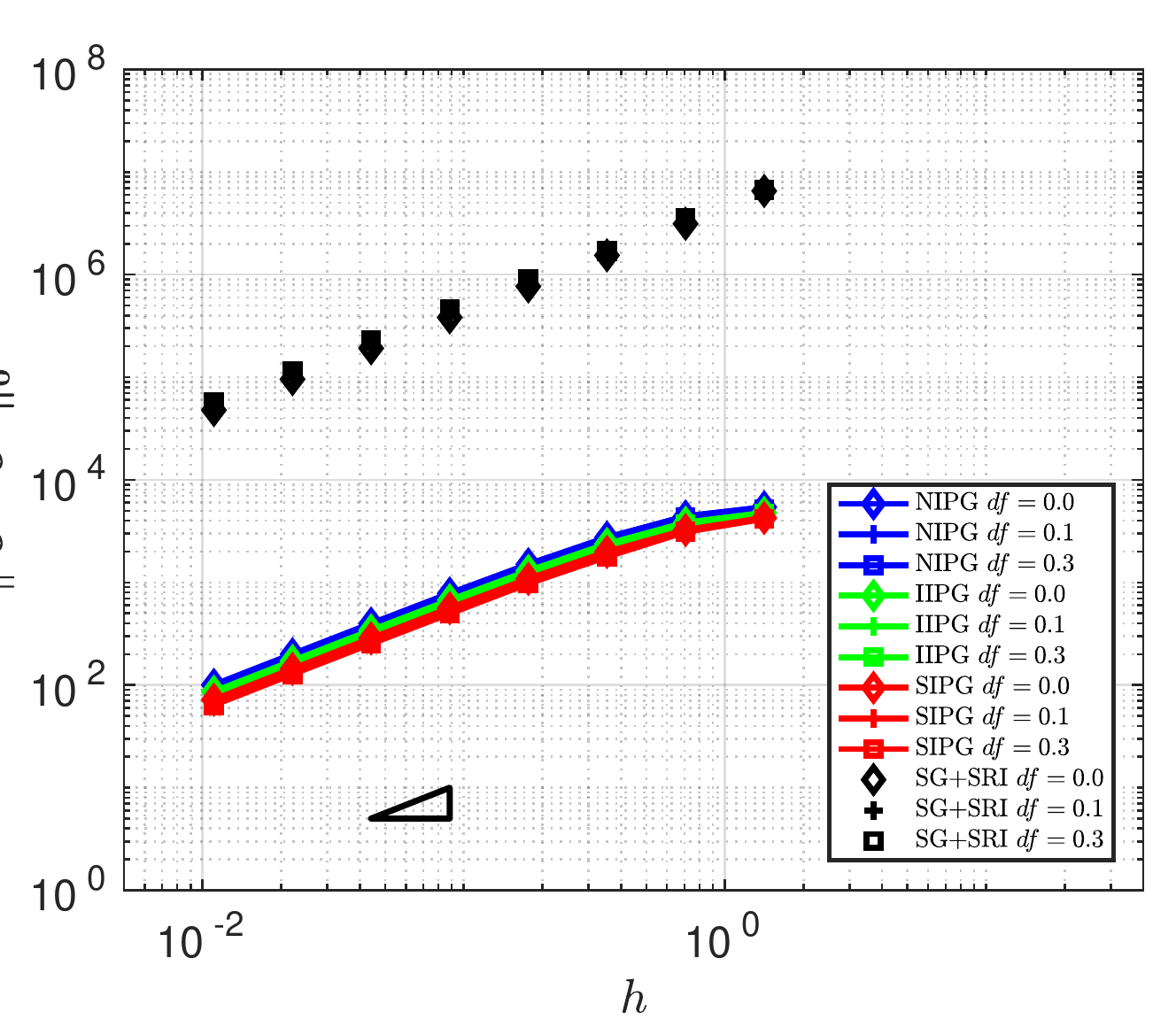}\label{CB_stress_error_newIPnu2}}
\subfloat[New NIPG]{\includegraphics[trim={3.5cm 9cm 4cm 9cm},clip,width=.49\columnwidth]{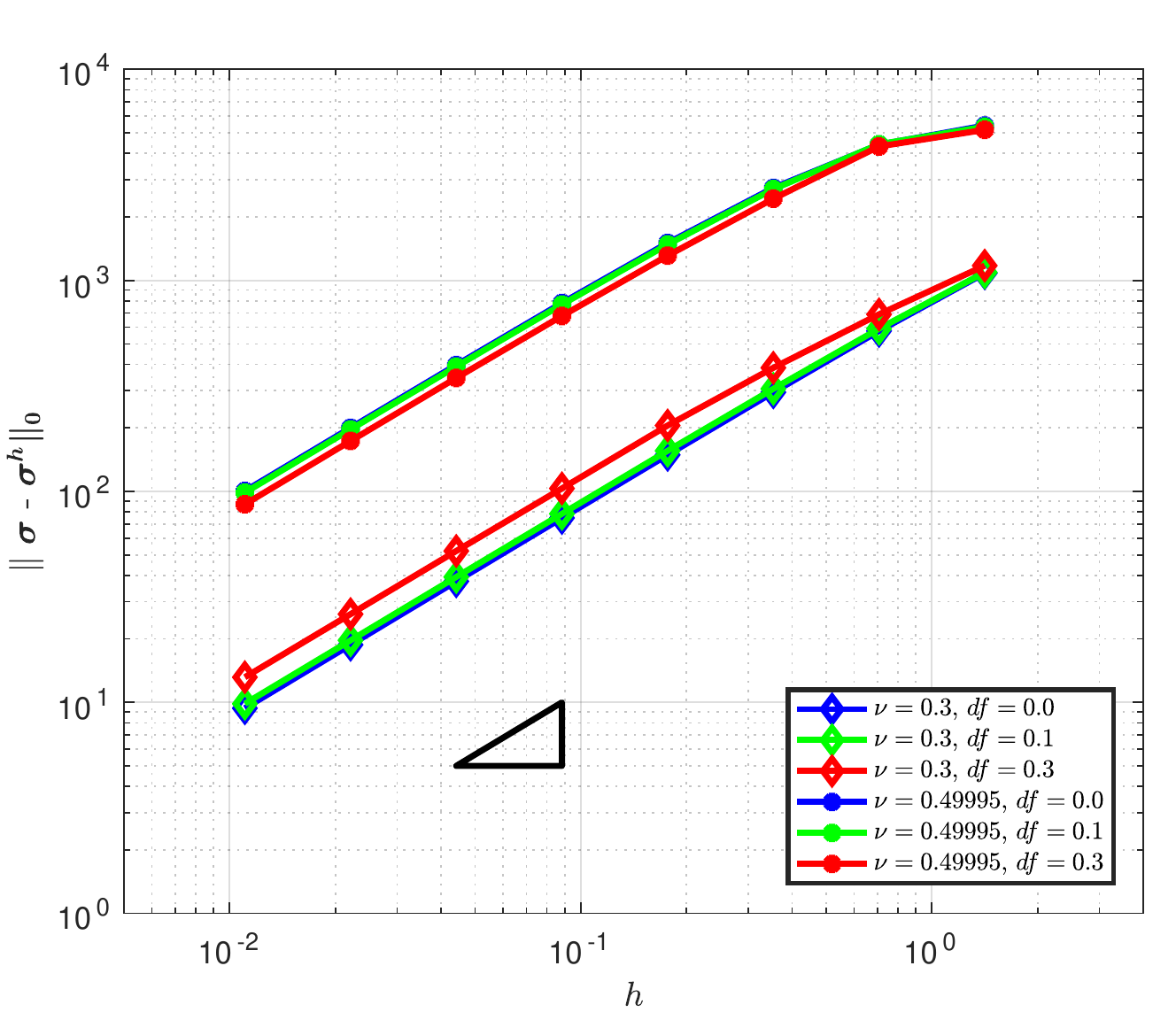}\label{CB_stress_error_NIPG}}
\caption{Cantilever beam: Stress $L^2$ error convergence. The hypotenuse of the triangle has a slope of 1 in each case.}
\label{CB_stress_error}
\end{figure}


\begin{figure}[!ht]
\centering
\subfloat[Exact solution]{\includegraphics[width=.20\columnwidth]{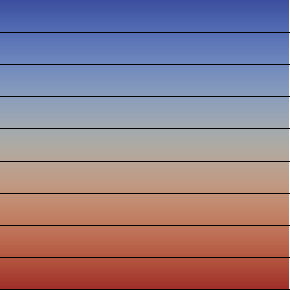}\label{CB_stress_xx_exact}} \hspace{10mm}
\subfloat[NIPG, \textit{df} = 0.0, mesh 5]{\includegraphics[width=.20\columnwidth]{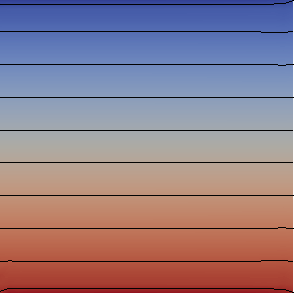}\label{CB_stress_xx_NIPG_00_m5}} \\
\subfloat[NIPG, \textit{df} = 0.1, mesh 5]{\includegraphics[width=.20\columnwidth]{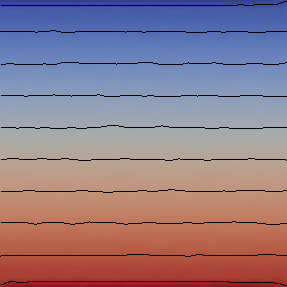}\label{CB_stress_xx_NIPG_01_m5}} \hspace{10mm}
\subfloat[NIPG, \textit{df} = 0.1, mesh 8]{\includegraphics[width=.20\columnwidth]{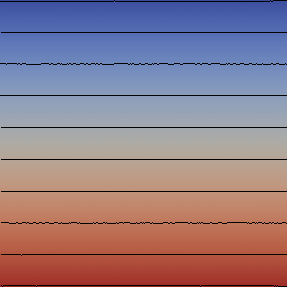}\label{CB_stress_xx_NIPG_01_m8}} \\
\subfloat[NIPG, \textit{df} = 0.3, mesh 5]{\includegraphics[width=.20\columnwidth]{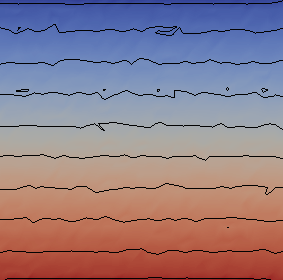}\label{CB_stress_xx_NIPG_03_m5}}\hspace{10mm} 
\subfloat[NIPG, \textit{df} = 0.3, mesh 8]{\includegraphics[width=.20\columnwidth]{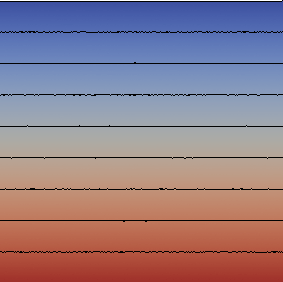}\label{CB_stress_xx_NIPG_03_m8}}\\\hspace{10mm}
\subfloat[NIPG, \textit{df} = 0.3, mesh 5, without contour lines]{\includegraphics[width=.2\columnwidth]{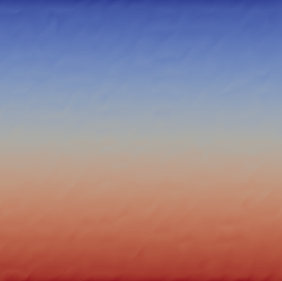}\label{CB_stress_xx_NIPG_03_m5_plain}} \hspace{10mm}
\subfloat[NIPG, \textit{df} = 0.3, mesh 8, without contour lines]{\includegraphics[width=.28\columnwidth]{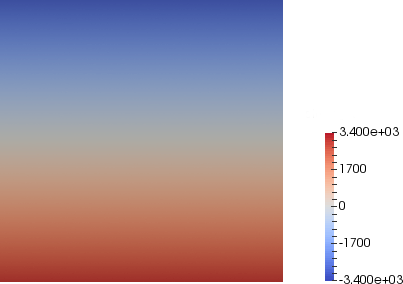}\label{CB_stress_xx_NIPG_03_m8_plain}}\caption{Cantilever beam: $\sigma_{xx}$, $\nu = 0.49995$}\label{CB_stress_xx}
\end{figure}

\begin{figure}[!ht]
\centering
\subfloat[$\sigma_{xy}$, mesh 5]{\includegraphics[width=.20\columnwidth]{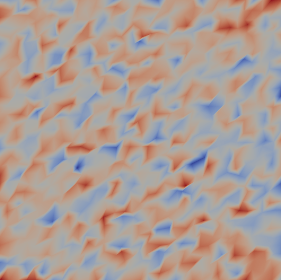}\label{CB_stress_xy_NIPG_03_m5}}\hspace{10mm}
\subfloat[$\sigma_{xy}$, mesh 7]{\includegraphics[width=.275\columnwidth]{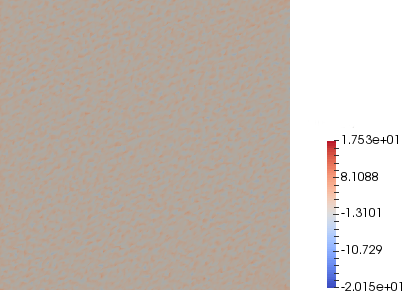}\label{CB_stress_xy_NIPG_03_m8}} \\
\subfloat[$\sigma_{yy}$, mesh 5]{\includegraphics[width=.20\columnwidth]{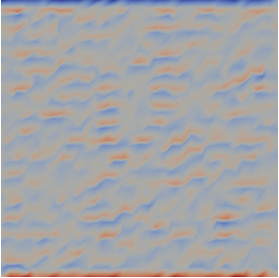}\label{CB_stress_yy_NIPG_03_m5}}\hspace{10mm}
\subfloat[$\sigma_{yy}$, mesh 8]{\includegraphics[width=.275\columnwidth]{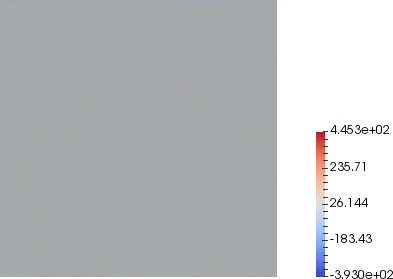}\label{CB_stress_yy_NIPG_03_m8}}
\caption{Cantilever beam: Stress, NIPG, \textit{df} = 0.3, $\nu = 0.49995$}\label{CB_stress_xy_and_yy_NIPG}
\end{figure}

For the SG method with $Q_1$ elements (Figure \ref{CB_stress_error_SG}), for $\nu = 0.3$ there is first-order convergence of the postprocessed stresses, with a slight increase of error magnitude as \textit{df} increases, while for $\nu = 0.49995$ the convergence rate is poor. For $\nu = 0.49995$, the original IP methods display poor stress convergence rates for coarse meshes; these rates improve with refinement, more quickly for higher \textit{df} in a given variant (Figure \ref{CB_stress_error_oldIPnu2}). In contrast, as for the displacement approximations, the new IP methods produce first-order convergence results uniformly with respect to \textit{df} for $\nu = 0.49995$ (Figure \ref{CB_stress_error_newIPnu2}), with significantly lower error magnitudes that those of the original IP method, and uniformly with respect to the compressibility parameter (NIPG results shown in Figure \ref{CB_stress_error_NIPG}). The errors of the SG method with SRI are shown in Figure \ref{CB_stress_error_newIPnu2} for comparison, and while the convergence rate is first-order, the errors are greater than those of the IP methods, unlike for the displacement approximations.

The contour plots in Figure \ref{CB_stress_xx} show, with the exact solution for comparison, the slight deterioration in accuracy of the NIPG post-processed stress field for the component $\sigma_{xx}$, for the near-incompressible case, as element regularity decreases. This is vastly improved by mesh refinement. Figures \ref{CB_stress_xx_NIPG_03_m5_plain} and \ref{CB_stress_xx_NIPG_03_m8_plain}, when compared to Figures \ref{CB_stress_xx_NIPG_03_m5} and \ref{CB_stress_xx_NIPG_03_m8}, illustrate how the contour lines highlight the discrepancies, while the overall stress field is fairly smooth, even for \textit{df} = 0.3 at mesh 5.

The exact stress solutions of the other two (in-plane) components are both 0. When low-regularity elements are used (\textit{df} = 0.3, here) at refinement level 5 the NIPG method produced a mottled $\sigma_{xy}$ field when $\nu = 0.49995$ (Figure \ref{CB_stress_xy_NIPG_03_m5}); however, the values are two orders of magnitude lower than those of the $\sigma_{xx}$ field, and shown on the same scale the $\sigma_{xy}$ field appears as 0. Even at the lower scale, refinement improves the field significantly (Figure \ref{CB_stress_xy_NIPG_03_m8}). Similarly, at refinement level 5, the NIPG $\sigma_{yy}$ field is mottled for \textit{df} = 0.3 when $\nu = 0.49995$ (Figure \ref{CB_stress_yy_NIPG_03_m5}), here at one order of magnitude below the values of the $\sigma_{xx}$ field, but appears as 0 when the same scale is used, and at the original scale is improved by refinement (Figure \ref{CB_stress_yy_NIPG_03_m8}).

The results of the other two IP methods are again similar to those of NIPG.

\subsection{Square plate}

The linear elastic unit square plate $[0,1]^2$ described in \cite{Brenner1993}, with $\mu = 1$, is fixed on all its edges and subjected to an internal body force $\ubf$, with components

\begin{align}
f_x &=  0.04 \,\pi^2 \lrsb{ 4 \sin 2 \pi y \lrb{-1 + 2 \cos  2 \pi x} - \cos \pi \lrb{x+y} + \frac{2}{1 + \lambda} \sin \pi x \sin \pi y }, \notag\\
f_y &=  0.04 \,\pi^2 \lrsb{ 4 \sin 2 \pi x \lrb{1 - 2 \cos  2 \pi y} - \cos \pi \lrb{x+y} + \frac{2}{1 + \lambda} \sin \pi x \sin \pi y }. \notag
\end{align}

Table \ref{mesh_info_CB_SP} details the number of elements and dofs for each method at the various refinement levels.

Displacement results of the original and new IP methods with square elements (i.e.\ with \textit{df} = 0.0) have been presented in \cite{Grieshaber2015} and are repeated here for comparison.

\subsubsection{Displacement approximation}


\begin{figure}[!ht]
\centering
\subfloat[SG $Q_1$]{\includegraphics[trim={3.5cm 9cm 4cm 9cm},clip,width=.49\columnwidth]{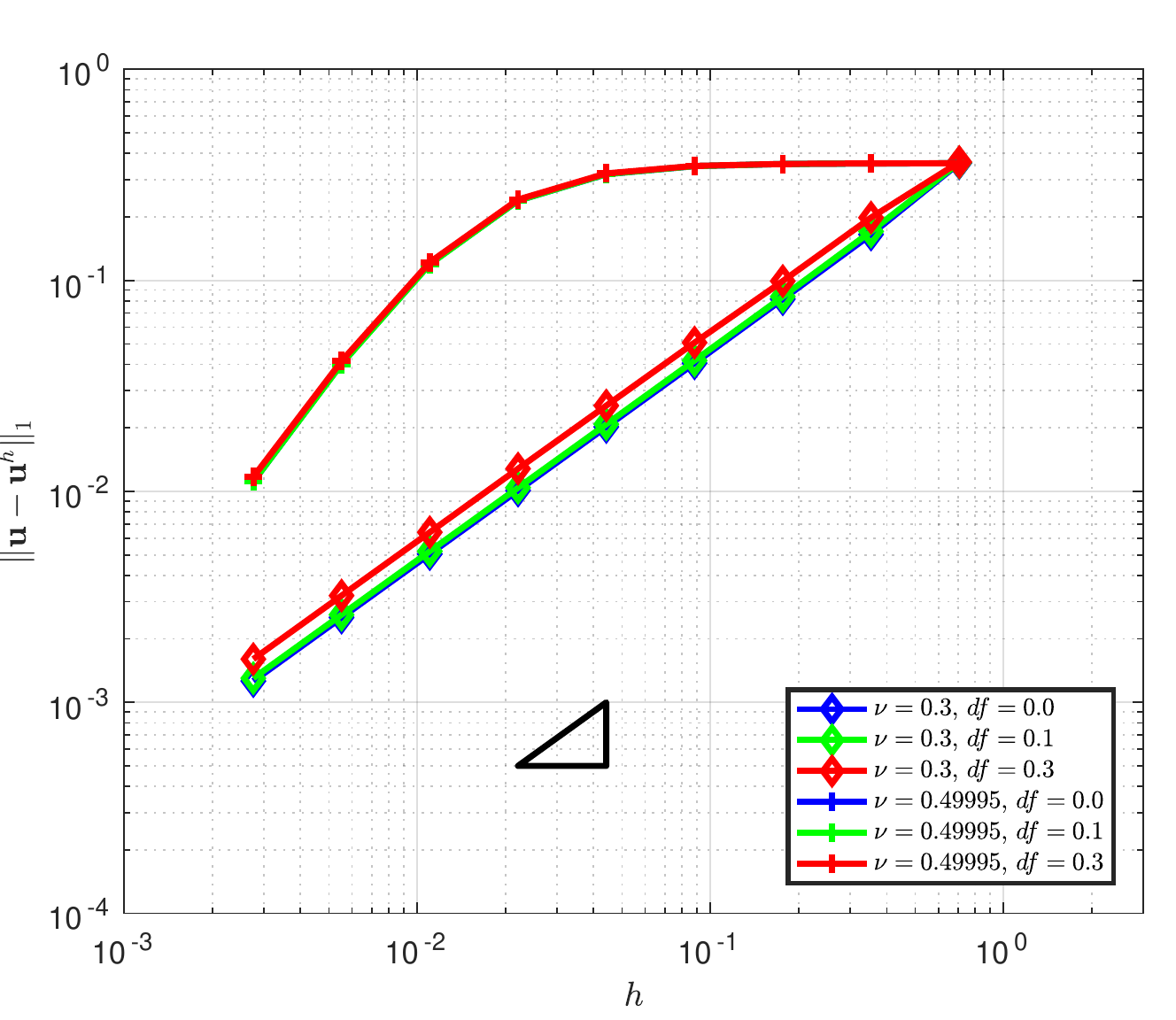}\label{SP_disp_error_SG1}}
\subfloat[SG $Q_1$ with SRI]{\includegraphics[trim={3.5cm 9cm 4cm 9cm},clip,width=.49\columnwidth]{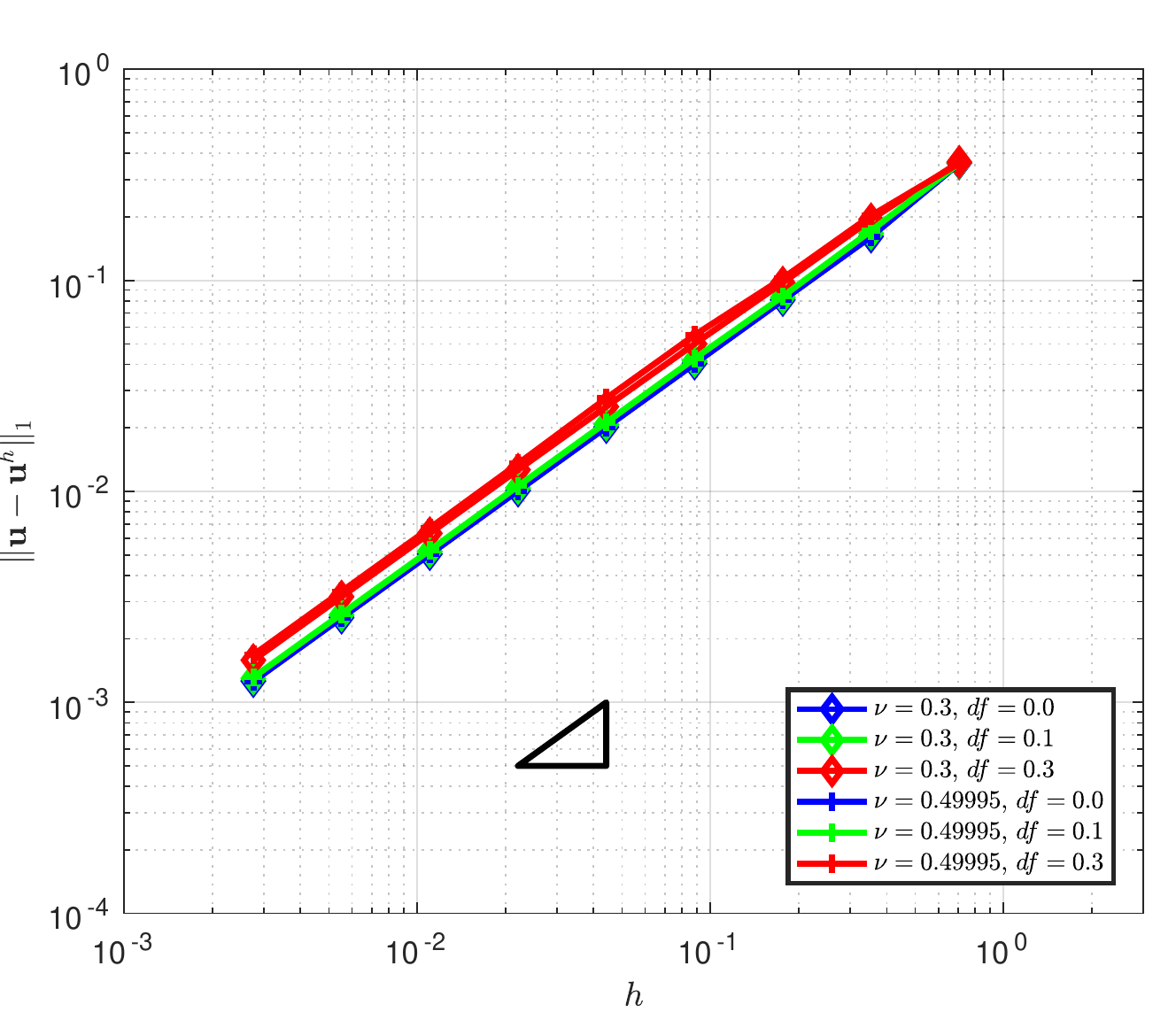}\label{SP_disp_error_SG_sri}}\\
\subfloat[SG $Q_2$]{\includegraphics[trim={3.5cm 9cm 4cm 9cm},clip,width=.49\columnwidth]{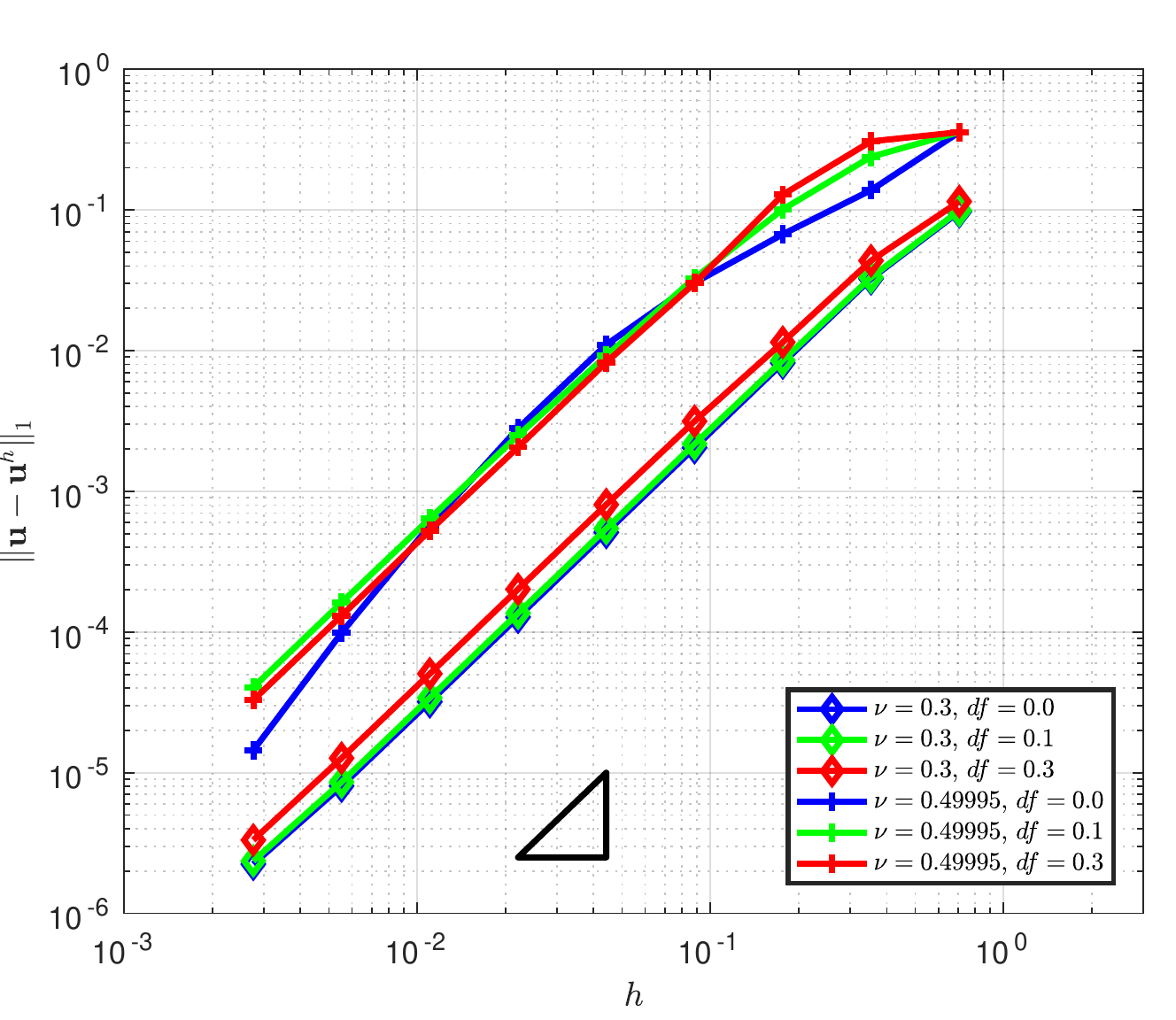}\label{SP_disp_error_SG2}}
\caption{Square plate: SG displacement $H^1$ error convergence. The hypotenuse of the triangle in each of (a) and (b) has a slope of 1, and in (c) has a slope of 2.}\label{SP_disp_error_SG}
\end{figure}

\begin{figure}[!ht]
\centering
\subfloat[Original IP, $\nu = 0.49995$]{\includegraphics[trim={3.5cm 9cm 4cm 9cm},clip,width=.49\columnwidth]{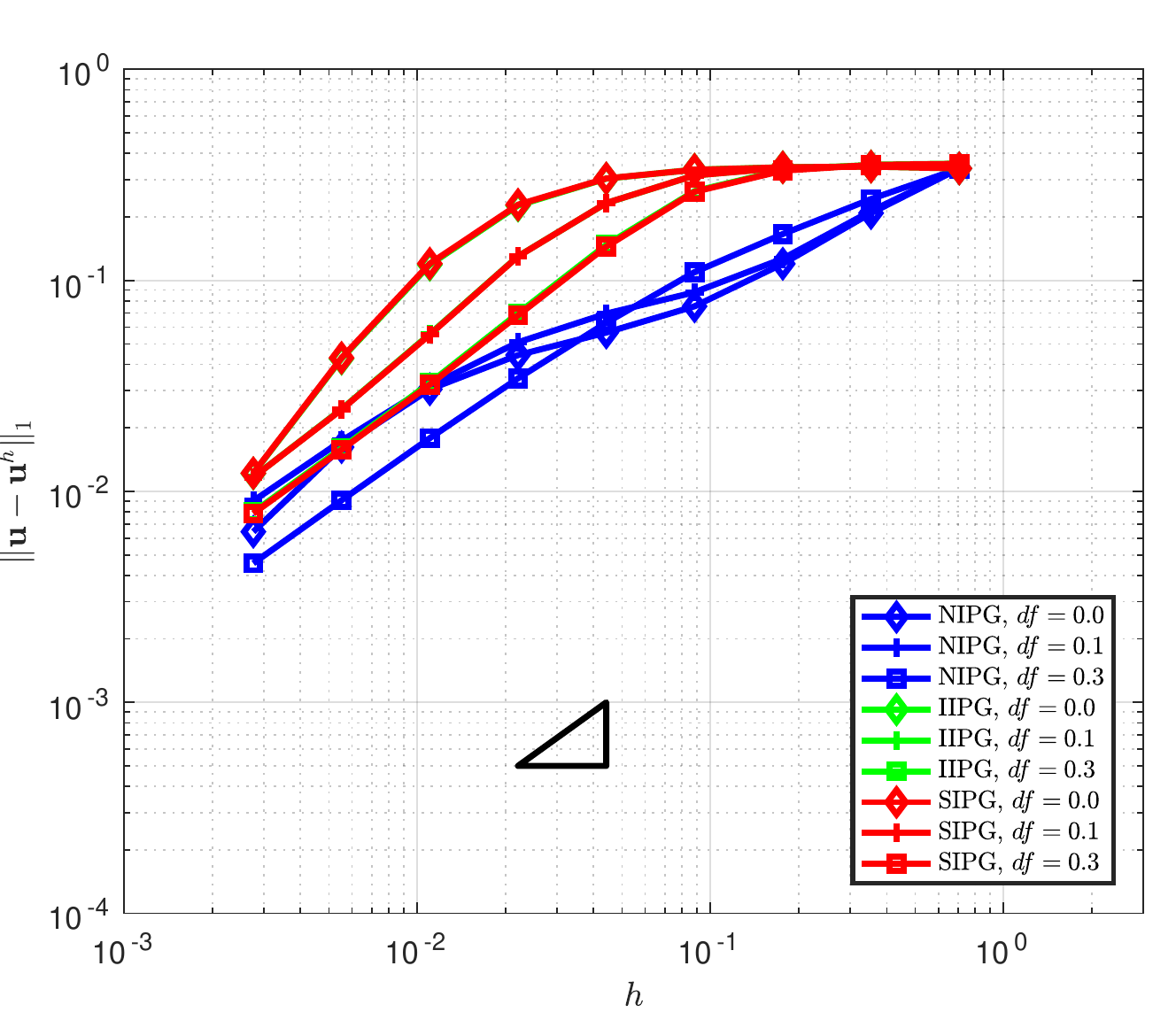}\label{SP_disp_error_oldIPnu2}}
\subfloat[New IP, $\nu = 0.49995$]{\includegraphics[trim={3.5cm 9cm 4cm 9cm},clip,width=.49\columnwidth]{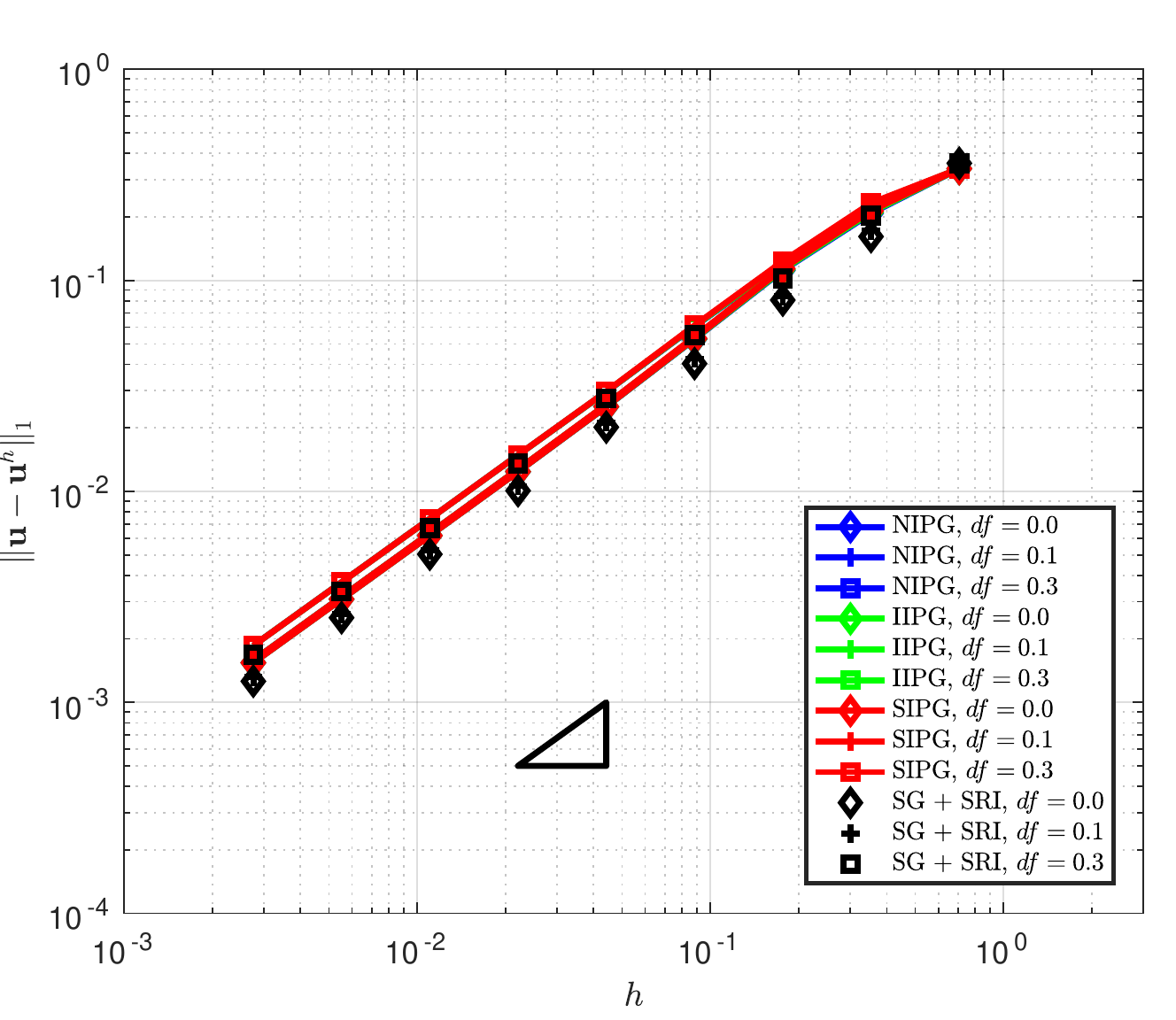}\label{SP_disp_error_newIPnu2}}\\
\subfloat[New NIPG]{\includegraphics[trim={3.5cm 9cm 4cm 9cm},clip,width=.49\columnwidth]{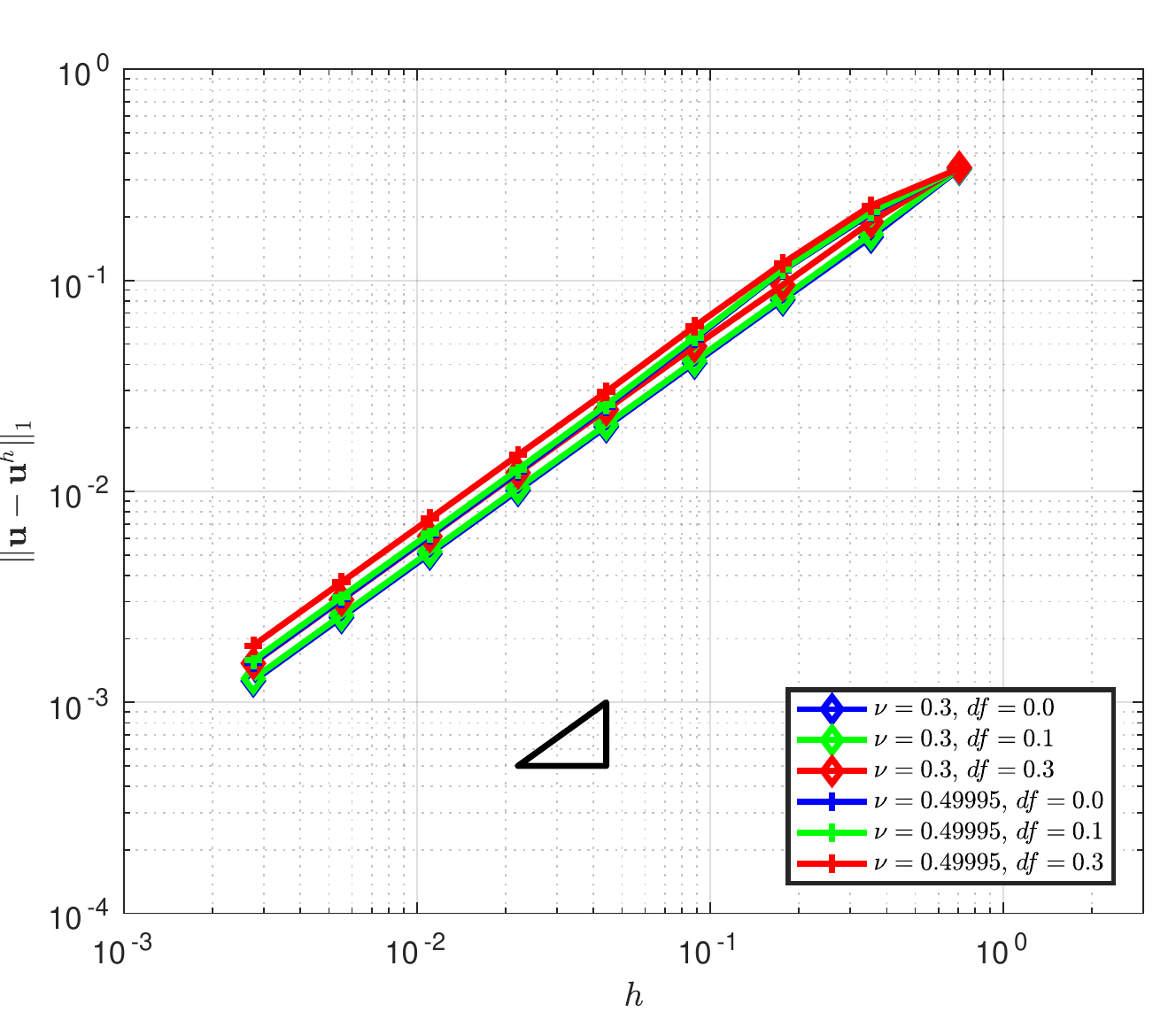}\label{SP_disp_error_NIPG}}
\caption{Square plate: IP displacement $H^1$ error convergence. The hypotenuse of the triangle has a slope of 1 in each case.}
\label{SP_disp_error_IP}
\end{figure}


\begin{figure}[!ht]
\centering
\subfloat[Exact solution]{\includegraphics[width=.20\columnwidth]{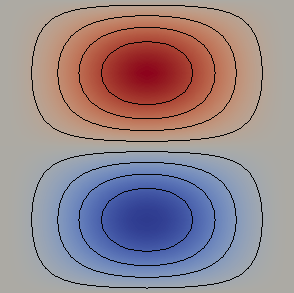}\label{SP_disp_exact}} \hspace{10mm}
\subfloat[SIPG, \textit{df} = 0.0, mesh 5]{\includegraphics[width=.20\columnwidth]{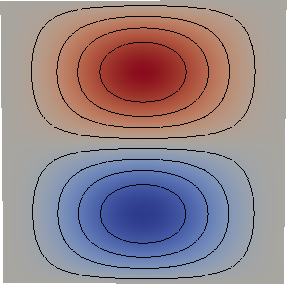}\label{SP_disp_SIPG_00}} \hspace{10mm}
\subfloat[SIPG, \textit{df} = 0.3, mesh 5]{\includegraphics[width=.27\columnwidth]{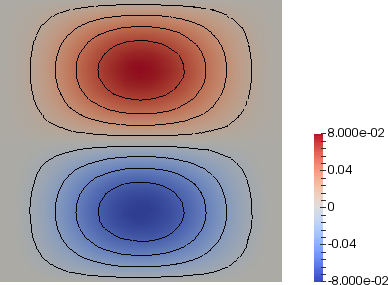}\label{SP_disp_SIPG_03}}
\caption{Square plate: displacement $u_x$, $\nu = 0.49995$}\label{SP_disp_x}
\end{figure}

The SG method with $Q_1$ elements (Figure \ref{SP_disp_error_SG1}) shows the same behaviour as for the beam example: optimal convergence when $\nu = 0.3$, with increasing error magnitude for increasing \textit{df}, and poor convergence, indicating locking, when $\nu = 0.49995$, irrespective of element shape regularity, until high refinement levels. In Figure \ref{SP_disp_error_SG_sri}, optimal convergence is seen for the same elements with SRI applied, for both values of $\nu$ and all \textit{df}, with error magnitudes increasing slightly as \textit{df} increases. The SG method with $Q_2$ elements (Figure \ref{SP_disp_error_SG2}) shows optimal convergence (second-order, in this case) when $\nu = 0.3$, irrespective of element regularity, though with decreasing accuracy as regularity decreases. With $\nu = 0.49995$, convergence is slower for the coarsest meshes but reaches (or exceeds) optimal convergence quickly with refinement.

Results from the original and new IP methods are shown in Figure \ref{SP_disp_error_IP} with performance similar to that in the beam example. In the near-incompressible case, the original IP methods converge poorly for coarse or medium-refinement meshes, depending on the variant, reaching optimal convergence earliest for the least regular meshes (Figure \ref{SP_disp_error_oldIPnu2}), while the new methods (Figure \ref{SP_disp_error_newIPnu2}) converge optimally even at low refinement levels. The error magnitudes are, overall, lower than those of the original methods. SG with SRI is again slightly more accurate in this example. The convergence of the new IP methods with respect to the compressibility parameter is uniform, as shown for NIPG in Figure \ref{SP_disp_error_NIPG}. 

The quality of the $x$-displacement approximation when $\nu = 0.49995$ is shown for SIPG in Figure \ref{SP_disp_x}, with a comparison to the exact solution, where it is evident that the accuracy of the method is not lost when general quadrilateral rather than rectangular elements are used. Performance for the component $u_y$ is the same, and these results extend to the other two IP methods.

\subsubsection{Post-processed stress}


\begin{figure}[!ht]
\centering
\subfloat[SG $Q_1$]{\includegraphics[trim={3.5cm 9cm 4cm 9cm},clip,width=.49\columnwidth]{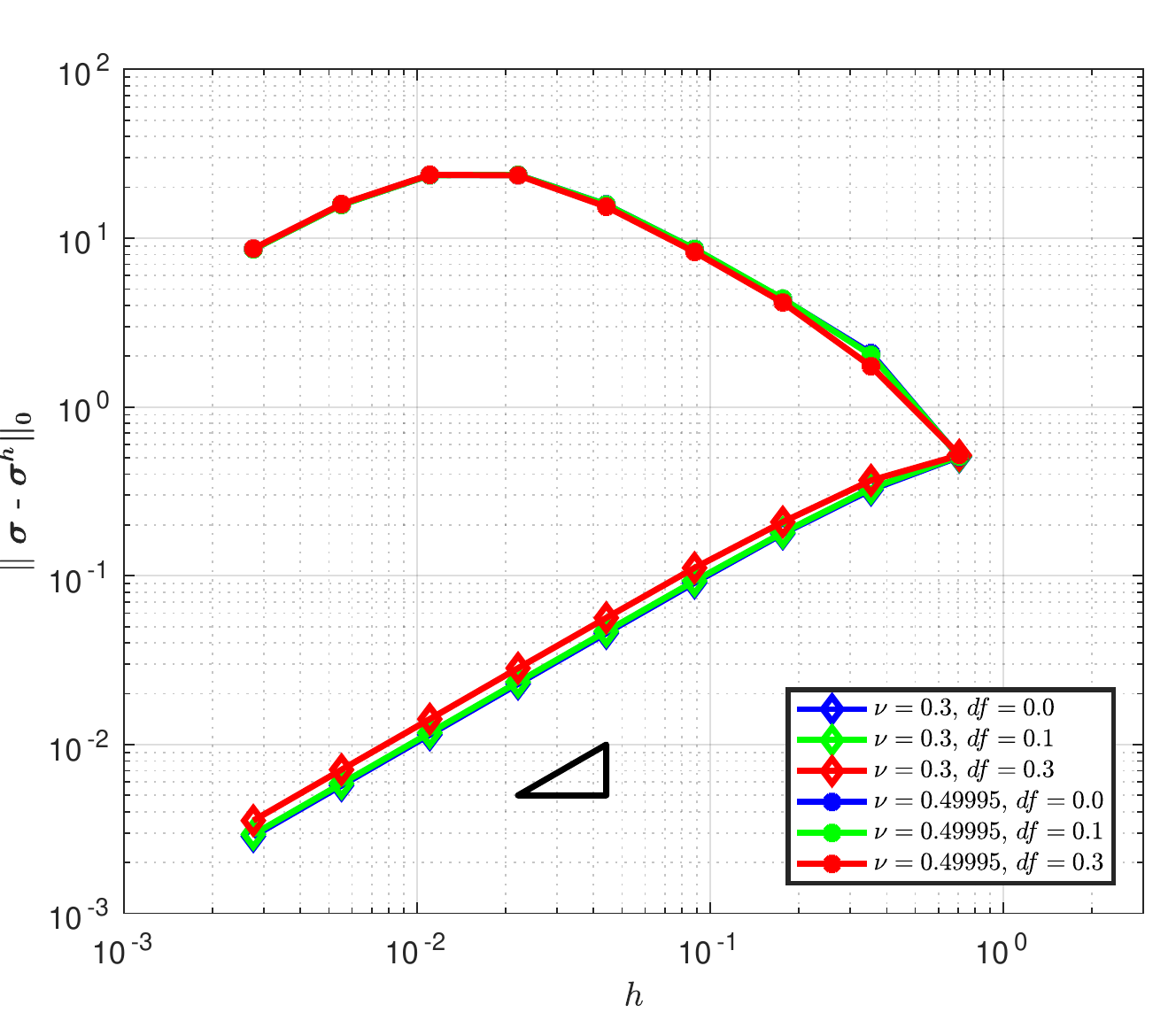}\label{SP_stress_error_SG1}}
\subfloat[SG $Q_1$ with SRI]{\includegraphics[trim={3.5cm 9cm 4cm 9cm},clip,width=.49\columnwidth]{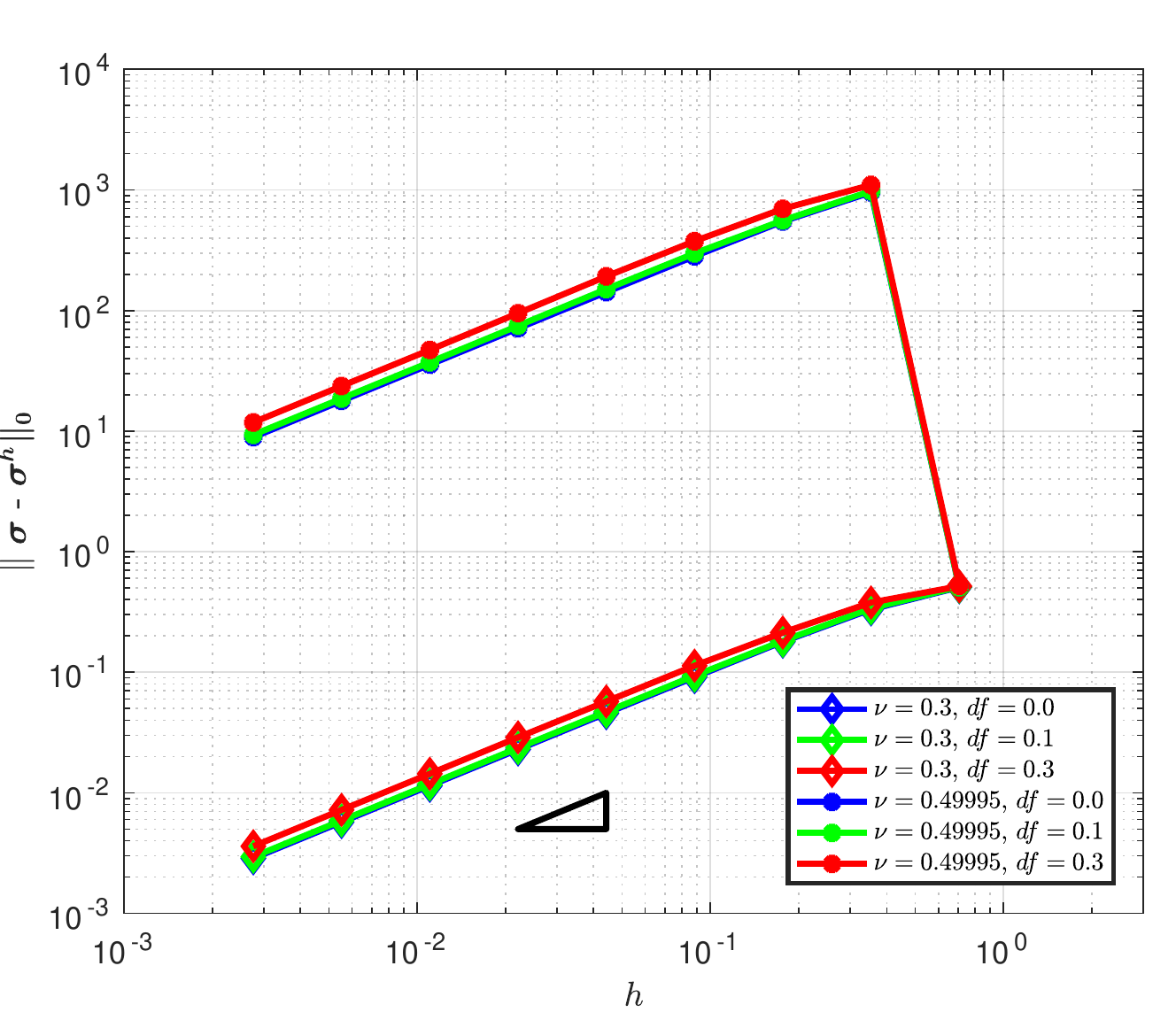}\label{SP_stress_error_SG_sri}}\\
\subfloat[SG $Q_2$]{\includegraphics[trim={3.5cm 9cm 4cm 9cm},clip,width=.49\columnwidth]{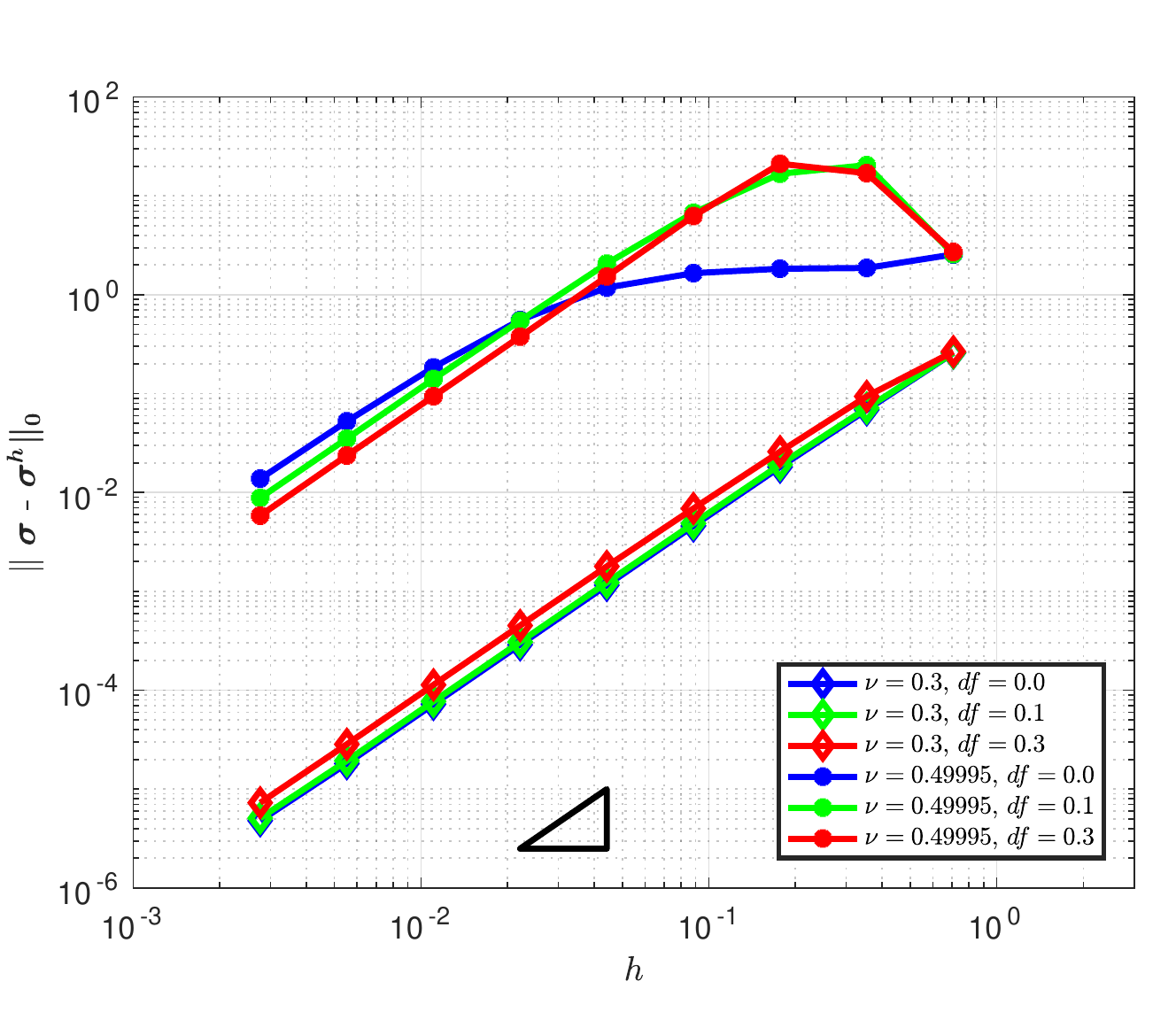}\label{SP_stress_error_SG2}}
\caption{Square plate: SG stress $L^2$ error convergence. The hypotenuse of the triangle in each of (a) and (b) has a slope of 1, and in (c) has a slope of 2.}\label{SP_stress_error_SG}
\end{figure}

\begin{figure}[!ht]
\centering
\subfloat[Original IP, $\nu = 0.49995$]{\includegraphics[trim={3.5cm 9cm 4cm 9cm},clip,width=.49\columnwidth]{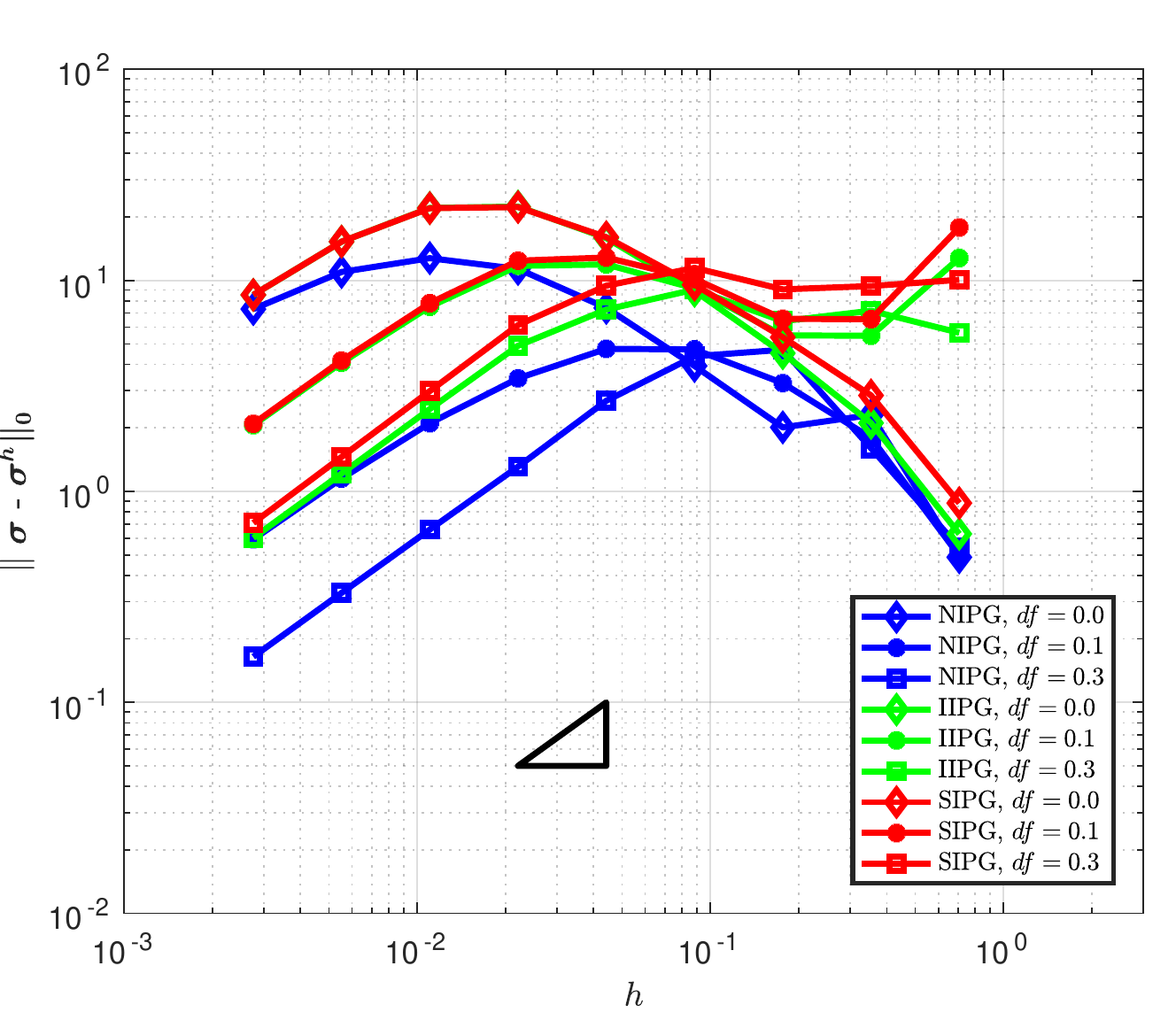}\label{SP_stress_error_oldIPnu2}}
\subfloat[New NIPG]{\includegraphics[trim={3.5cm 9cm 4cm 9cm},clip,width=.49\columnwidth]{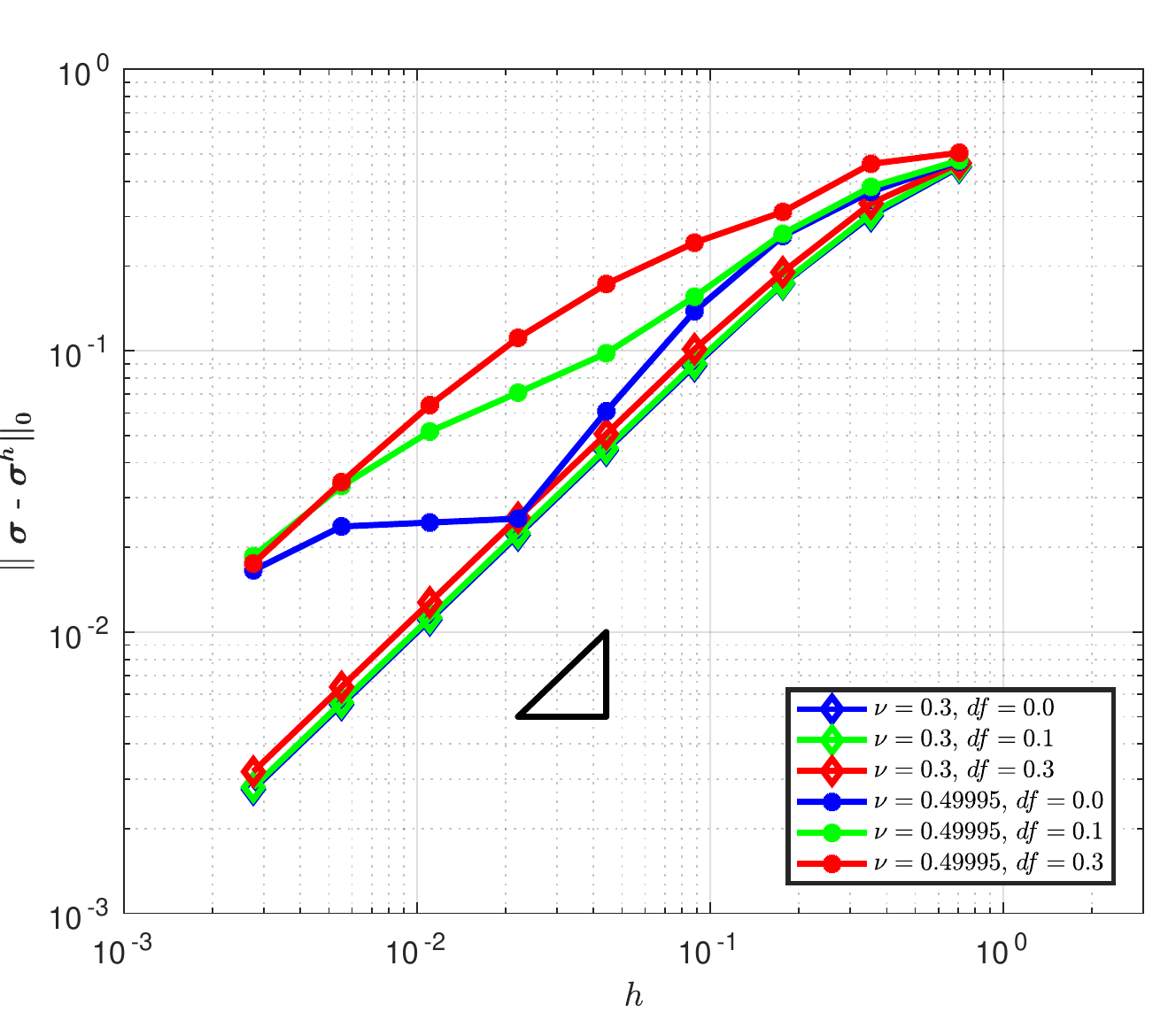}\label{SP_stress_error_NIPG}}\\
\subfloat[New IIPG]{\includegraphics[trim={3.5cm 9cm 4cm 9cm},clip,width=.49\columnwidth]{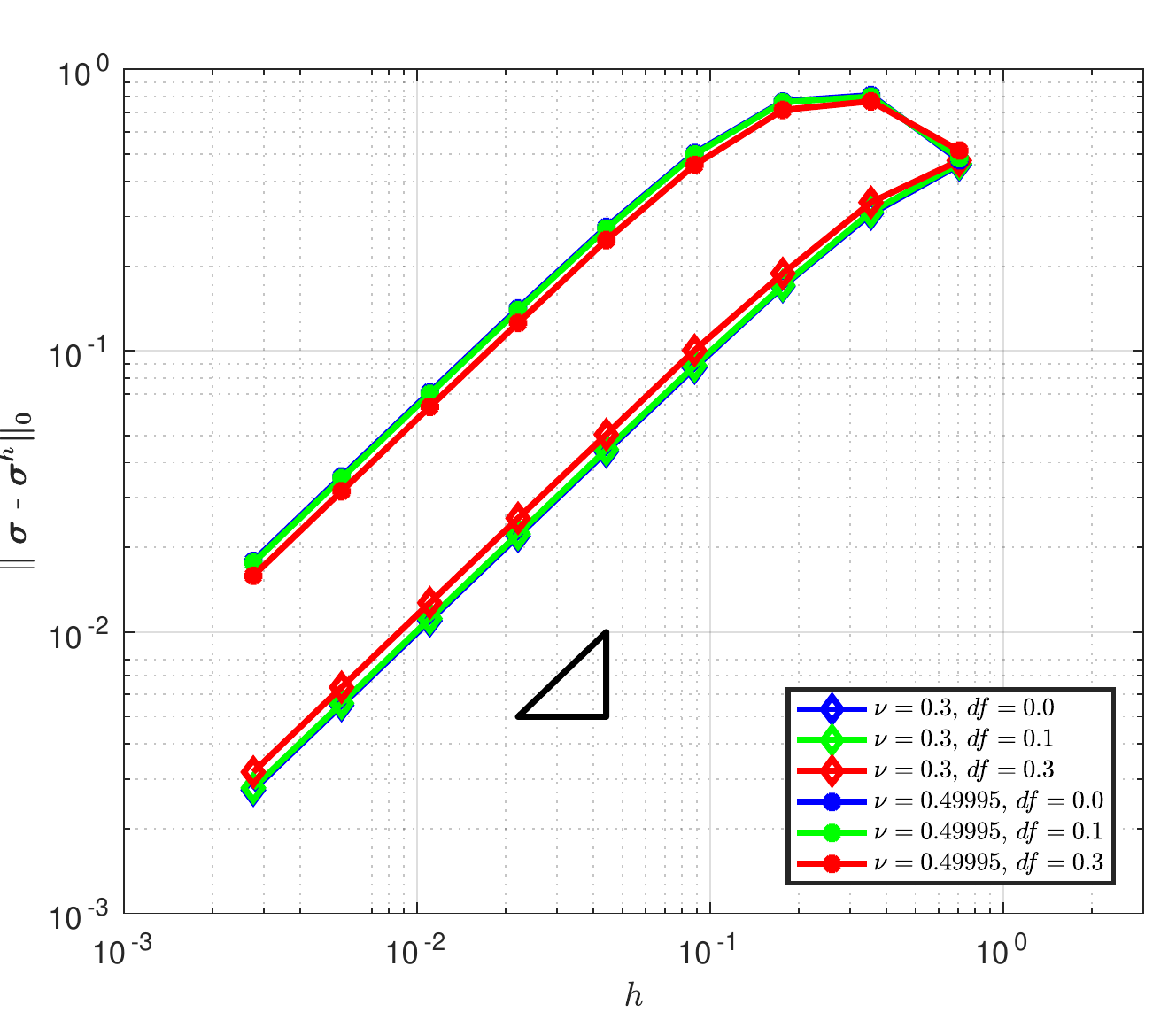}\label{SP_stress_error_IIPG}}
\subfloat[New IP, $\nu = 0.49995$]{\includegraphics[trim={3.5cm 9cm 4cm 9cm},clip,width=.49\columnwidth]{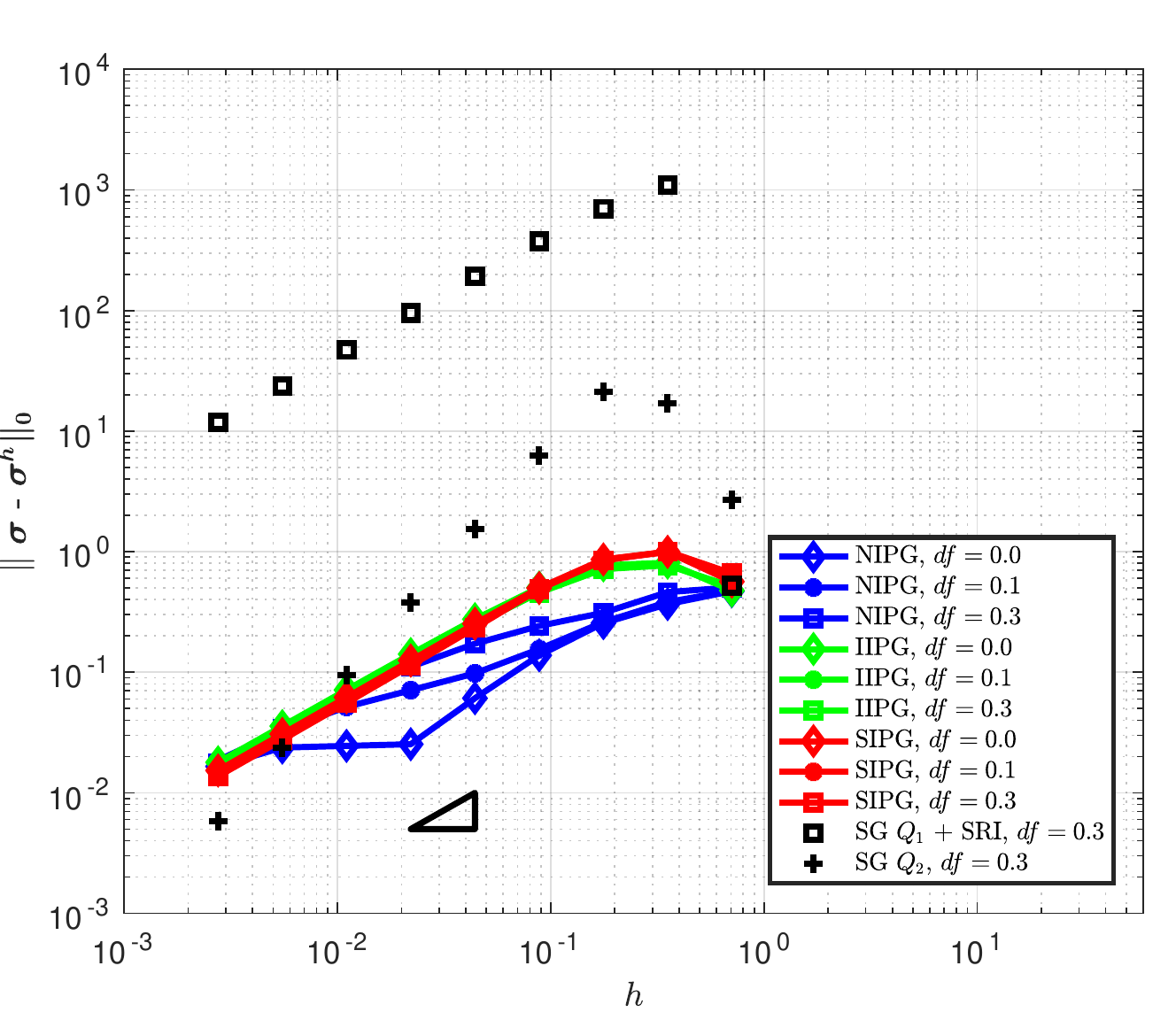}\label{SP_stress_error_newIPnu2}}
\caption{Square plate: IP stress $L^2$ error convergence. The hypotenuse of the triangle has a slope of 1 in each case.}
\end{figure}


\begin{figure}[!ht]
\centering
\subfloat[$\sigma_{xx}$]{\includegraphics[width=.27\columnwidth]{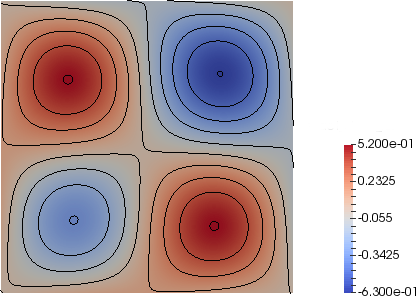}\label{SP_stress_xx_exact}}  \hspace{10mm}
\subfloat[$\sigma_{xy}$]{\includegraphics[width=.27\columnwidth]{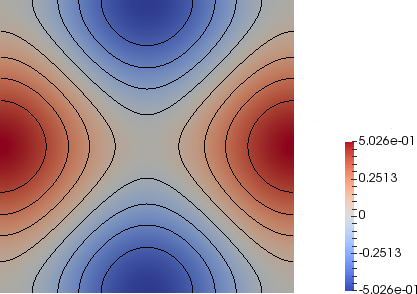}\label{SP_stress_xy_exact}}
\caption{Square plate: exact solution, $\nu = 0.49995$}\label{SP_stress_exact}
\end{figure}

\begin{figure}[!ht]
\centering
\subfloat[\textit{df} = 0.0, mesh 5]{\includegraphics[width=.18\columnwidth]{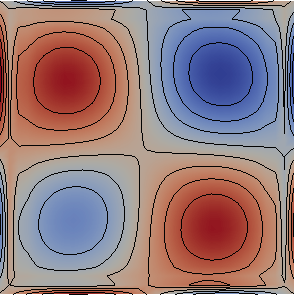}\label{SP_stress_xx_SIPG_00_m5}} \hspace{10mm}
\subfloat[\textit{df} = 0.1, mesh 5]{\includegraphics[width=.18\columnwidth]{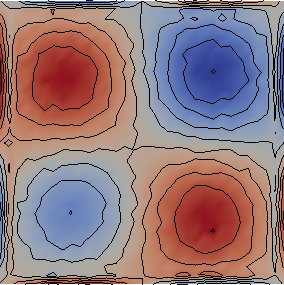}\label{SP_stress_xx_SIPG_01_m5}} \\
\subfloat[\textit{df} = 0.3, mesh 5]{\includegraphics[width=.18\columnwidth]{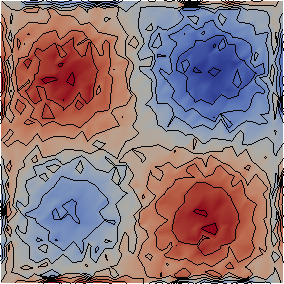}\label{SP_stress_xx_SIPG_03_m5}}\hspace{10mm}
\subfloat[\textit{df} = 0.3, mesh 9]{\includegraphics[width=.18\columnwidth]{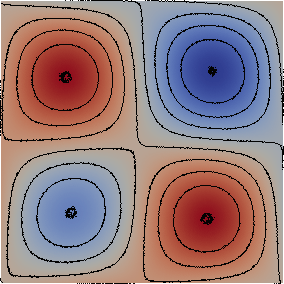}\label{SP_stress_xx_SIPG_03_m9}}\\
\hspace*{9mm}\subfloat[ \textit{df} = 0.3, mesh 5, without contour lines]{\includegraphics[width=.18\columnwidth]{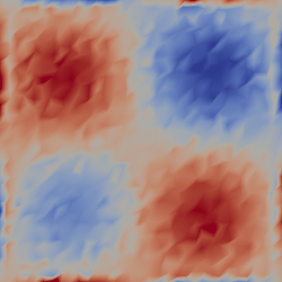}\label{SP_stress_xx_SIPG_03_m5_plain}}\hspace{10mm}
\subfloat[\textit{df} = 0.3, mesh 9, without contour lines]{\includegraphics[width=.23\columnwidth]{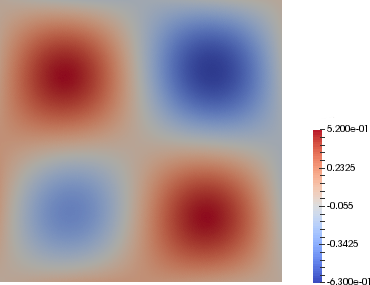}\label{SP_stress_xx_SIPG_03_m9_plain}}
\caption{Square plate: $\sigma_{xx}$, SIPG, $\nu = 0.49995$}\label{SP_stress_xx_SIPG}
\end{figure}

\begin{figure}[!ht]
\centering
\subfloat[\textit{df} = 0.0, mesh 5]{\includegraphics[width=.20\columnwidth]{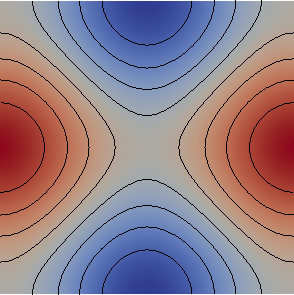}\label{SP_stress_xy_SIPG_00_m5}} \hspace{10mm}
\subfloat[\textit{df} = 0.3, mesh 5]{\includegraphics[width=.20\columnwidth]{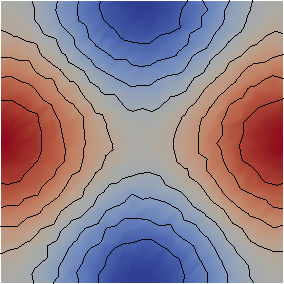}\label{SP_stress_xy_SIPG_03_m5}} \\
\hspace*{9mm}\subfloat[\textit{df} = 0.3, mesh 5, without contour lines]{\includegraphics[width=.20\columnwidth]{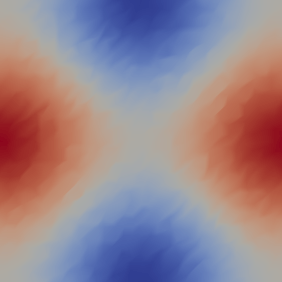}\label{SP_stress_xy_SIPG_03_m5}}\hspace{10mm}\subfloat[\textit{df} = 0.3, mesh 9]{\includegraphics[width=.26\columnwidth]{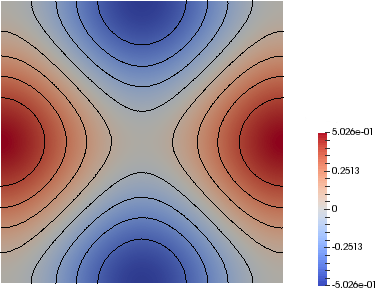}\label{SP_stress_xy_SIPG_03_m9}}
\caption{Square plate: $\sigma_{xy}$, SIPG, $\nu = 0.49995$}\label{SP_stress_xy_SIPG}
\end{figure}

The stress $L^2$ errors for the SG method with $Q_1$ elements (Figure \ref{SP_stress_error_SG1}) show similar convergence behaviour to those in the beam example: with $\nu = 0.3$, convergence is first-order, and for near-incompressibility convergence is poor until high refinements. With SRI (Figure \ref{SP_stress_error_SG_sri}), convergence is first-order irrespective of the compressibility level, but with $\nu = 0.49995$ the error magnitudes are significantly greater than with $\nu = 0.3$. SG with $Q_2$ elements (Figure \ref{SP_stress_error_SG2}) gives second-order convergence for $\nu = 0.3$; for $\nu = 0.49995$ convergence is also close to second-order except on the coarsest meshes. In most cases, the SG methods show slightly lower accuracy as \textit{df} increases.

For near-incompressibility, the original IP methods (Figure \ref{SP_stress_error_oldIPnu2}) give poor convergence for coarse meshes, and the level of refinement at which first-order convergence is attained varies, as in the beam example, with method and element regularity. The new SIPG and IIPG methods (Figure \ref{SP_stress_error_newIPnu2}) reach first-order convergence at low refinement levels when $\nu = 0.49995$, while  NIPG has rates that vary with element regularity, but with errors smaller than the other two methods because of faster convergence for the coarser meshes. This figure also shows the comparative magnitude for \textit{df} = 0.3 of the errors of SG $Q_1$ with SRI and of SG $Q_2$, the former significantly larger than the errors of the IP methods and the latter larger until high refinement. Figures \ref{SP_stress_error_NIPG} and \ref{SP_stress_error_IIPG} show enlarged plots of the NIPG and IIPG results for near-incompressibility, as well as the results for $\nu = 0.3$, where in both cases the first-order convergence is clear. 

The accuracy of the individual components $\sigma_{xx}$ and $\sigma_{xy}$ of the postprocessed stress for the near-incompressible case can be seen in Figures \ref{SP_stress_xx_SIPG} and \ref{SP_stress_xy_SIPG} for SIPG, with the exact solution in Figure \ref{SP_stress_exact} as a comparison. (The stress pattern for the component $\sigma_{yy}$ is qualitatively similar to that of $\sigma_{xx}$ and the SIPG results are likewise comparable.) In $\sigma_{xx}$, here as for the beam example the smoothness of the field decreases as \textit{df} increases, but mesh refinement has a great effect in smoothing the field. The $\sigma_{xy}$ field, in contrast, shows very little deterioration in accuracy, as \textit{df} increases, the smoothness again increasing with refinement. These different behaviours stem from a lack of smoothness in the field of $tr\,\bvarepsilon$ for this boundary value problem: although the order of magnitude of $tr\,\bvarepsilon$ is low, it is amplified by the factor of $\lambda$ in calculating the direct stresses (see equation (\ref{const_rel})). For near-incompressibility, $\lambda$ is large ($O(10^3)$) and thus the contribution of this term (which doesn't appear in the shear component) is significant.

IIPG performs similarly to SIPG, while NIPG is similar but additionally does not show the inaccuracy around the edges of the domain in the direct stress approximations, as the other two methods do for medium refinement (eg. mesh 5).

\subsection{L-shaped domain}

\begin{figure}[!ht]
\centering
\subfloat[Domain]{\includegraphics[width=.44\columnwidth]{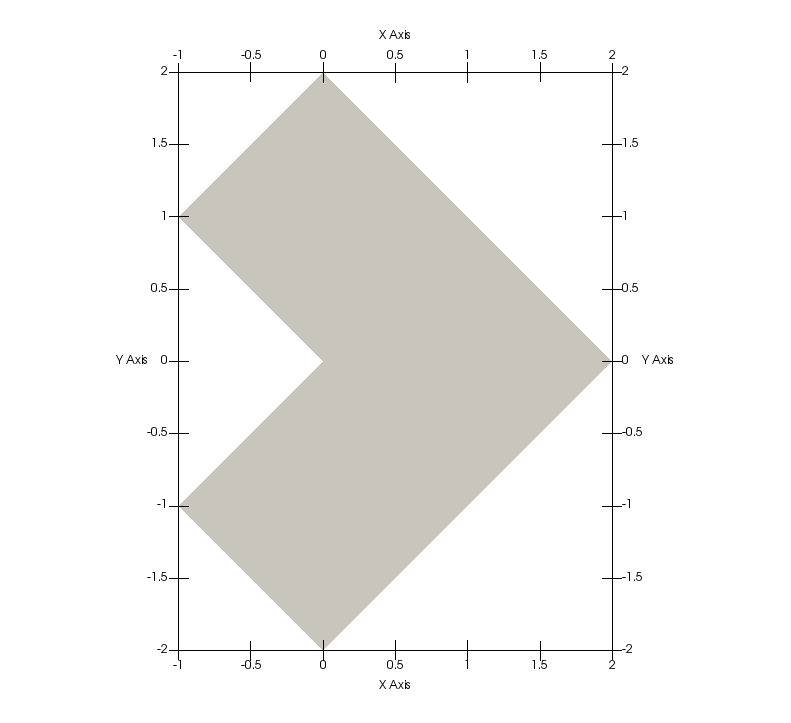}
\label{LS_domain}}
\subfloat[Domain and solution (deformation scaled $\times 3000$)]{\includegraphics[width=.44\columnwidth]{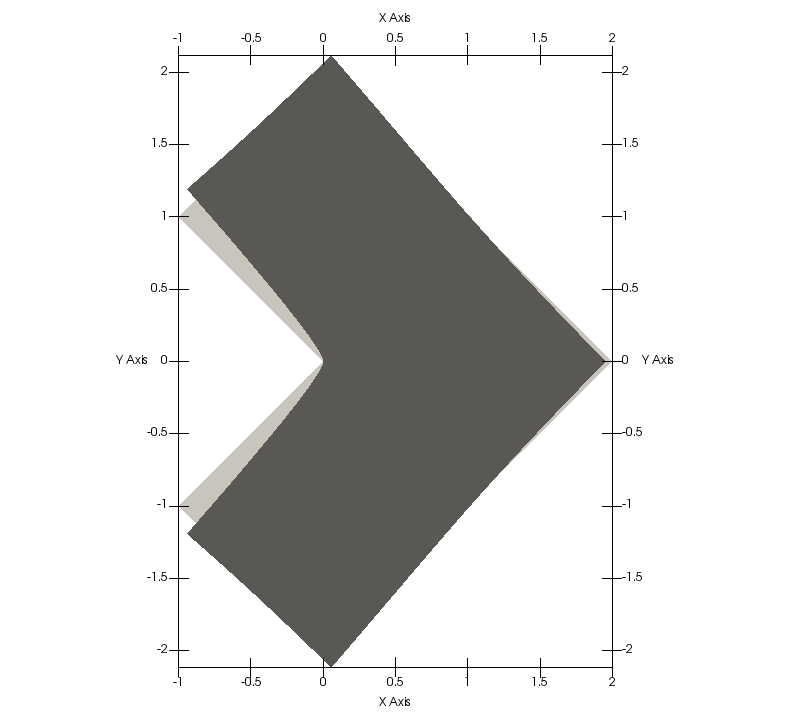}\label{LS_domain_and_solution}}
\caption{L-shaped domain}
\end{figure}

The L-shaped domain shown in Figure \ref{LS_domain} is displaced to the position shown in Figure \ref{LS_domain_and_solution} by imposing Dirichlet boundary conditions as detailed in \cite{Alberty2002}. The Young's modulus is $E = 100000$. The prescription of Dirichlet boundary conditions on the entire boundary precludes the possibility of volumetric locking in the numerical solution of this problem, but the effect of a high compressibility parameter on the quality of the displacement approximation is nevertheless of interest.

In this problem, with the SIPG method the stabilization parameter values $k_{\mu} = 10$ and $k_{\lambda} = 50$ are used. This combination is introduced here for its effect on the stability of the stress results for the new SIPG method.

\begin{figure}[!ht]
\centering
\subfloat[Mesh]{\includegraphics[width=.49\columnwidth]{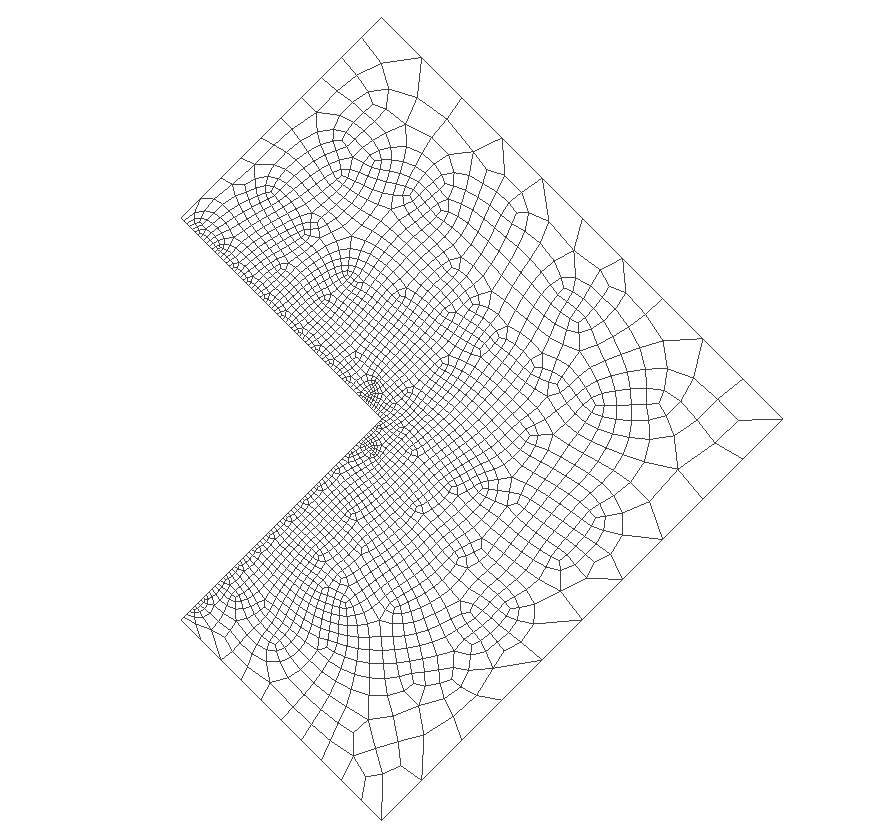} \label{LS_mesh1}}
\subfloat[Measure of deviation of shape of element from that of a square element]{\includegraphics[width=.49\columnwidth]{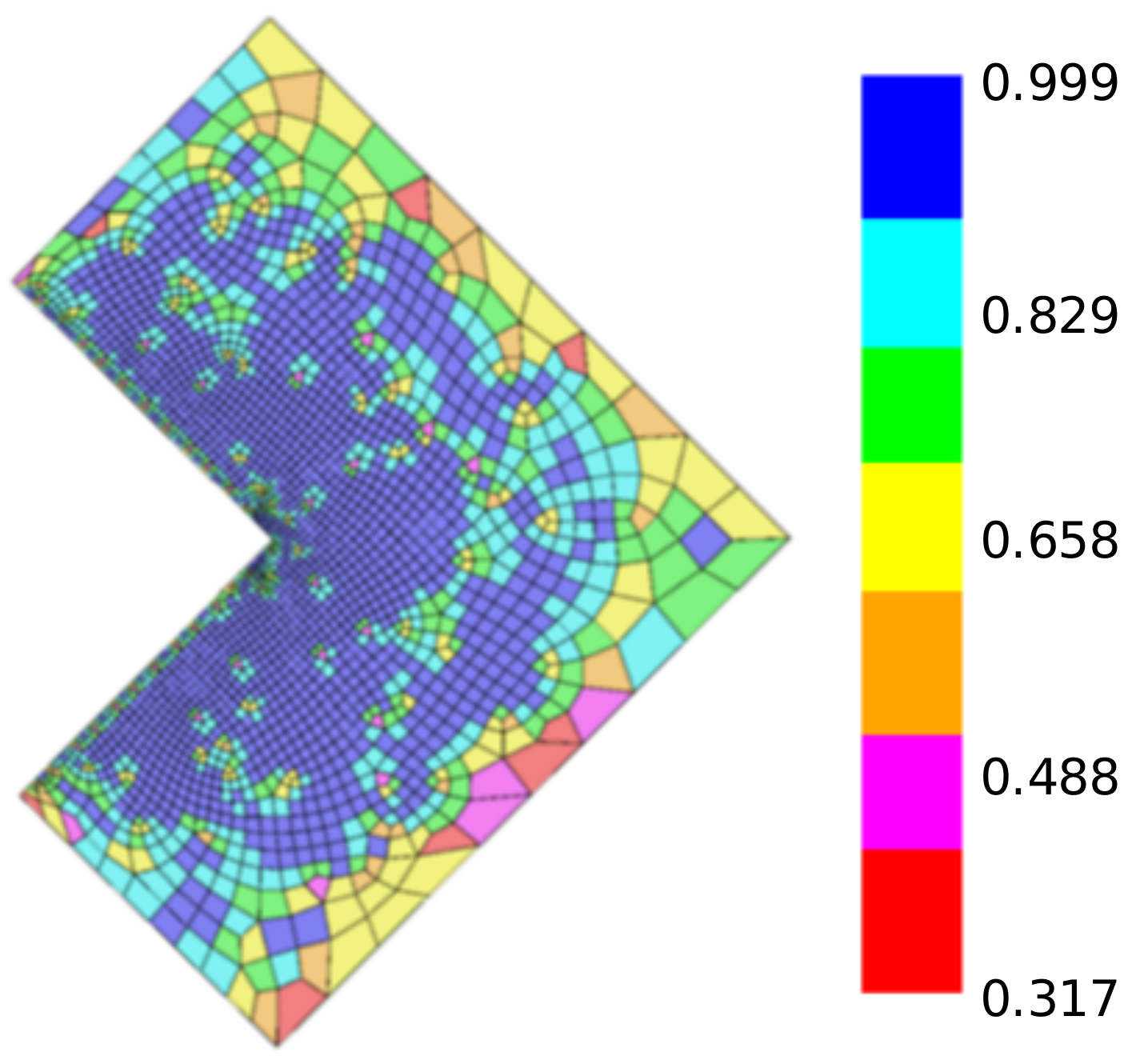}
\label{LS_mesh_quality}}
\caption{L-shaped domain: mesh 1}
\label{LS_mesh}
\end{figure}

The initial mesh used is unstructured and graded, and finest around the re-entrant corner (Figure \ref{LS_mesh1}). While a large portion of the quadrilateral elements in this mesh are close to square, many are of lower regularity and some are very poor in shape (Figure \ref{LS_mesh_quality})\footnote{The measure of deviation in Figure \ref{LS_mesh_quality} is calculated as 2 divided by the condition number of the weighted Jacobian matrix, by the CUBIT metric ``Shape" for quadrilateral elements \cite{CubitManual}.}. The initial mesh is then refined, and Table \ref{mesh_info_LS} details the number of elements and dofs for each method at the various refinement levels.

\begin{table}[!ht]
\begin{center}
\begin{tabular}{|c|c|c|c|c|} \hline
$n$ & no.\ els & \multicolumn{3}{|c|}{no.\ dofs} \\ \hline 
& & SG $Q_1$ & SG $Q_2$ & IP \\ \hline
1 & 2303 & 4796 & 18802 & 18424\\ \hline
2 & 9212 & 18802 & 74450 & 73696\\ \hline
3 & 36848 & 74450 & 296290 & 294784\\ \hline
4 & 147392 & 296290 & 1182146 & 1179136\\ \hline
5 & 589568 & 1182146 & 4722562 & 4716544\\ \hline
\end{tabular}
\end{center}
\caption{Mesh details for meshes 1 to 5, for standard Galerkin (SG) and interior penalty (IP) methods for the L-shaped domain example: number of elements and number of degrees of freedom are shown.}\label{mesh_info_LS}\end{table}

\textbf{Note on domain regularity:} The analytical results of \cite{Grieshaber2015} for rectangular elements rely on a regularity estimate that assumes a degree of domain regularity not satisfied by the nonconvex polygon of this example (see \cite{Brenner2002}). (Note that Wihler \cite{Wihler2004} has developed an extension of this result for general polygonal domains within the framework of weighted Sobolev spaces.) For conforming finite elements, the theoretically predicted rate of convergence for a uniform mesh is $\frac{2}{3}$, based on the interior angle of $\frac{3 \pi}{2}$ (see \cite{Carstensen2001} and references therein). We, however, use a graded mesh and do not have a theoretical prediction for the convergence rate.

At the origin, there is a stress singularity in the exact solution, which makes this a challenging benchmark boundary value problem for testing computational approximation methods.

\subsubsection{Displacement approximation}


\begin{figure}[!ht]
\centering
\subfloat[SG]{\includegraphics[trim={3.5cm 9cm 4cm 9cm},clip,width=.49\columnwidth]{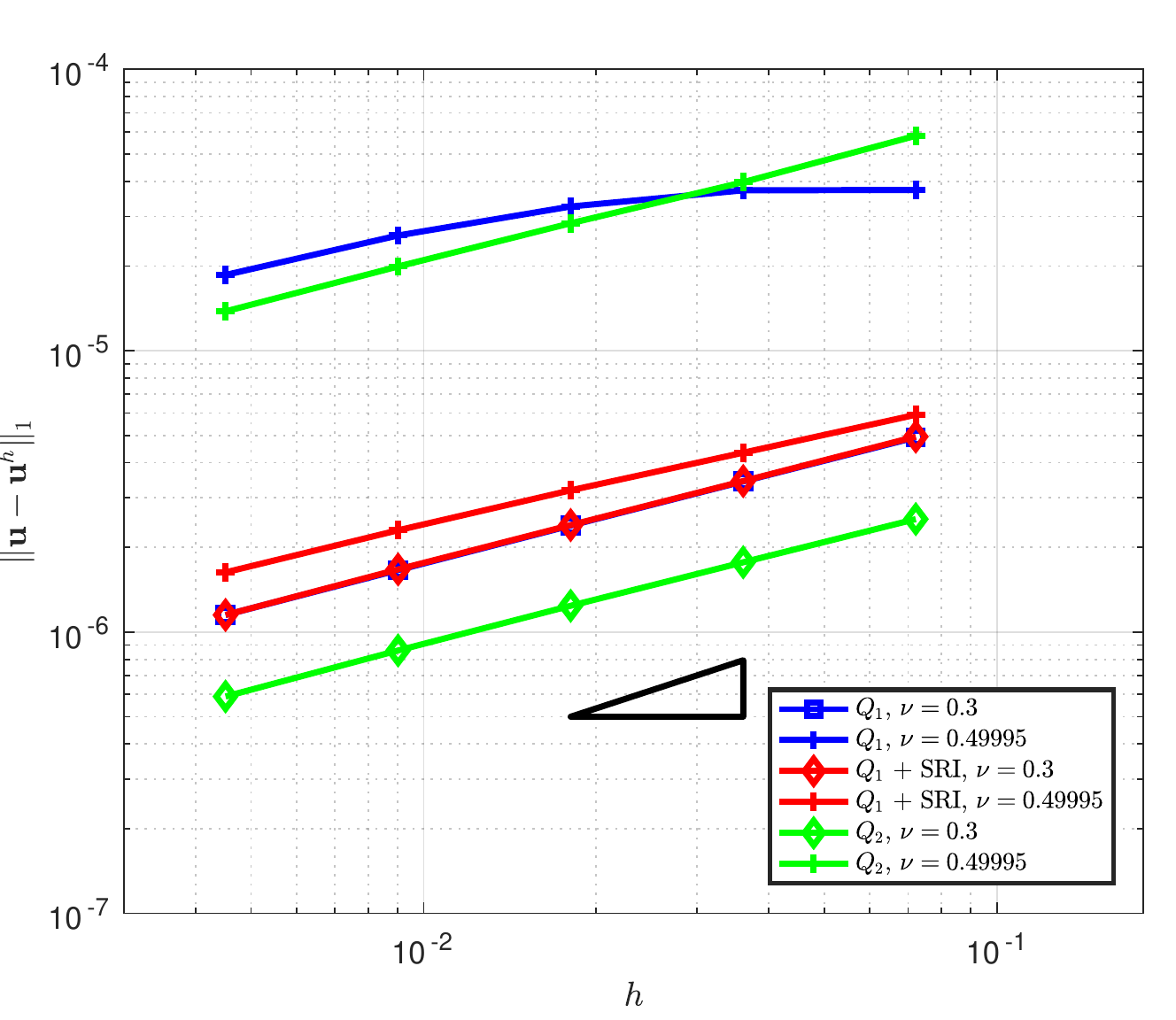}
\label{LS_disp_error_SG}}
\subfloat[Original IP]{\includegraphics[trim={3.5cm 9cm 4cm 9cm},clip,width=.49\columnwidth]{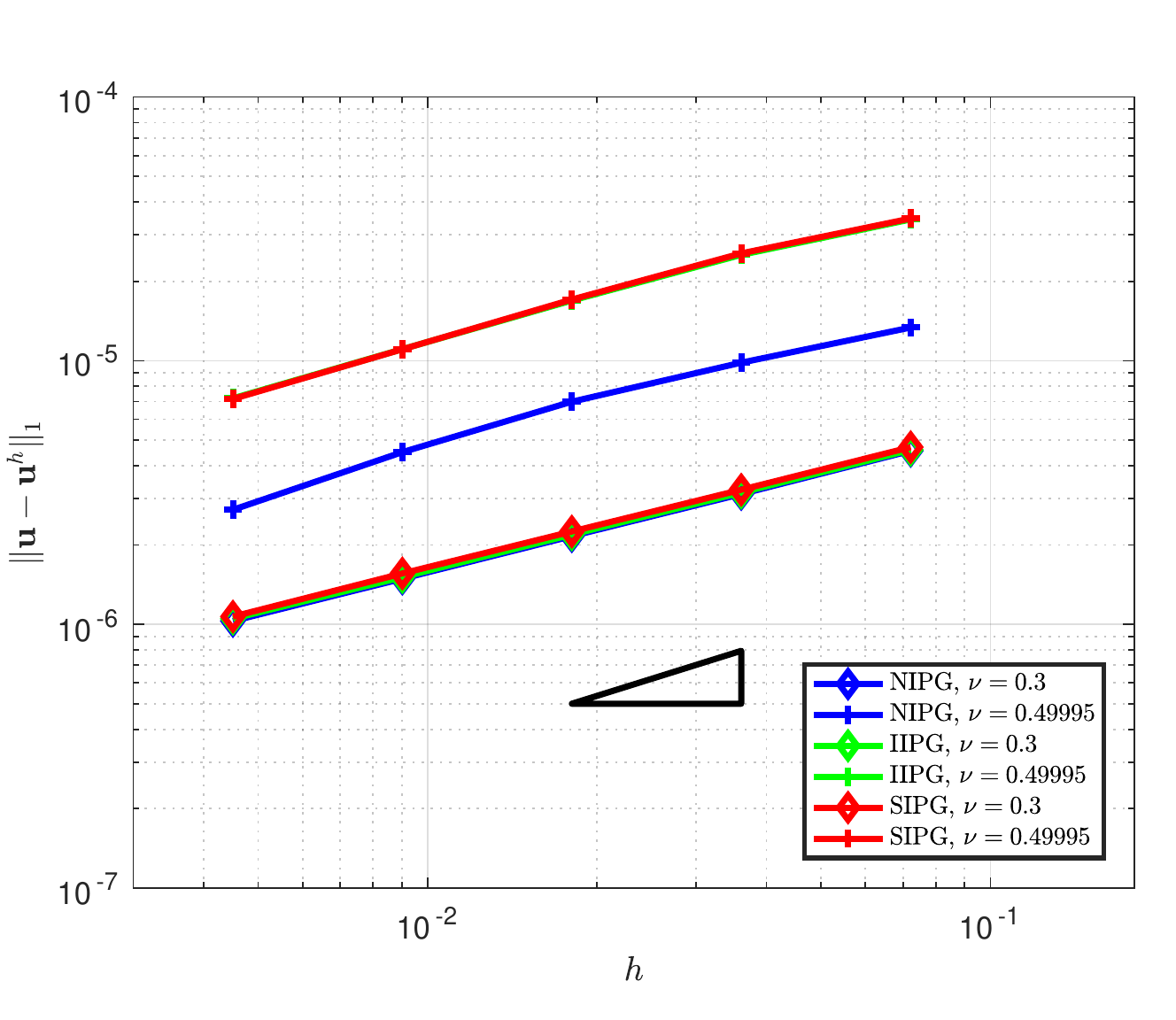}\label{LS_disp_error_IPold}}\\
\subfloat[New IP]{\includegraphics[trim={3.5cm 9cm 4cm 9cm},clip,width=.49\columnwidth]{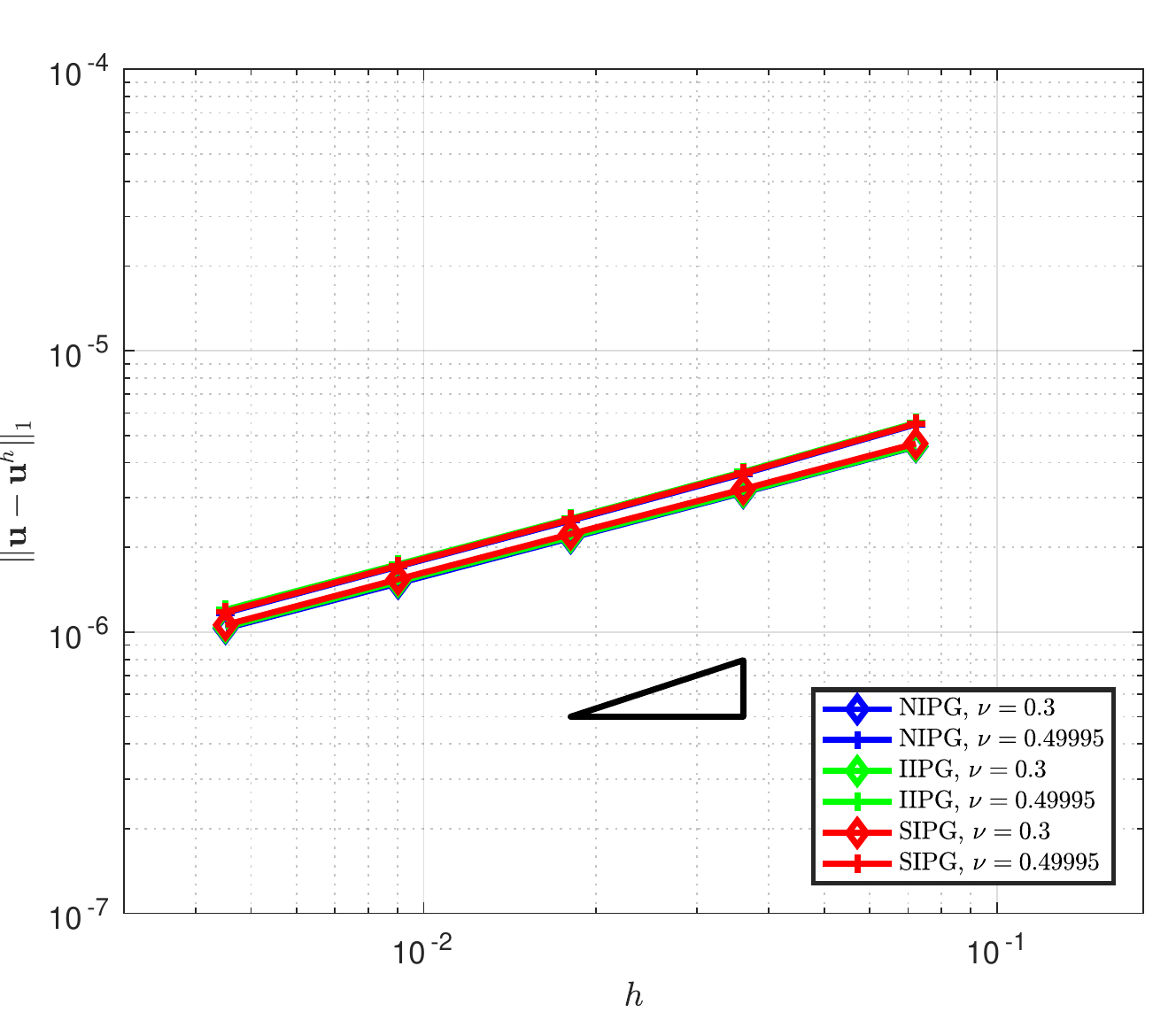}\label{LS_disp_error_IPnew}}
\subfloat[Various  IP and SG, $\nu = 0.49995$]{\includegraphics[trim={3.5cm 9cm 4cm 9cm},clip,width=.49\columnwidth]{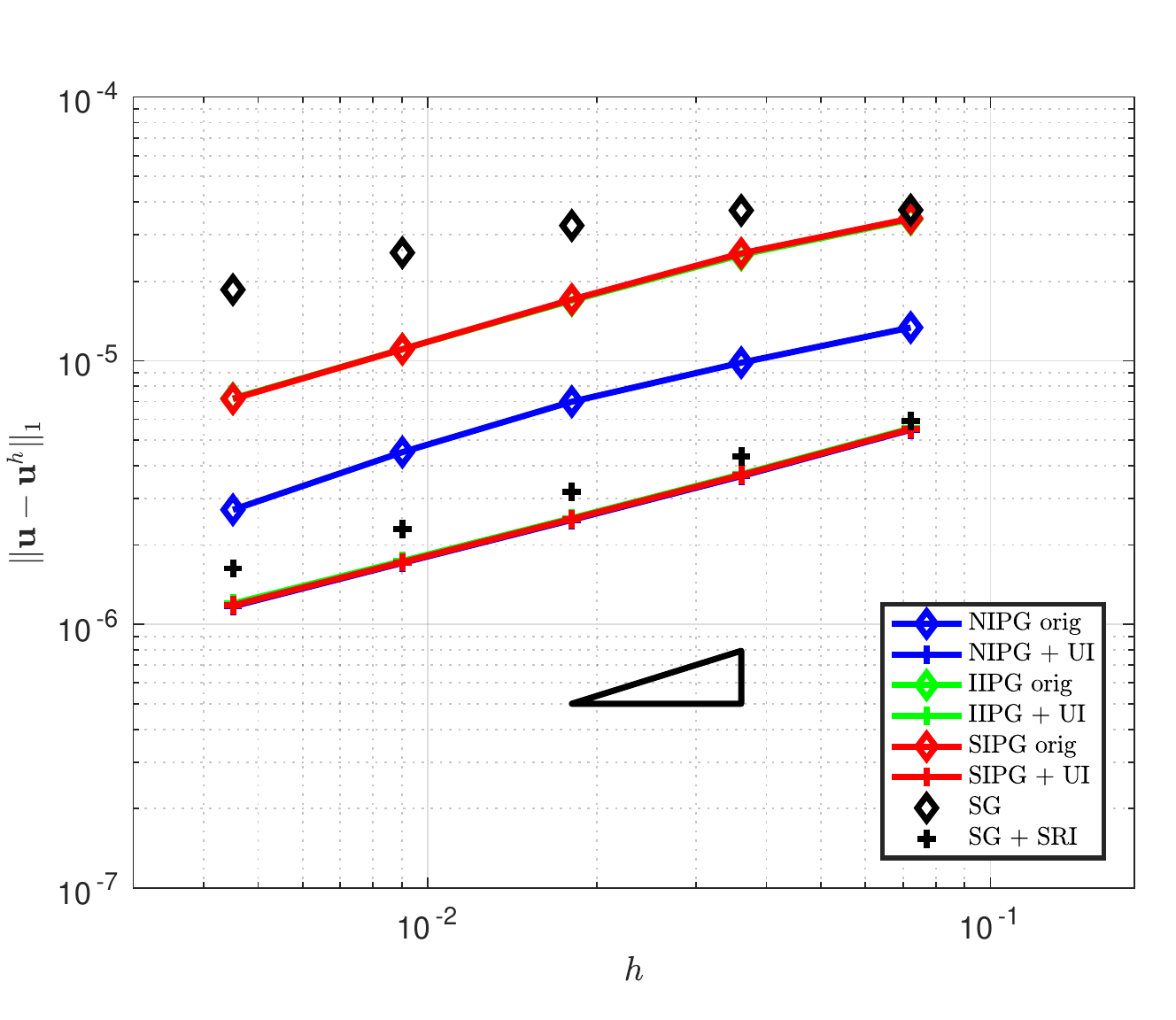}\label{LS_disp_error_various_nu2}}
\caption{L-shaped domain: Displacement $H^1$ error convergence. The hypotenuse of the triangle has a slope of $2/3$ in each case.}
\end{figure}


\begin{figure}[!ht]
\centering
\subfloat[Exact solution]{\includegraphics[width=.2\columnwidth]{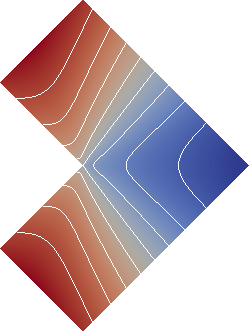}
\label{LS_disp_x_exact}}\hspace{15mm}
\subfloat[SG $\Q_1$, mesh 1]{\includegraphics[width=.2\columnwidth]{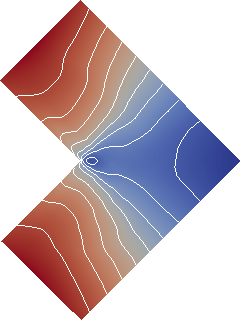}
\label{LS_disp_x_SG1_m1}}\hspace{15mm}
\subfloat[SIPG, mesh 1]{\includegraphics[width=.3\columnwidth]{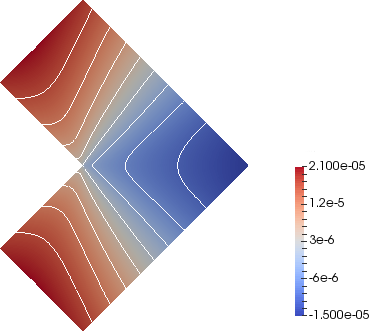}
\label{LS_disp_x_N_m1}}
\caption{L-shaped domain: Displacement $u_x$}
\end{figure}

Figure \ref{LS_disp_error_SG} shows the results of various cases of the SG method. With $Q_1$ elements, there is a convergence rate of 0.54 for the compressible case. The results are indistinguishable for $Q_1$ with SRI, and with $Q_2$ elements the errors are smaller but the convergence rate is the same. For near-incompressibility, the rates decrease significantly for $Q_1$ elements, and very slightly for SRI. With $Q_2$ elements, the rate is essentially unaffected but the error magnitudes are higher.

The original IP methods (Figure \ref{LS_disp_error_IPold}) attain a convergence rate of 0.54, when $\nu = 0.3$, and with $\nu = 0.49995$ the rate is slightly higher although the errors are greater. The new IP methods (Figure \ref{LS_disp_error_IPnew}) have a consistent rate of 0.54, for both values of $\nu$, with comparatively little increase in error for the near-incompressible case. Figure \ref{LS_disp_error_various_nu2} shows the higher accuracy of the new methods than of the original for $\nu = 0.49995$, as well as that of the SG method with and without SRI. 

The poor displacement approximation of the SG method ($x$-displacement shown in Figure \ref{LS_disp_x_SG1_m1}) when $\nu= 0.49995$ is not an example of locking but of the inaccuracy of the deformation over the domain, and particularly in the region of the origin, as can be seen by comparison to the exact solution in Figure \ref{LS_disp_x_exact}. The three  IP approximations in this case are excellent (SIPG result shown in Figure \ref{LS_disp_x_N_m1}). The corresponding $y$-displacements display similar accuracy in each case.

\subsubsection{Postprocessed stress}


\begin{figure}[!ht]
\centering
\subfloat[SG]{\includegraphics[trim={3.5cm 9cm 4cm 9cm},clip,width=.49\columnwidth]{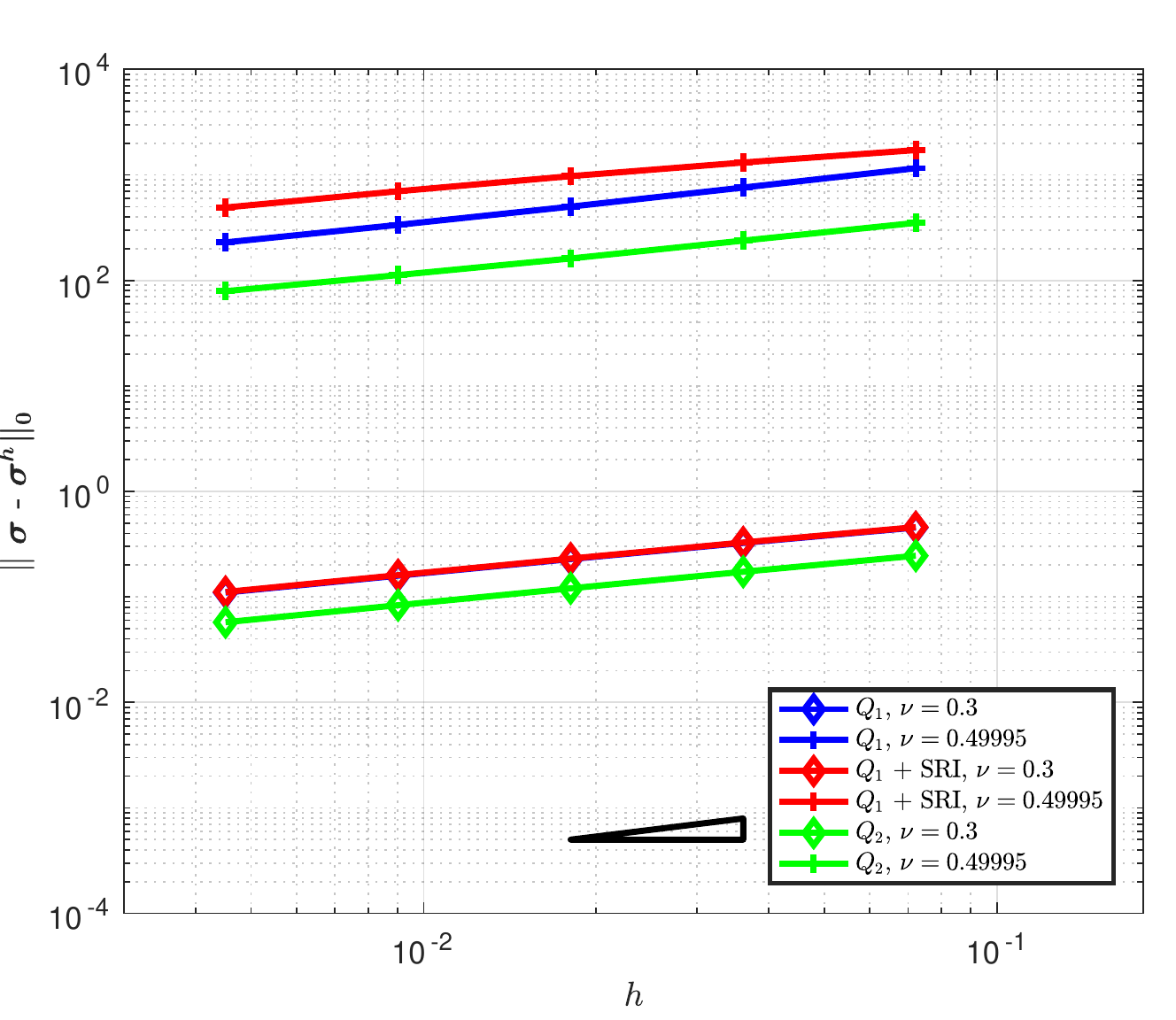}
\label{LS_stress_error_SG}}
\subfloat[Original IP]{\includegraphics[trim={3.5cm 9cm 4cm 9cm},clip,width=.49\columnwidth]{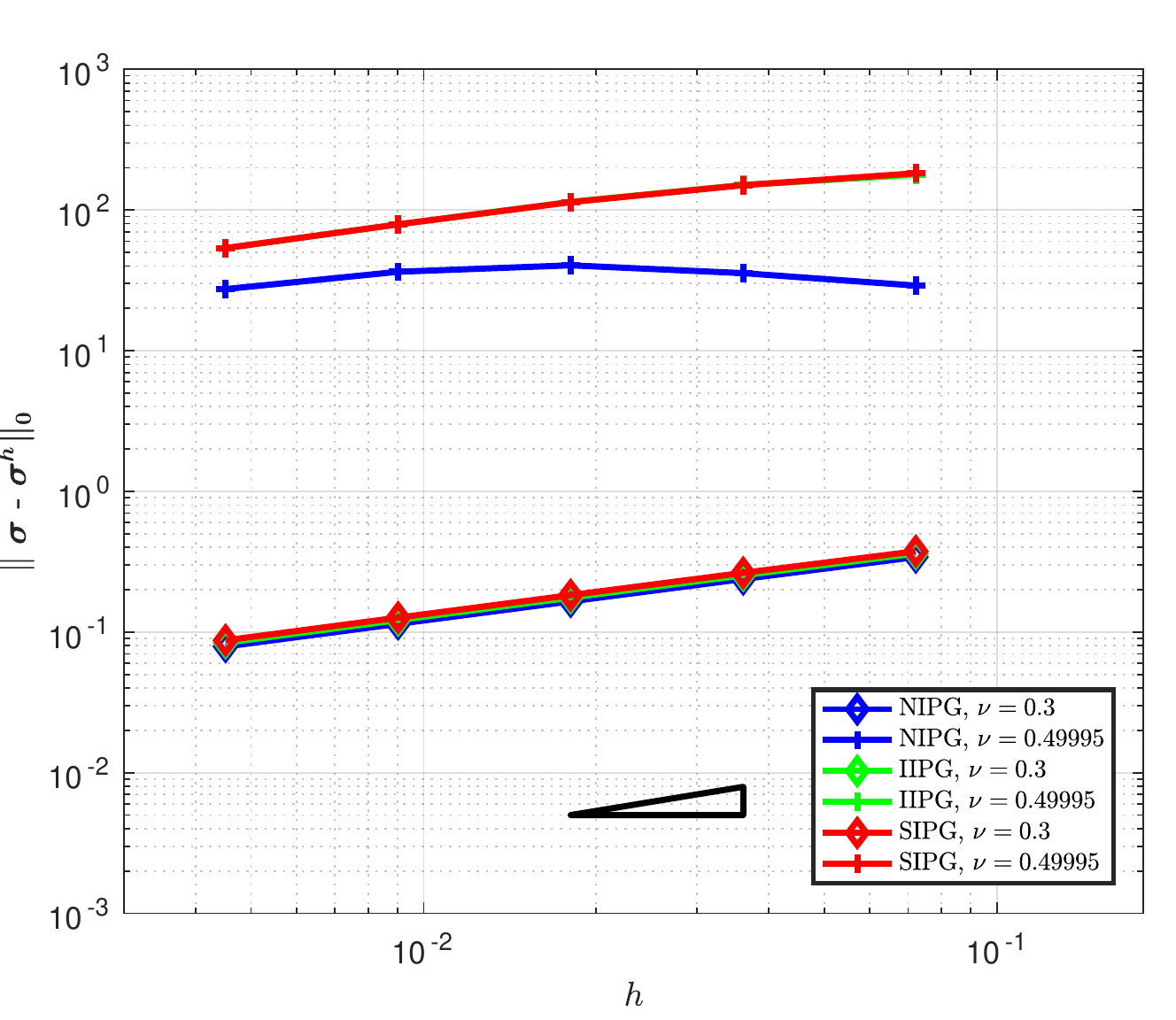}\label{LS_stress_error_IPold}}\\
\subfloat[New IP]{\includegraphics[trim={3.5cm 9cm 4cm 9cm},clip,width=.49\columnwidth]{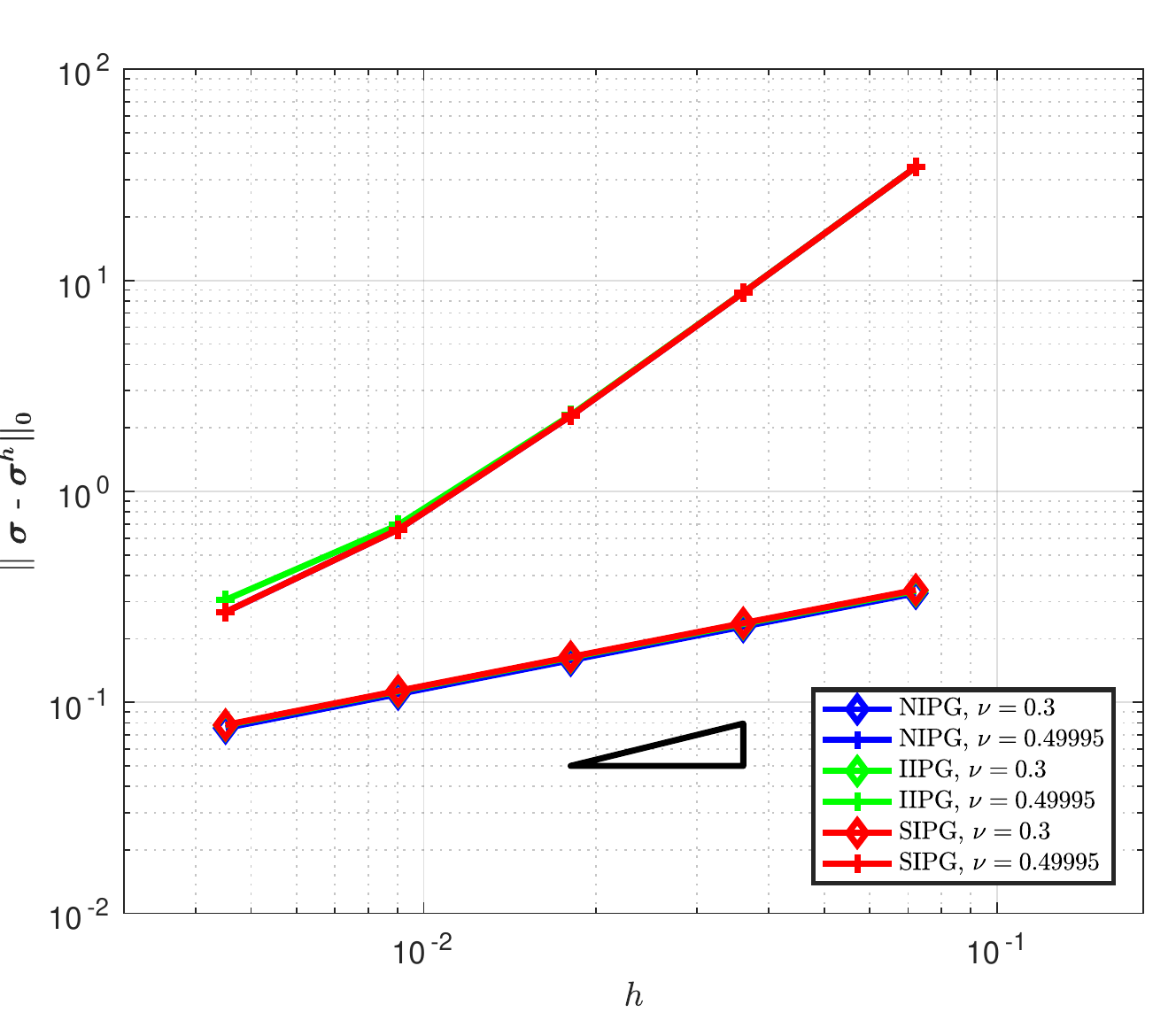}\label{LS_stress_error_IPnew}}
\subfloat[Various IP and SG, $\nu = 0.49995$]{\includegraphics[trim={3.5cm 9cm 4cm 9cm},clip,width=.49\columnwidth]{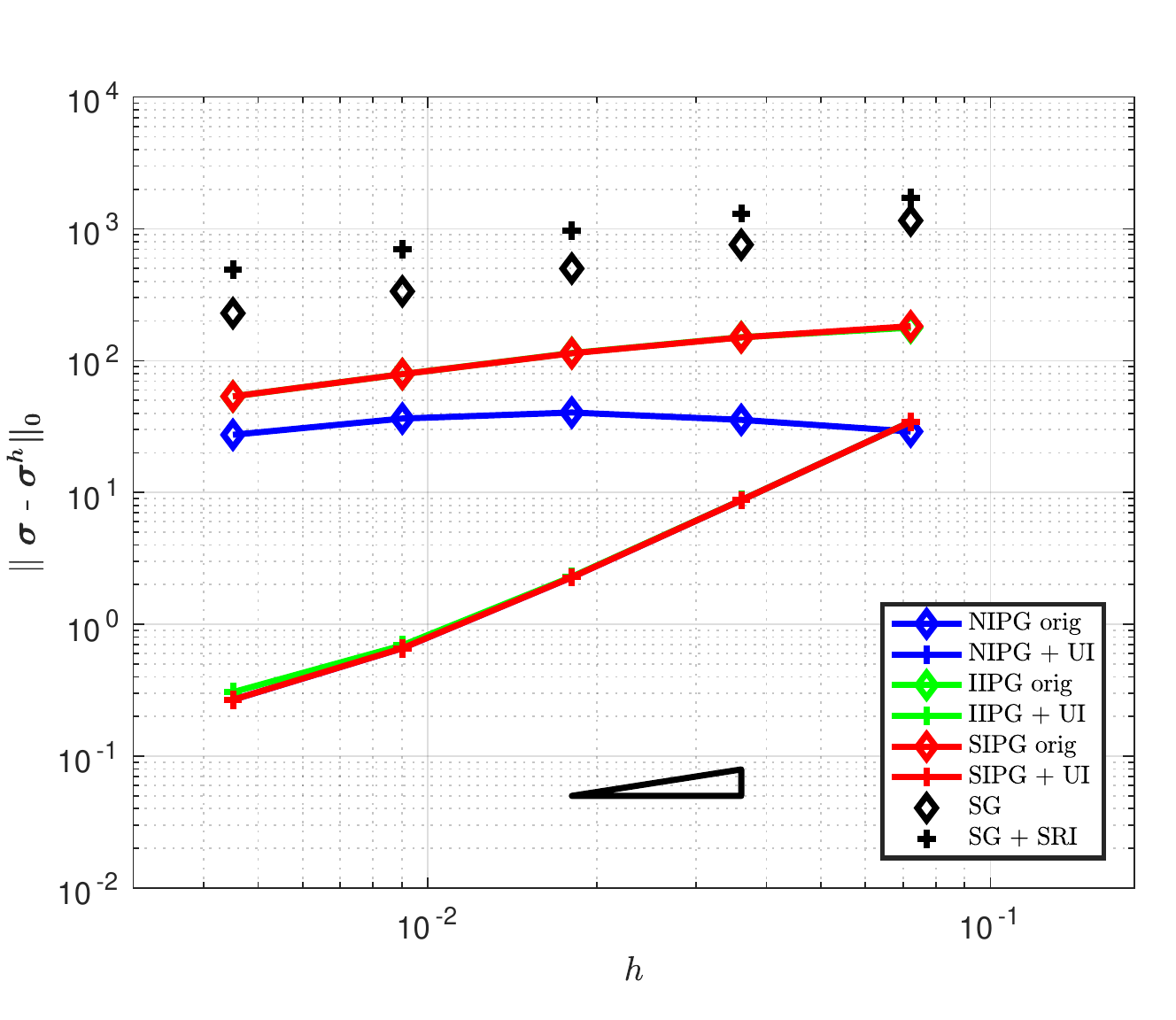}\label{LS_stress_error_various_nu2}}
\caption{L-shaped domain: Stress $L^2$ error convergence. The hypotenuse of the triangle has a slope of $2/3$ in each case.}
\end{figure}


\begin{figure}[!ht]
\centering
\subfloat[]{\includegraphics[width=.30\columnwidth]{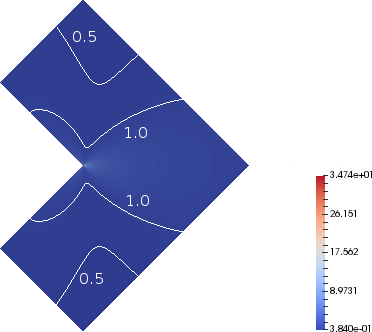}\label{SP_stress_xx_exact_full}} \hspace{10mm}
\subfloat[]{\includegraphics[width=.30\columnwidth]{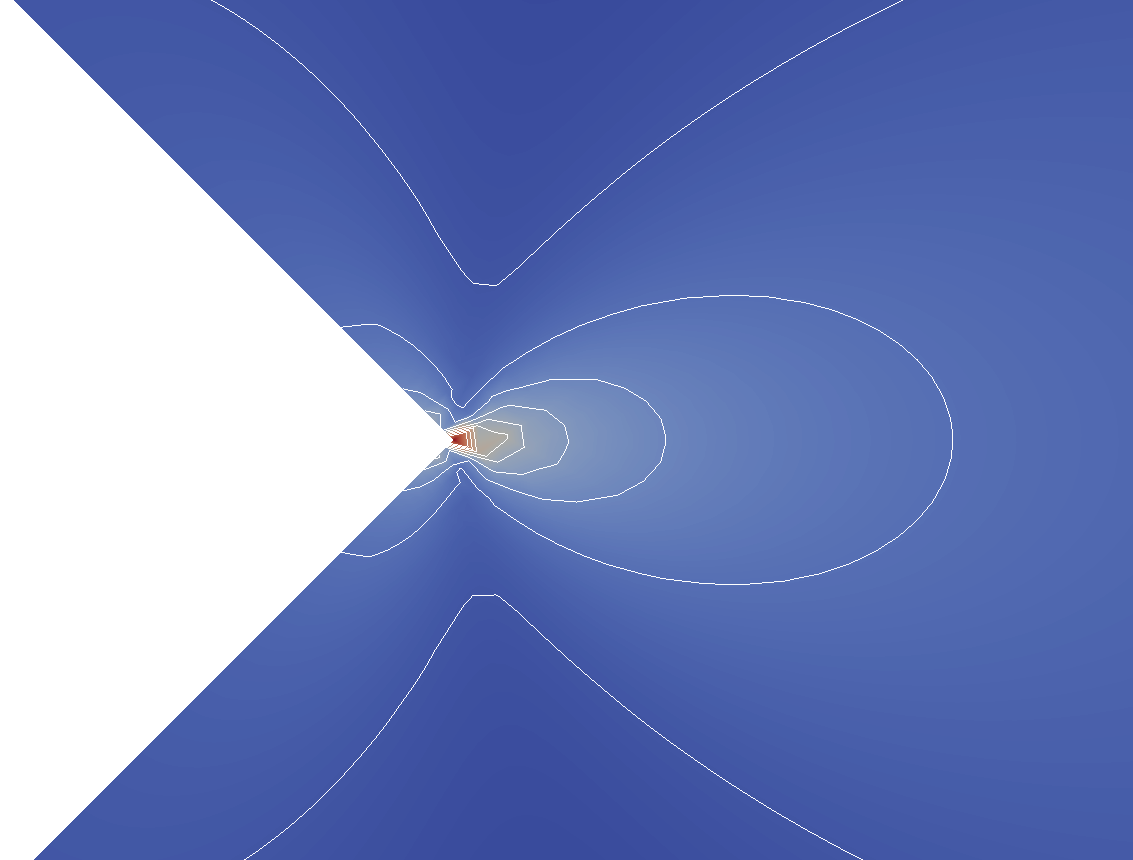}\label{LS_stress_xx_exact_zoom}}
\caption{L-shaped domain: $\sigma_{xx}$, projected exact solution, $\nu = 0.49995$}\label{LS_stress_exact_xx}
\end{figure}

\begin{figure}[!ht]
\centering
\subfloat[Mesh 1]{\includegraphics[width=.20\columnwidth]{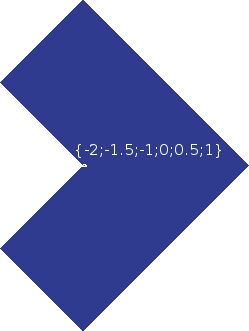}\label{SP_stress_xx_SIPG_full_m1}} \hspace{10mm}
\subfloat[Mesh 1]{\includegraphics[width=.30\columnwidth]{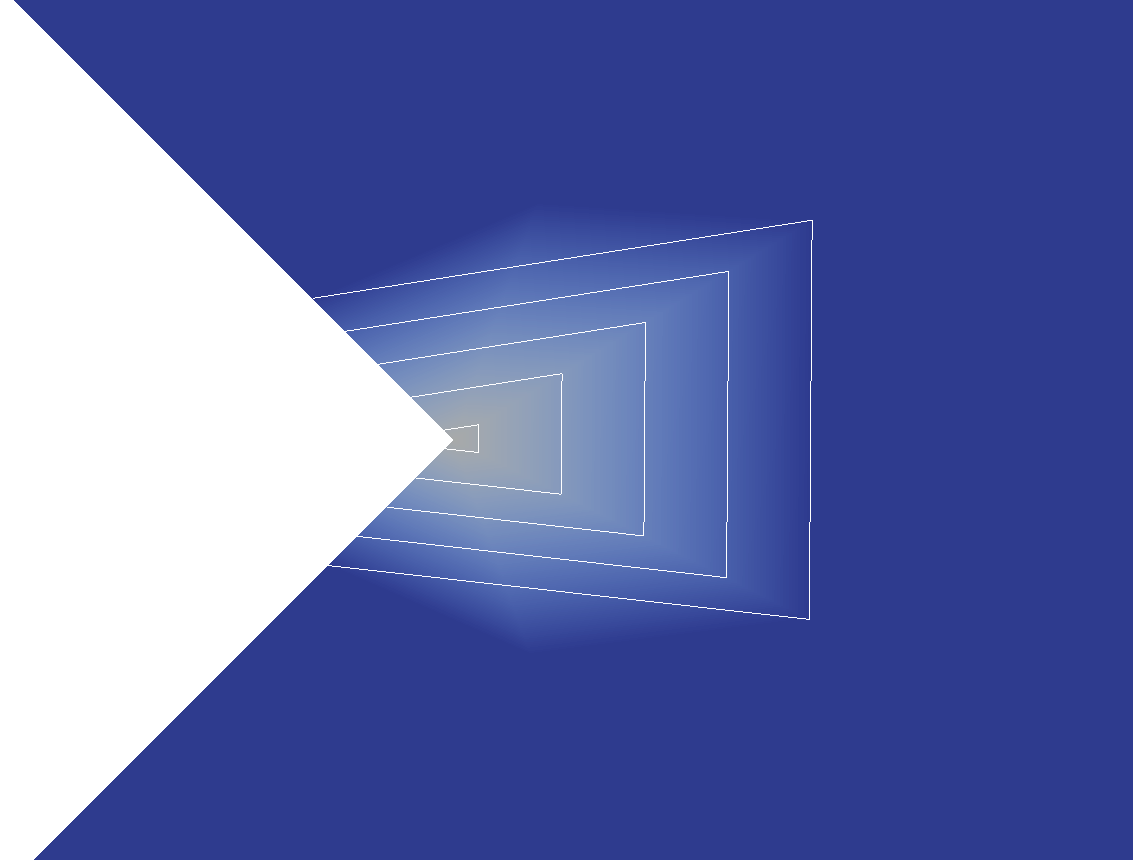}\label{LS_stress_xx_SIPG_zoom_m1}} \\
\subfloat[Mesh 2 (unlabelled contours are at values 0.5 and 1, in sequence towards the singularity)]{\includegraphics[width=.20\columnwidth]{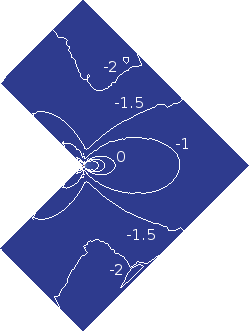}\label{SP_stress_xx_SIPG_full_m2}} \hspace{10mm}
\subfloat[Mesh 2]{\includegraphics[width=.30\columnwidth]{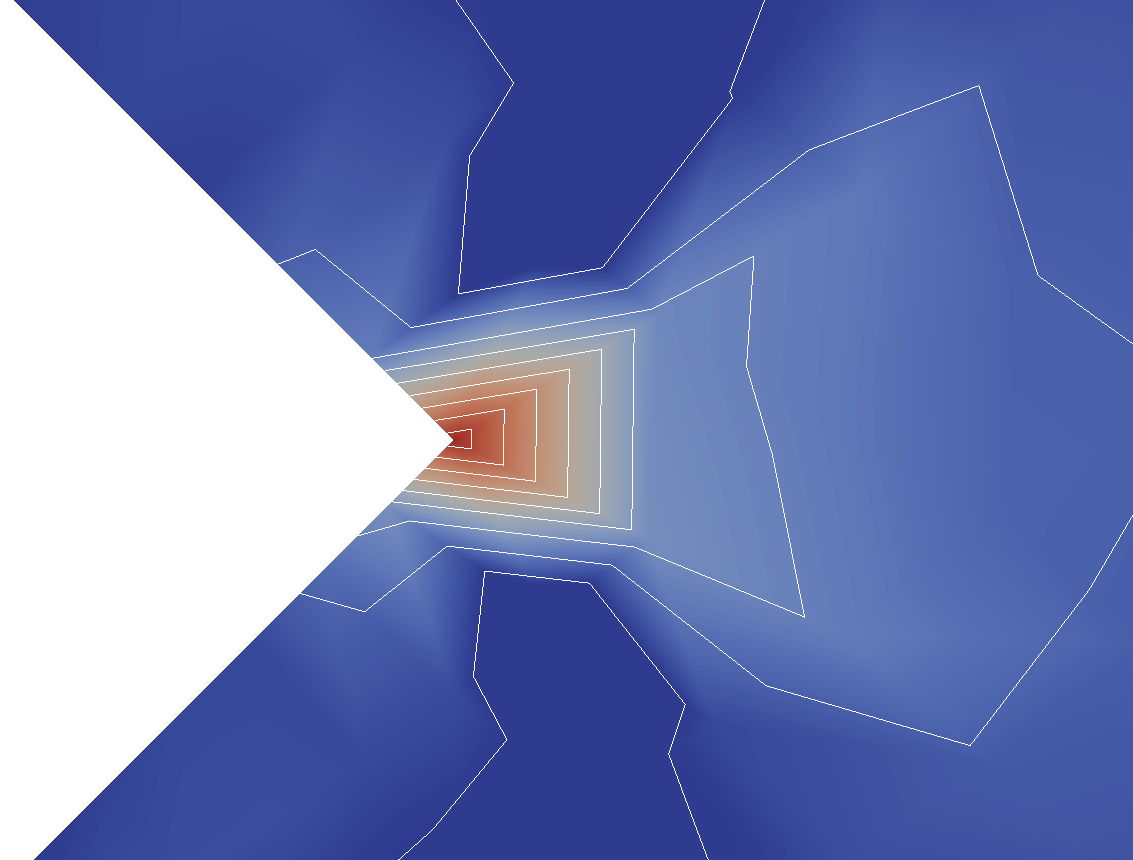}\label{LS_stress_xx_SIPG_zoom_m2}} \\
\subfloat[Mesh 5]{\includegraphics[width=.275\columnwidth]{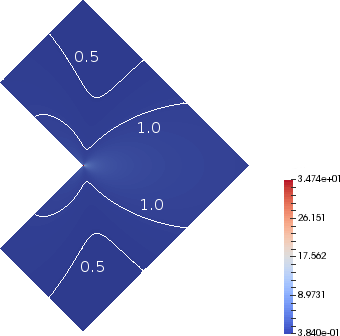}\label{SP_stress_xx_SIPG_full_m5}} \hspace{10mm}
\subfloat[Mesh 5]{\includegraphics[width=.30\columnwidth]{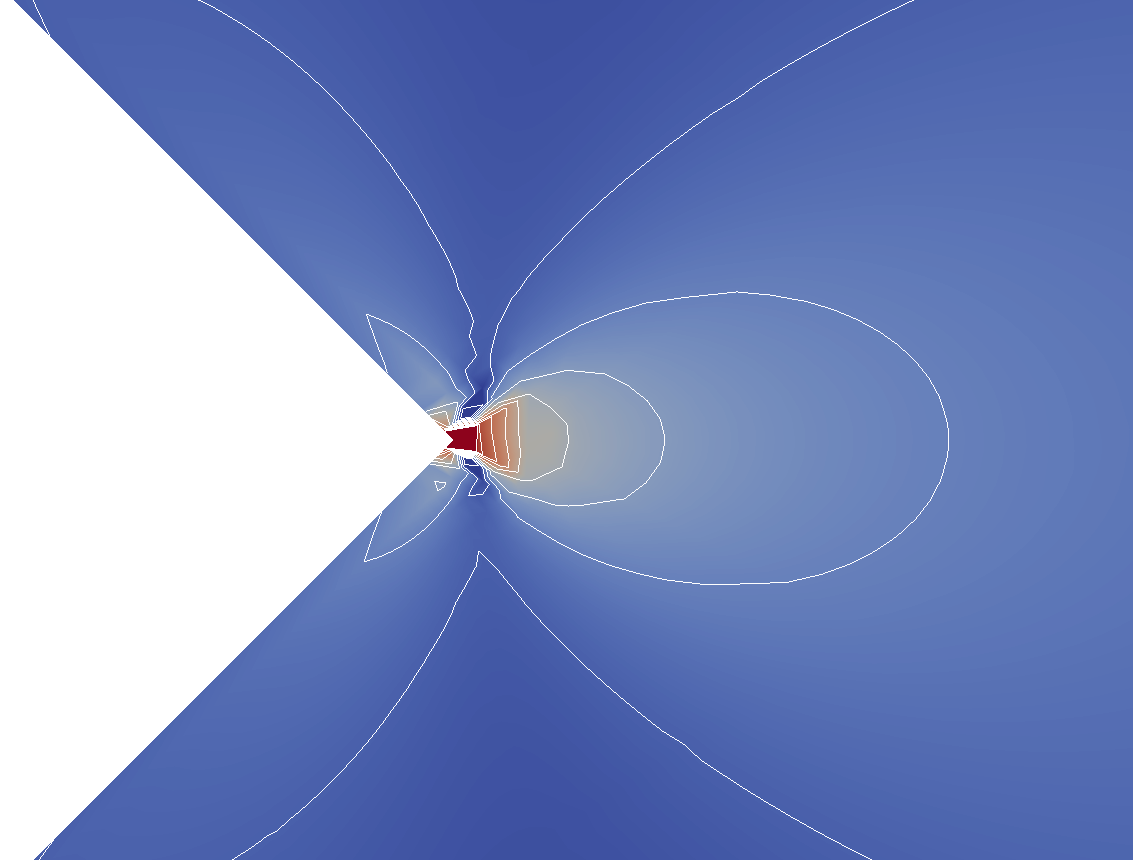}\label{LS_stress_xx_SIPG_zoom_m5}}
\caption{L-shaped domain: $\sigma_{xx}$, SIPG, $\nu = 0.49995$}\label{LS_stress_xx_SIPG}
\end{figure}

\begin{figure}[!ht]
\centering
\subfloat[]{\includegraphics[width=.30\columnwidth]{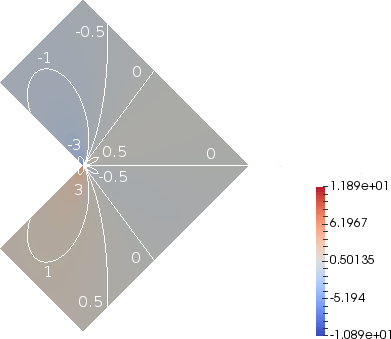}\label{SP_stress_xy_exact_full}} \hspace{10mm}
\subfloat[]{\includegraphics[width=.30\columnwidth]{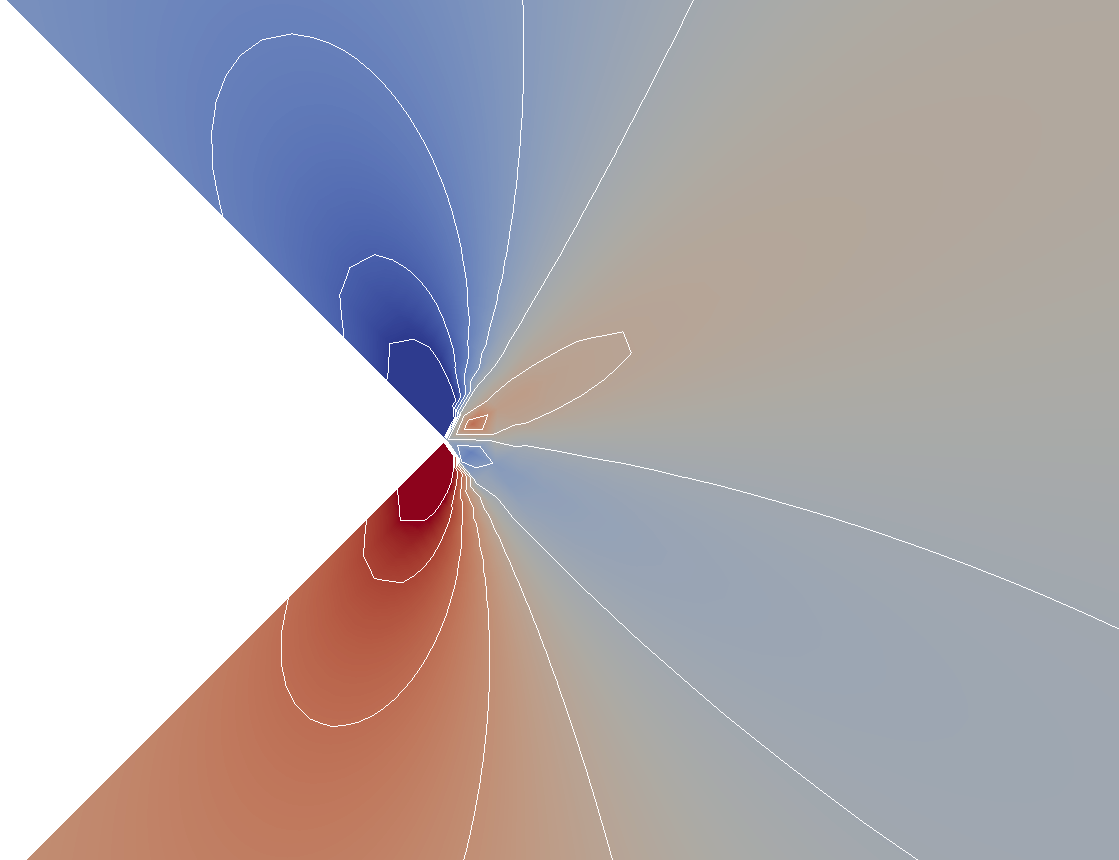}\label{LS_stress_xy_exact_zoom}}
\caption{L-shaped domain: $\sigma_{xy}$, projected exact solution, $\nu = 0.49995$}\label{LS_stress_xy_exact}
\end{figure}

\begin{figure}[!ht]
\centering
\subfloat[Mesh 2]{\includegraphics[width=.20\columnwidth]{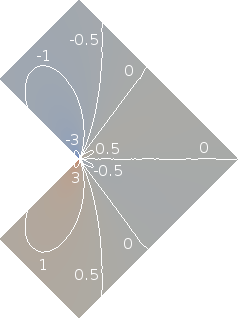}\label{SP_stress_xy_SIPG_full_m2}} \hspace{27mm}
\subfloat[Mesh 2]{\includegraphics[width=.30\columnwidth]{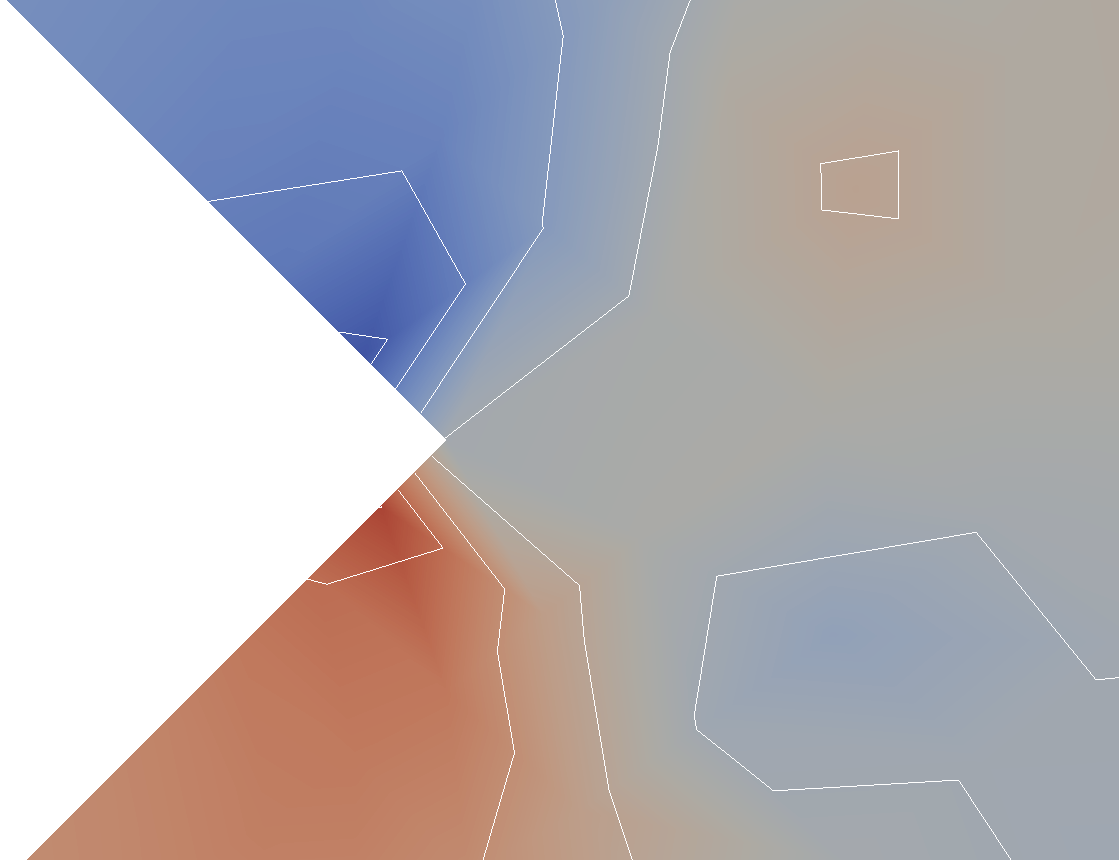}
\label{LS_stress_xy_SIPG_zoom_m2}} \\
\subfloat[Mesh 5]{\includegraphics[width=.30\columnwidth]{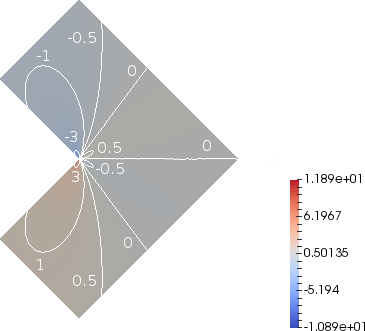}\label{SP_stress_xy_SIPG_full_m5}} \hspace{10mm}
\subfloat[Mesh 5]{\includegraphics[width=.30\columnwidth]{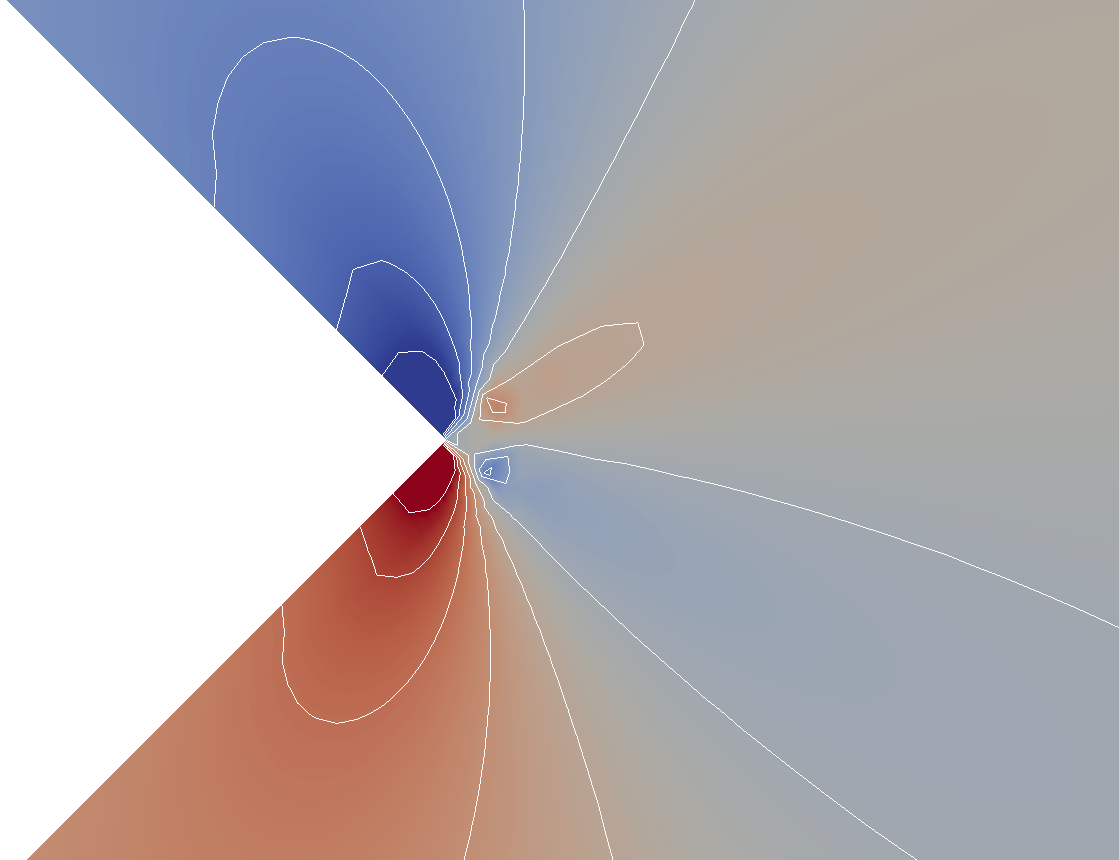}\label{LS_stress_xy_SIPG_zoom_m5}}
\caption{L-shaped domain: $\sigma_{xy}$, SIPG, $\nu = 0.49995$}\label{LS_stress_xy_SIPG}
\end{figure}

\begin{figure}[!ht]
\centering
\subfloat[]{\includegraphics[width=.30\columnwidth]{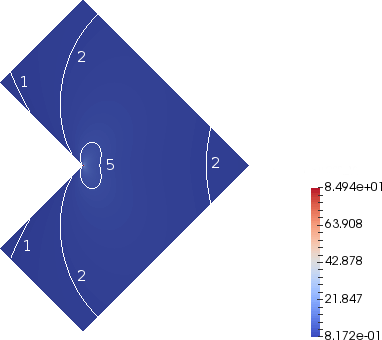}\label{LS_stress_yy_exact_full}} \hspace{10mm}
\subfloat[]{\includegraphics[width=.30\columnwidth]{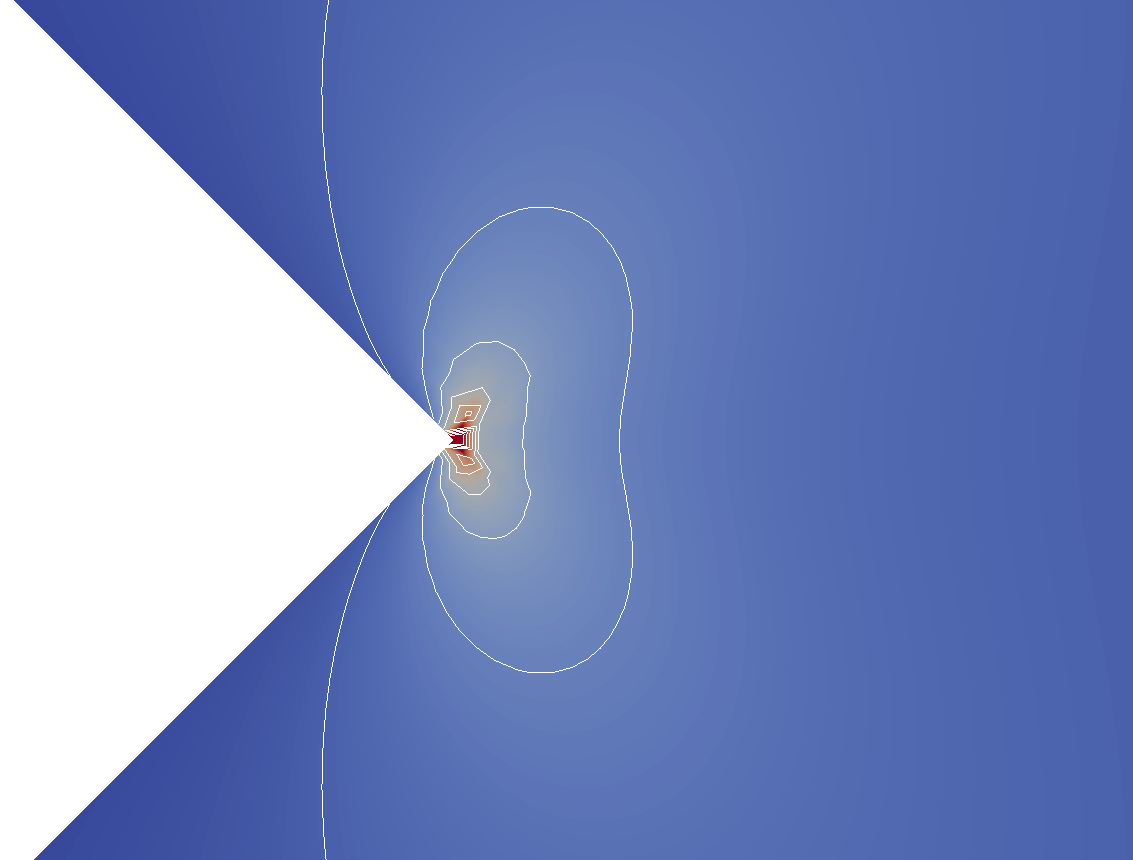}\label{LS_stress_yy_exact_zoom}}
\caption{L-shaped domain: $\sigma_{yy}$, projected exact solution, $\nu = 0.49995$}\label{LS_stress_yy_exact}
\end{figure}

\begin{figure}[!ht]
\centering
\subfloat[Mesh 1 (unlabelled contours have values $\{-3;0;1;2;5\}$, in sequence towards the singularity)]{\includegraphics[width=.18\columnwidth]{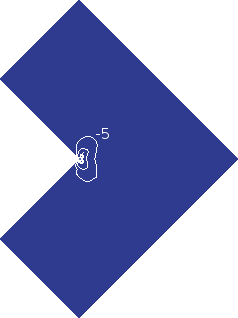}\label{SP_stress_yy_SIPG_full_m1}} \hspace{10mm}
\subfloat[Mesh 1]{\includegraphics[width=.30\columnwidth]{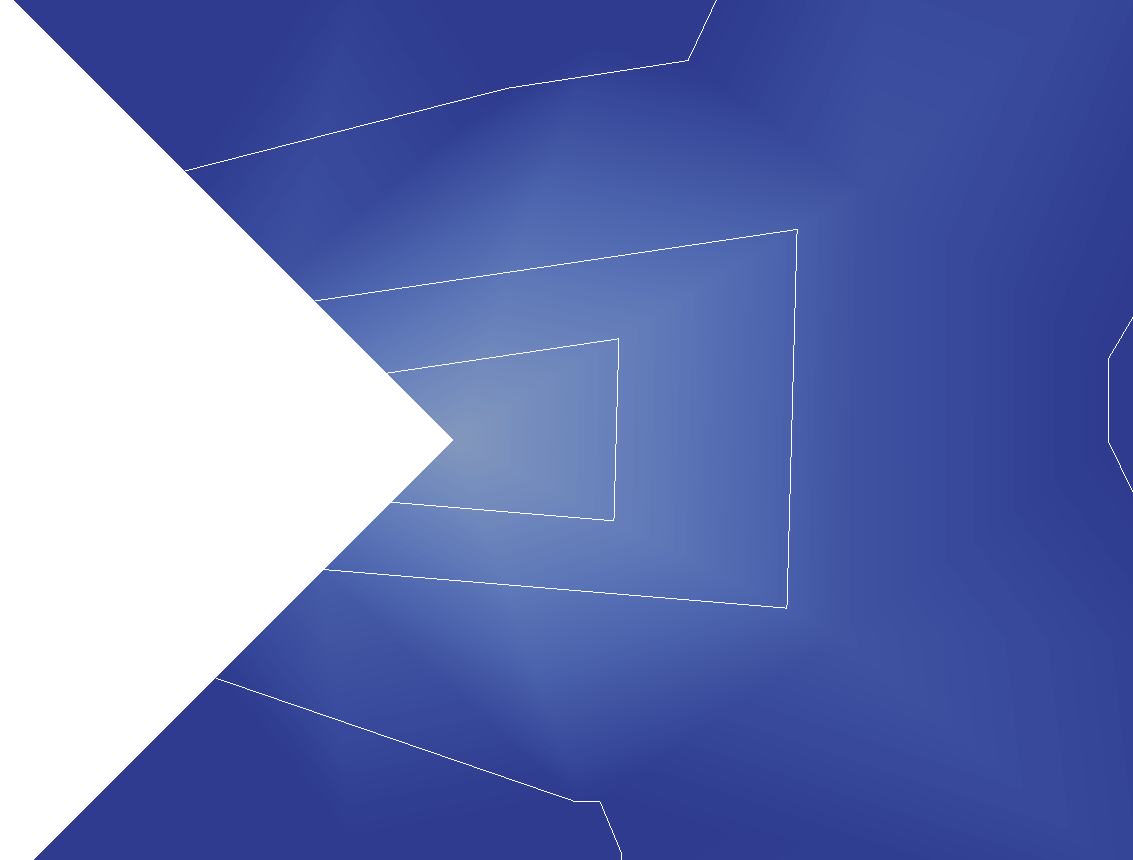}\label{LS_stress_yy_SIPG_zoom_m1}} \\\subfloat[Mesh 2]{\includegraphics[width=.18\columnwidth]{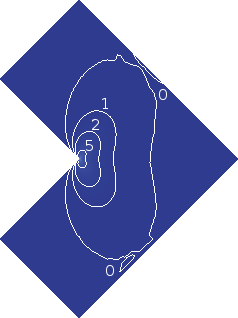}\label{SP_stress_yy_SIPG_full_m2}} \hspace{10mm}
\subfloat[Mesh 2]{\includegraphics[width=.30\columnwidth]{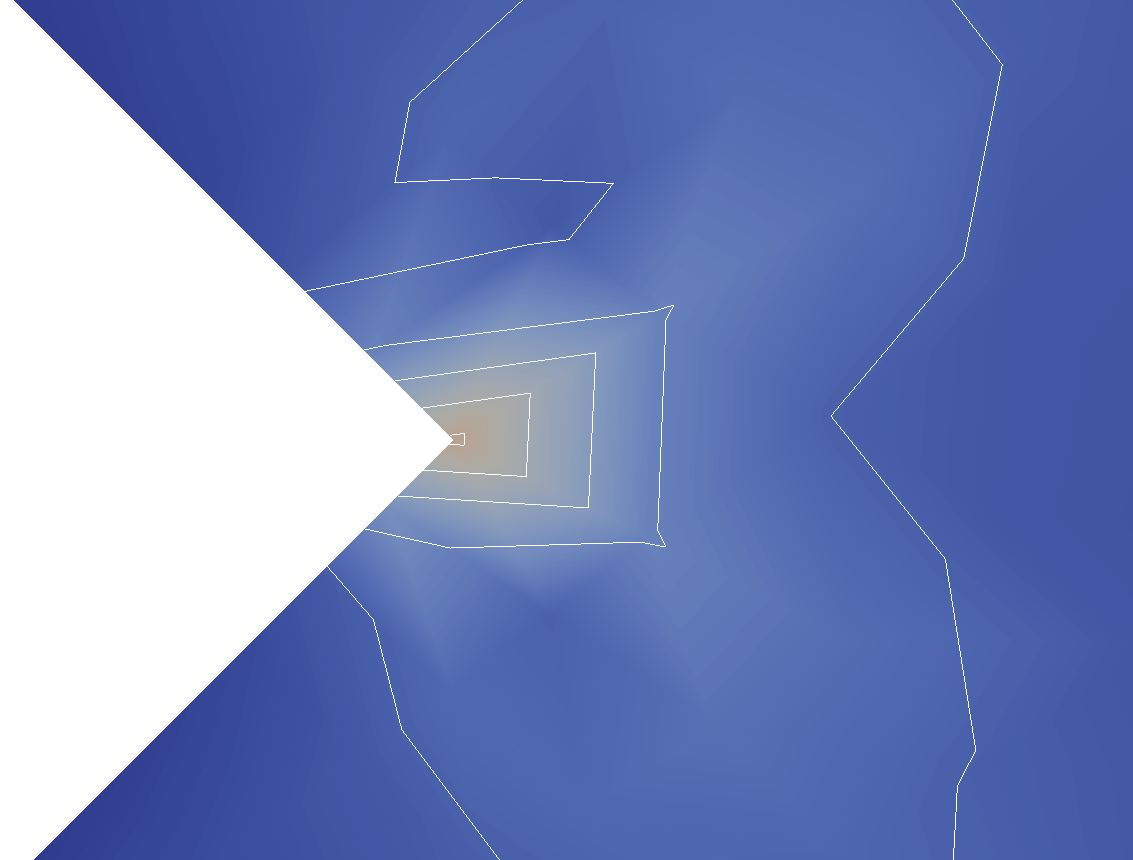}\label{LS_stress_yy_SIPG_zoom_m2}} \\
\subfloat[Mesh 5]{\includegraphics[width=.30\columnwidth]{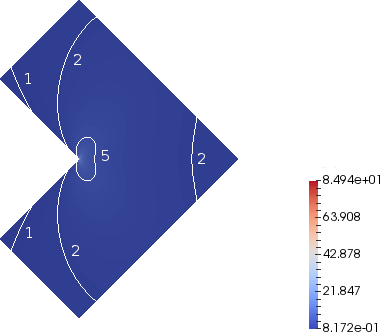}\label{SP_stress_yy_SIPG_full_m5}} \hspace{10mm}
\subfloat[Mesh 5]{\includegraphics[width=.30\columnwidth]{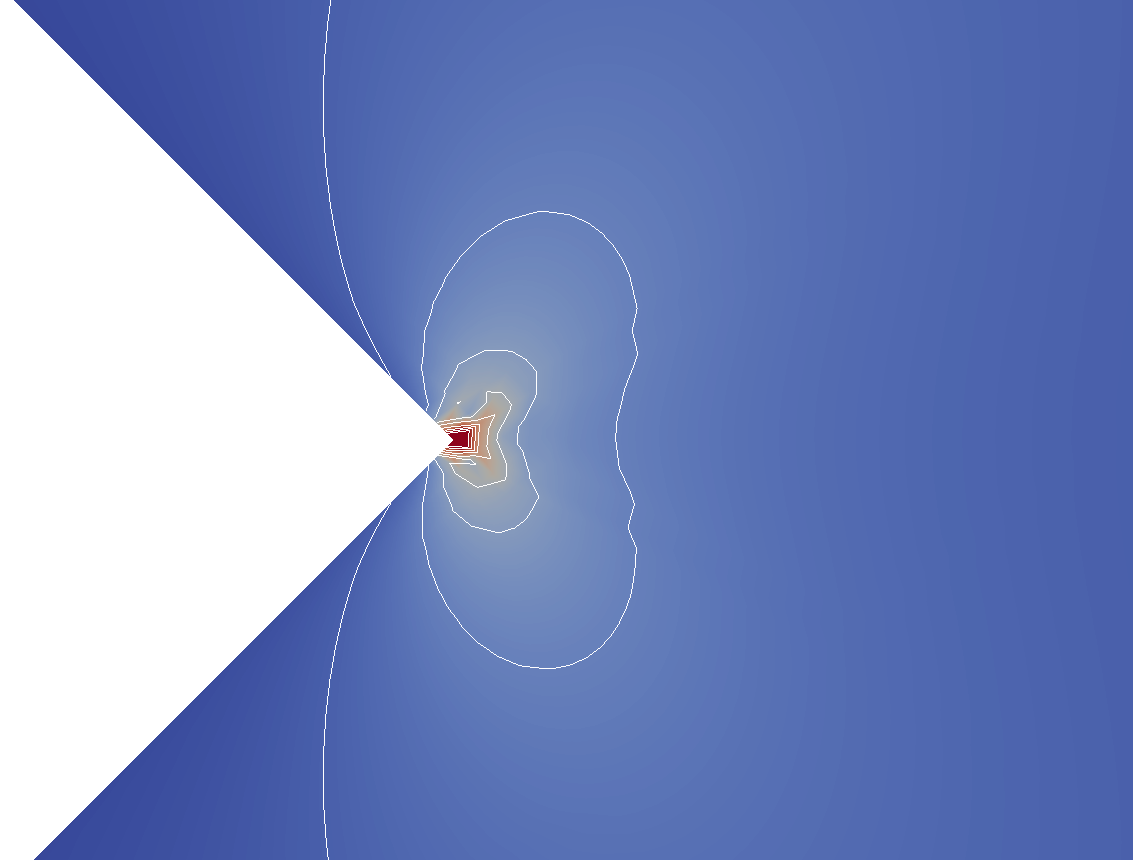}\label{LS_stress_yy_SIPG_zoom_m5}}
\caption{L-shaped domain: $\sigma_{yy}$, SIPG, $\nu = 0.49995$}\label{LS_stress_yy_SIPG}
\end{figure}

With the SG method, convergence for $\nu = 0.3$ is at a rate of about 0.54 for all 3 variants (Figure \ref{LS_stress_error_SG}). For $\nu = 0.49995$, there is no noticeable decrease in convergence rate (with $Q_1$ elements there is a slight increase) but the error is several orders of magnitude larger. With the original IP methods (Figure \ref{LS_stress_error_IPold}) the convergence rate is likewise 0.54 for $\nu = 0.3$, and there is a decrease in rate for NIPG though not for SIPG and IIPG, but in all three methods there is an increase in error magnitude for $\nu = 0.49995$. The new IP methods by contrast have a significantly higher convergence rate for $\nu = 0.49995$ than the rate of 0.54 produced for $\nu = 0.3$ (Figure \ref{LS_stress_error_IPnew}). Figure \ref{LS_stress_error_various_nu2} compares the error convergence and magnitudes of the new methods for near-incompressibility to the original as well as to the SG method with and without SRI, showing better performance in all cases.

In evaluating the qualitative accuracy of the postprocessed stress generated by the IP methods when $\nu = 0.49995$, we consider contour plots both of the full domain (with values of the contours, carefully chosen to show the important range of variation away from the origin, indicated on the plots) and of the region very close to the origin (an area of 0.05 units $\times$ 0.075 units, and with contours of ten values linearly dividing the range of the stress). Both presentations have been chosen to elucidate the nature of inaccuracies in the approximations, while the error plots give averaged representations of the overall accuracy. 

Because there is a stress singularity at the origin, an exact stress solution evaluated at nodal values cannot be plotted. Instead, we calculate the exact stress solution at element quadrature points, and find a smooth stress field by projection to the nodal points. This results in the projected exact solution increasing in accuracy with each refinement, most notably with the value at the origin becoming greater with each refinement. Mesh 5 (see Table \ref{mesh_info_LS}) is used for the contour plots.

The projected exact solution for the $\sigma_{xx}$ component is shown in Figure \ref{LS_stress_exact_xx}, with a high stress value appearing at the origin. Figure \ref{LS_stress_xx_SIPG} shows the corresponding SIPG results for various refinement levels. On mesh 1 the results are inaccurate, whether the full domain or the region of the singularity is considered. Mesh 2 produces an improved approximation across the domain but with lower stresses than the exact solution, and captures the broad features of the exact stress field in the region of the singularity. With the highly refined mesh 5, the full domain is approximated visually accurately, and in the region of the singularity the contour features of the projected exact solution are closely approximated, though small regions of lower stress appear. There are two very small patches of negative stress, above and below the high-stress patch, which do not match the true, all-positive stress field -- these patches decrease in area but increase in magnitude with refinement. 

The projected exact solution of the $\sigma_{xy}$ component is shown in Figure \ref{LS_stress_xy_exact}. The results of the SIPG method are shown in Figure \ref{LS_stress_xy_SIPG}, where the behaviour across the full domain is good with both meshes 2 and 5, but a significant increase in accuracy close to the origin can be seen in the contour features when mesh 5 instead of mesh 2 is used. The maximum and minimum values reached using mesh 5 are nevertheless only under two-thirds the magnitude of the corresponding quantities in the projected exact solution on mesh 5.

The performance of the SIPG method for the $\sigma_{yy}$ component (Figure \ref{LS_stress_yy_SIPG}) is similar to that for $\sigma_{xx}$. Considering the full domain, by comparison to the projected exact solution in Figure \ref{LS_stress_yy_exact_full}, mesh 1 produces an inaccurate field, mesh 2 gives a significantly improved field and mesh 5 a visually nearly-accurate one. In the region of the stress singularity, comparison to the projected exact solution in Figure \ref{LS_stress_yy_exact_zoom} shows that the approximation using mesh 1 is inaccurate, that using mesh 2 gives an improvement, and that using mesh 5 reflects the correct features with high accuracy.

The NIPG and IIPG methods display qualitative behaviour of comparable visual accuracy to the SIPG method in all three components.

The effect of mesh refinement on the postprocessed stress field is marked in this example. However, across the domain broadly, the change in values is small, and the inaccuracies have been highlighted in presentation of the results by the choice of contour values to enhance understanding of the behaviour of the methods. In the immediate region of the origin, the maximum and minimum values of stress are produced, and the highly magnified images show the mediocre performance of the IP methods in attaining these extremes using mesh 1, and the vastly improved performance using more refined meshes.

\subsection{Cube with trigonometric body force}

The linear elastic unit cube $[0,1]^3$, with $E = 1 500 000$, is fixed on all its faces and subjected to an internal body force $\ubf$, where

\begin{align}
f_x &= 0.1 \, \mu \,\pi^2 \lrsb{ \lrb{9 \cos 2 \pi x -5} \lrb{\sin 2 \pi y \sin \pi z - \sin \pi y \sin 2 \pi z} + \frac{3}{1 + \lambda} \sin \pi x \sin \pi y \sin \pi z } \notag\\
&\hspace{5mm}+ 0.1 \,\frac{\mu + \lambda}{1 + \lambda} \,\pi^2 \lrb{\sin \pi x \sin \pi y \sin \pi z - \cos \pi x \cos \pi y \sin \pi z - \cos \pi x \sin \pi y \cos \pi z }, \notag\\
f_y &=  0.1 \,\mu \,\pi^2 \lrsb{ \lrb{9 \cos 2 \pi y -5} \lrb{\sin 2 \pi z \sin \pi x - \sin \pi z \sin 2 \pi x} + \frac{3}{1 + \lambda} \sin \pi x \sin \pi y \sin \pi z } \notag\\
&\hspace{5mm}+ 0.1 \,\frac{\mu + \lambda}{1 + \lambda} \,\pi^2 \lrb{\sin \pi x \sin \pi y \sin \pi z - \cos \pi x \cos \pi y \sin \pi z - \sin \pi x \cos \pi y \cos \pi z }, \notag
\end{align}
\begin{align}
f_z &=  0.1 \,\mu \,\pi^2 \lrsb{ \lrb{9 \cos 2 \pi z -5} \lrb{\sin 2 \pi x \sin \pi y - \sin \pi x \sin 2 \pi y} + \frac{3}{1 + \lambda} \sin \pi x \sin \pi y \sin \pi z } \notag\\
&\hspace{5mm}+ 0.1 \,\frac{\mu + \lambda}{1 + \lambda} \,\pi^2 \lrb{\sin \pi x \sin \pi y \sin \pi z - \cos \pi x \sin \pi y \cos \pi z - \sin \pi x \cos \pi y \cos \pi z }.  \notag
\end{align}

The exact solution is given in \cite{Grieshaber2015}.

In this problem, with the SIPG method the stabilization parameter values $k_{\mu} = 50$ and $k_{\lambda} = 250$ are used. These values are introduced here for their effect on the stability of the displacement results for the original and new SIPG methods.

Table \ref{mesh_info_Cube} details the number of elements and dofs for each method at the various refinement levels.

Displacement results of the original and new IP methods with cubic elements (i.e.\ with \textit{df} = 0.0) for the first four meshes have been presented in \cite{Grieshaber2015} and are repeated here for comparison.

\begin{table}[H]
\begin{center}
\begin{tabular}{|c|c|c|c|c|c|c|} \hline
$n$ & no.\ els & \multicolumn{2}{|c|}{no.\ dofs} \\ \hline 
& & SG $Q_1$ &  IP \\ \hline
1 & 8 & 81  & 192\\ \hline
2 & 64 & 375 & 1536\\ \hline
3 & 512 & 2187 & 12288\\ \hline
4 & 4096 & 14739  & 98304\\ \hline
5 & 32768 & 107811 & 786432\\ \hline
\end{tabular}
\end{center}
\caption{Mesh details for standard Galerkin (SG) and interior penalty (IP) methods for the cube example: mesh $n$ has $2^n$ elements per edge; number of elements and number of degrees of freedom are shown.}\label{mesh_info_Cube}\end{table}

\subsubsection{Displacement approximation}


\begin{figure}[!ht]
\centering
\subfloat[SG $Q_1$]{\includegraphics[trim={3.5cm 9cm 4cm 9cm},clip,width=.40\columnwidth]{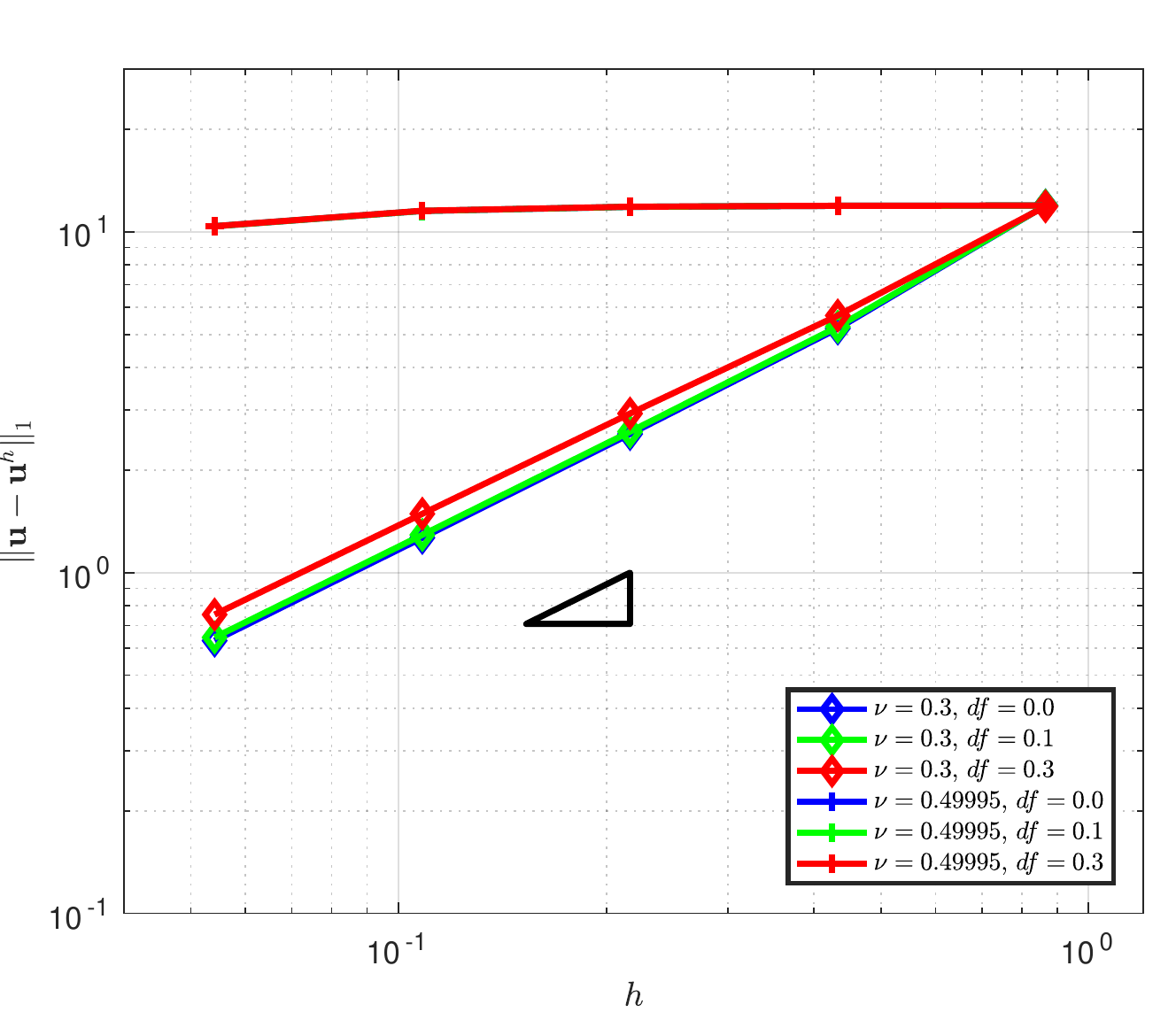}\label{Cube_disp_error_SG}}
\subfloat[SG $Q_1$ with SRI]{\includegraphics[trim={3.5cm 9cm 4cm 9cm},clip,width=.40\columnwidth]{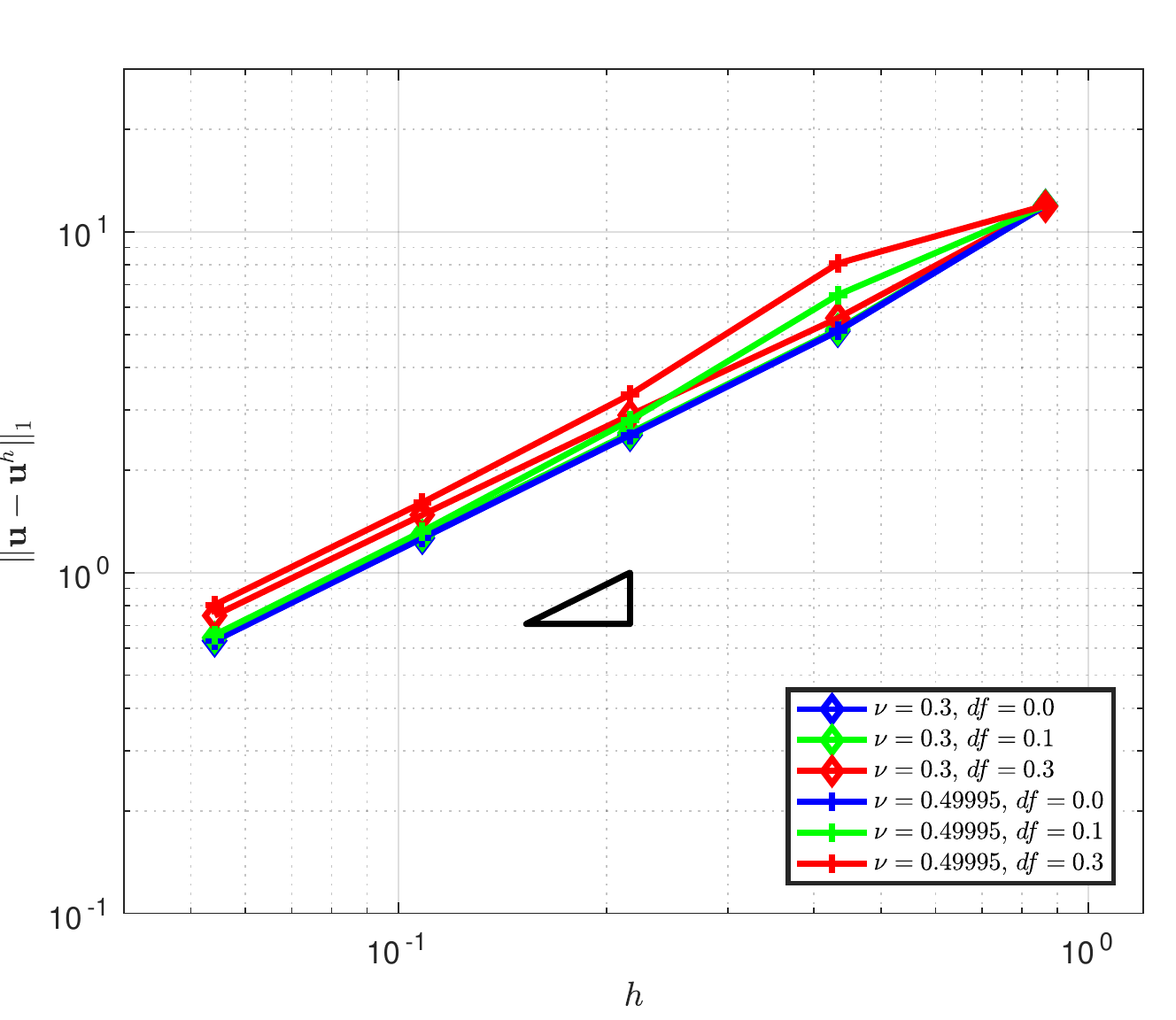}\label{Cube_disp_error_SG_sri}}\\
\subfloat[Original IP]{\includegraphics[trim={3.5cm 9cm 4cm 9cm},clip,width=.40\columnwidth]{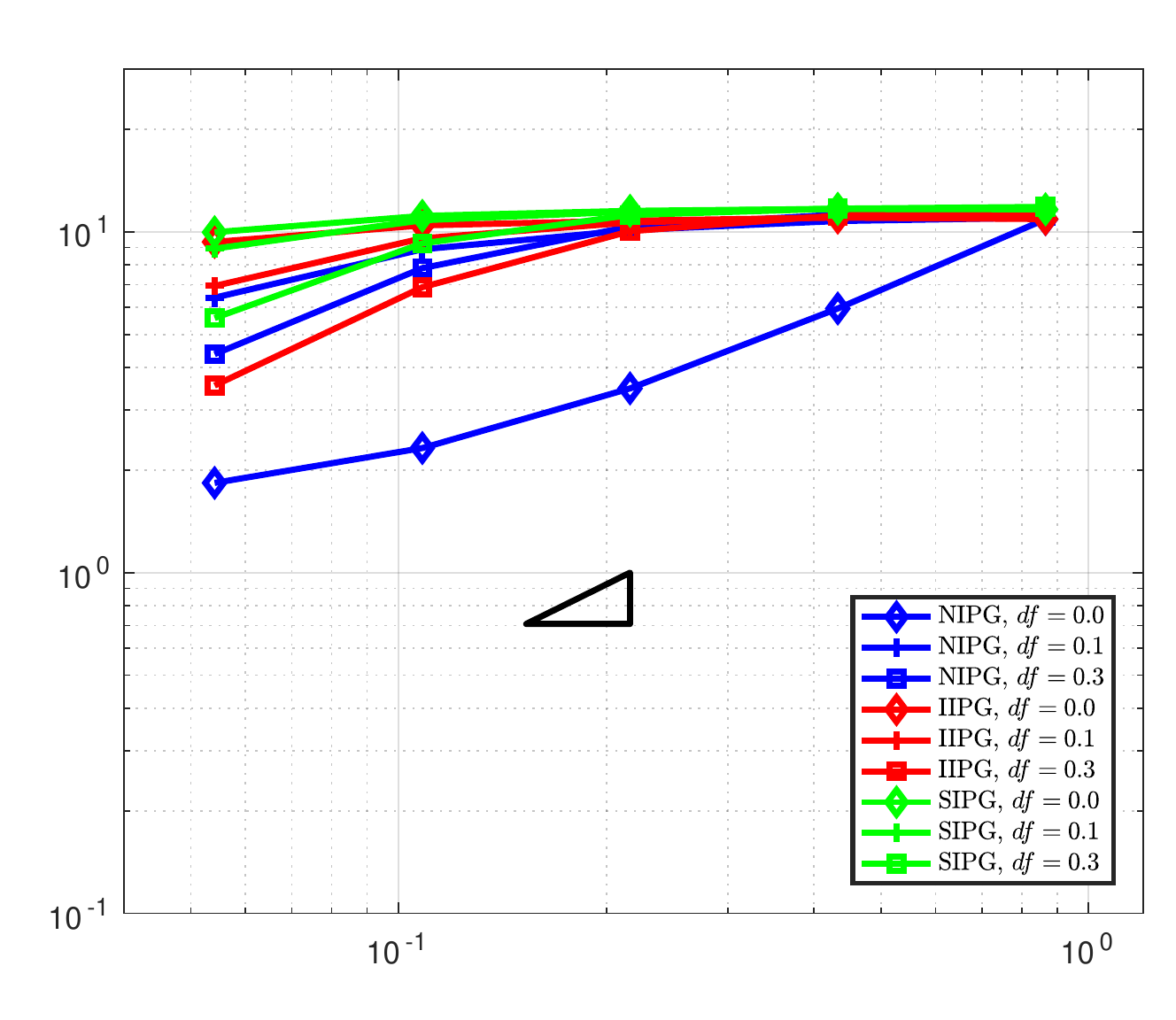}\label{Cube_disp_error_IPorig_nu2}}
\subfloat[New NIPG]{\includegraphics[trim={3.5cm 9cm 4cm 9cm},clip,width=.40\columnwidth]{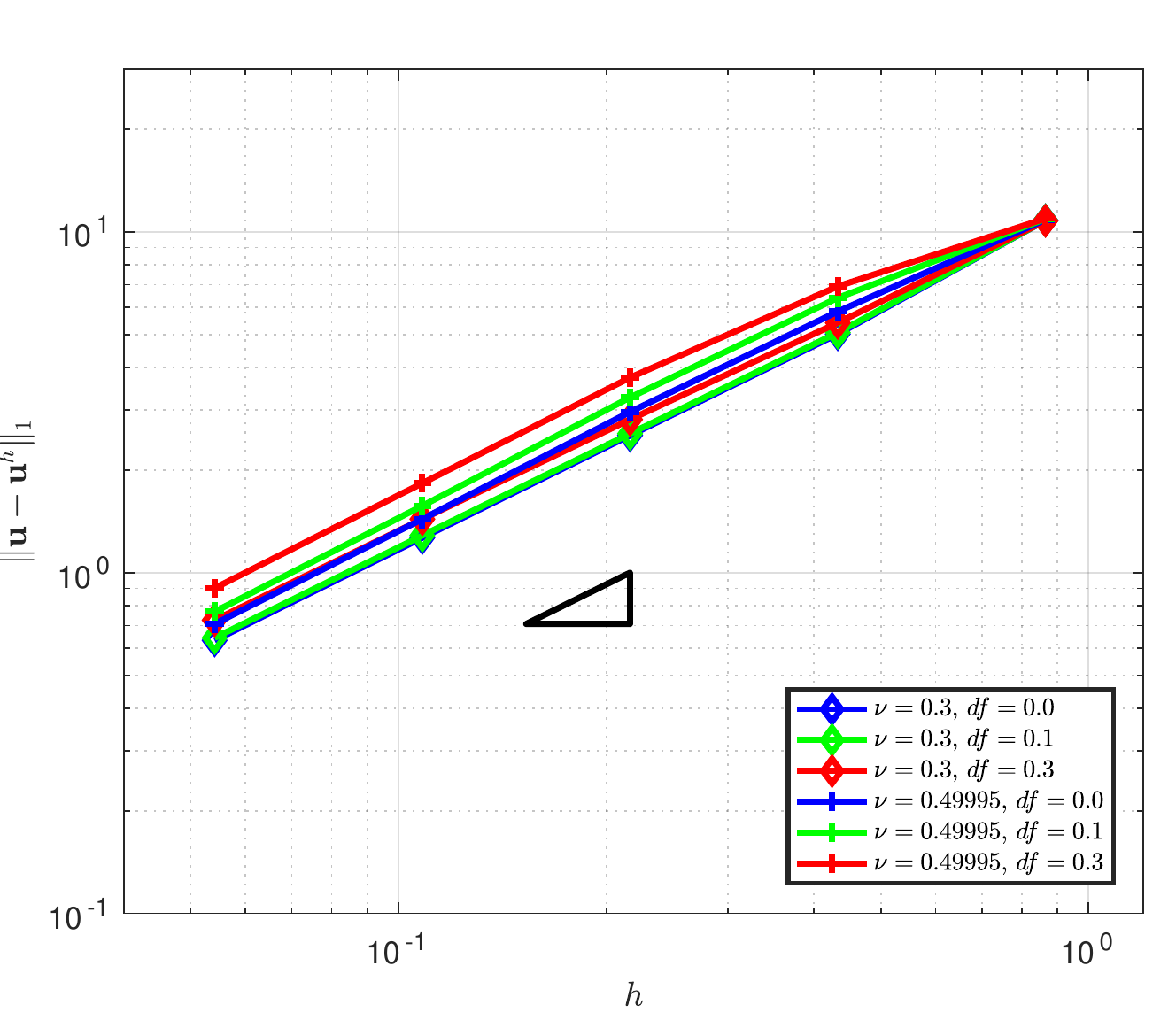}
\label{Cube_disp_error_NIPG}}\\
\subfloat[New IIPG]{\includegraphics[trim={3.5cm 9cm 4cm 9cm},clip,width=.40\columnwidth]{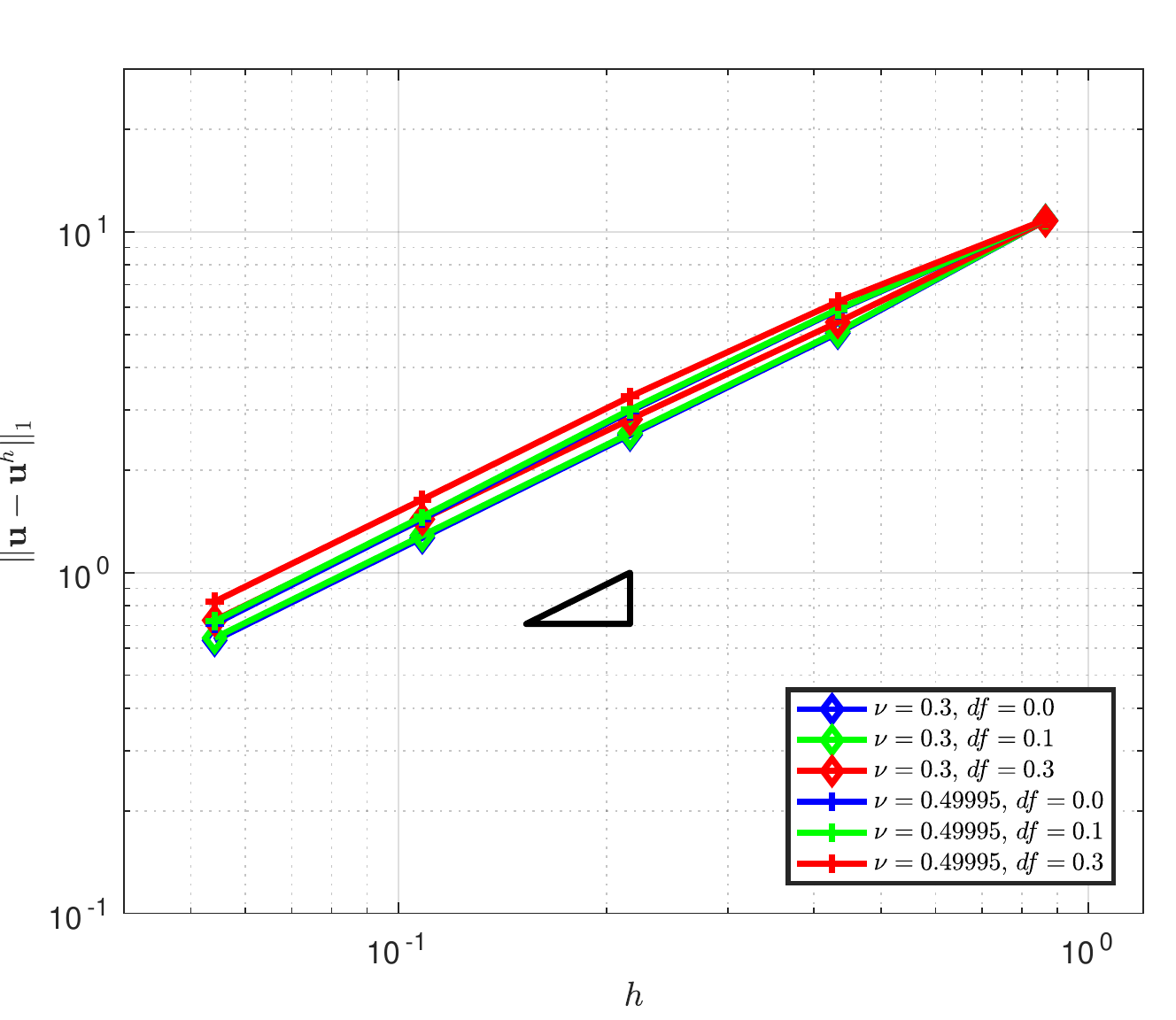}\label{Cube_disp_error_IIPG}}
\subfloat[New SIPG]{\includegraphics[trim={3.5cm 9cm 4cm 9cm},clip,width=.40\columnwidth]{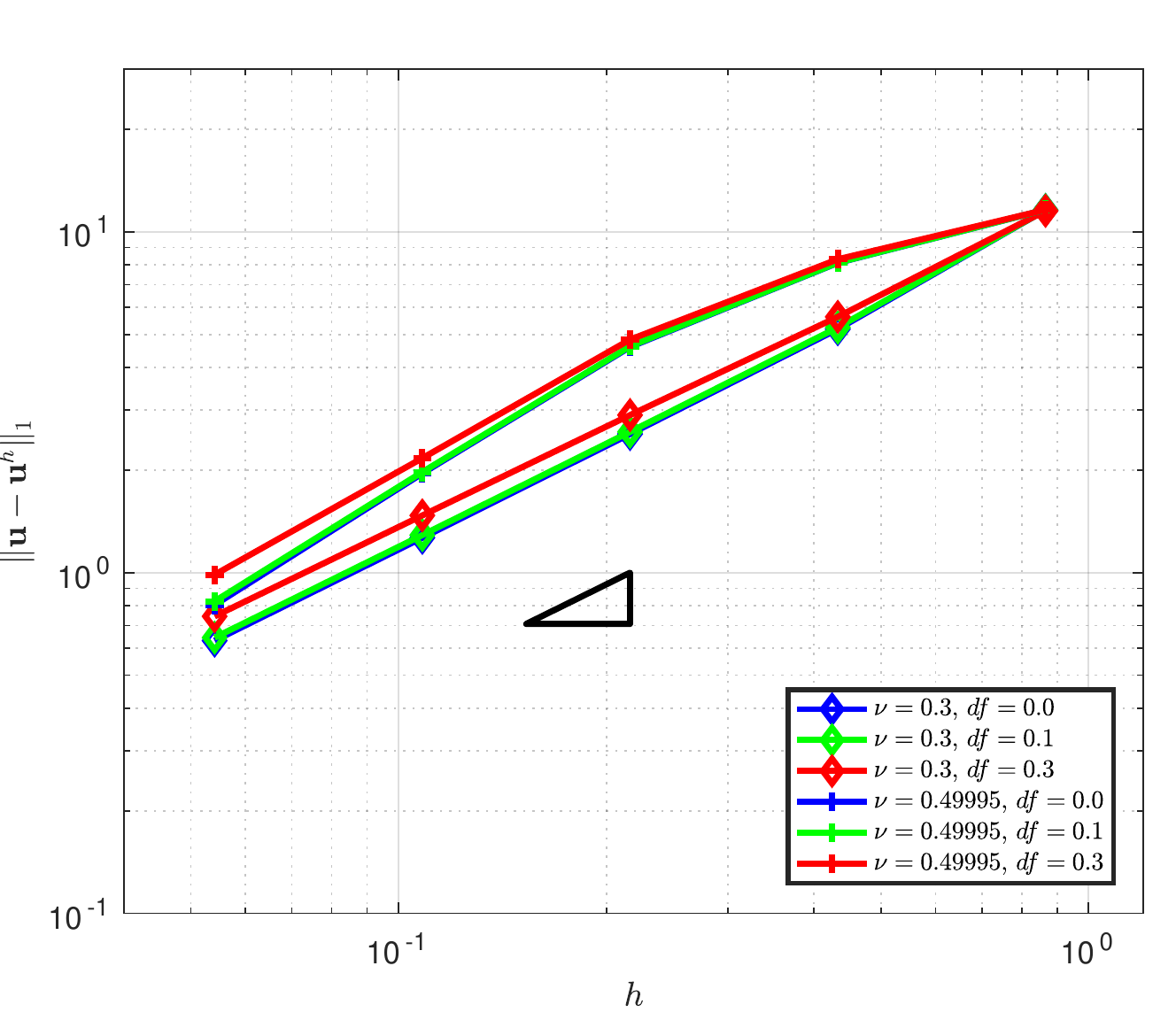}
\label{Cube_disp_error_SIPG_50_250}}
\caption{Cube: Displacement $H^1$ error convergence. The hypotenuse of the triangle has a slope of 1 in each case.}
\end{figure}

\begin{figure}[!ht]
\centering
\includegraphics[trim={3.5cm 9cm 4cm 9cm},clip,width=.49\columnwidth]{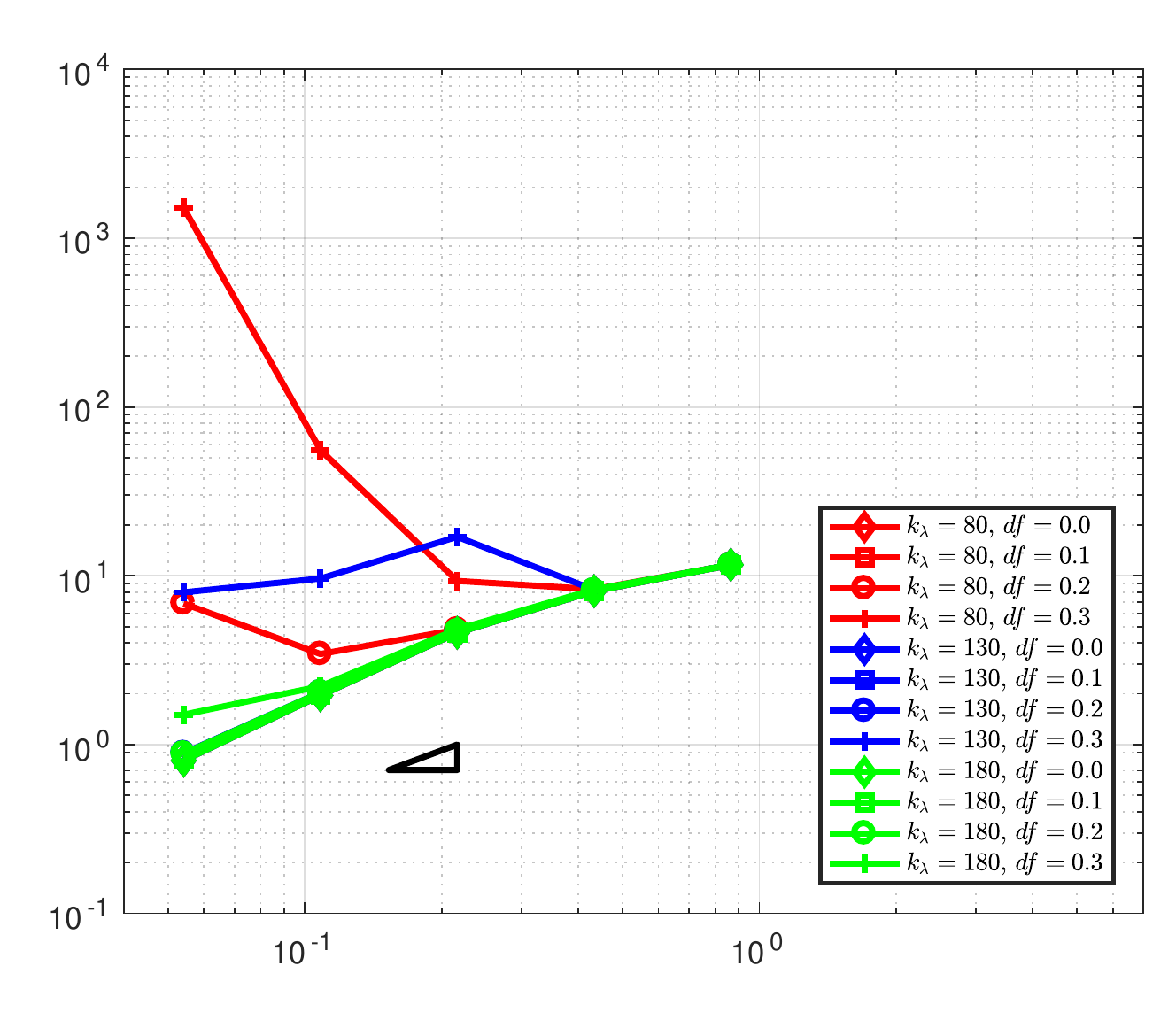}
\caption{Cube: SIPG displacement $H^1$ error convergence with $k_{\mu} = 50$ and varying values of $k_{\lambda}$}\label{Cube_disp_error_SIPG_var_stab}
\end{figure}

\begin{figure}[!ht]
\centering
\subfloat[Exact solution]{\includegraphics[width=.22\columnwidth]{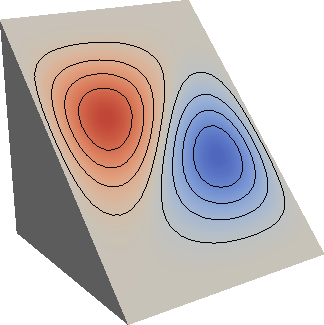}\label{Cube_disp_x_exact}} \hspace{10mm}
\subfloat[IIPG, \textit{df} = 0.0, mesh 5]{\includegraphics[width=.22\columnwidth]{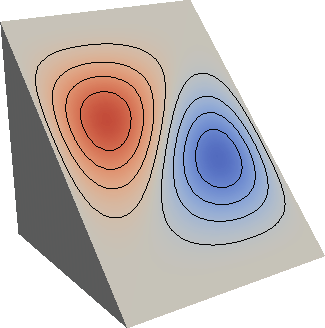}\label{Cube_disp_x_IIPG_00}} \hspace{10mm}
\subfloat[IIPG, \textit{df} = 0.3, mesh 5]{\includegraphics[width=.3\columnwidth]{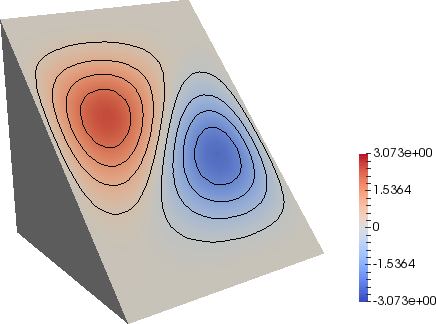}\label{Cube_disp_x_IIPG_03}}
\caption{Cube: Displacement $u_x$, $\nu = 0.49995$}\label{Cube_disp_x}
\end{figure}

The SG method with $Q_1$ elements shows the same behaviour in displacement $H^1$ error convergence here as for the first two examples, as does the same method with SRI (Figures \ref{Cube_disp_error_SG} and \ref{Cube_disp_error_SG_sri}): without SRI, convergence is optimal when $\nu = 0.3$ and poor when $\nu = 0.49995$, while with SRI convergence is optimal throughout. The original IP methods (Figure \ref{Cube_disp_error_IPorig_nu2}) display varied performance for near-incompressibility, depending on method and element regularity, with convergence rates improving significantly with refinement for higher \textit{df} but poor for low \textit{df}. NIPG when \textit{df} = 0.0 has an initially high rate, which deteriorates slightly in the final refinement here. The new IP methods (Figures \ref{Cube_disp_error_NIPG}, \ref{Cube_disp_error_IIPG}, and \ref{Cube_disp_error_SIPG_50_250}) show optimal convergence uniformly with respect to the compressibility parameter both for cubic elements (\textit{df} = 0.0) and for general hexahedra (\textit{df} $>$ 0.0), with SIPG showing a slightly lower convergence rate at the coarse refinement levels, for the near-incompressible case.

Stabilization parameters higher than the usual choices were required for stability of the SIPG method with $\nu = 0.49995$, in this example. While the question of how to choose stabilization parameters is one that is relevant to IP methods broadly and is an area of study in its own right, we note in the context of this work that increasingly high values of $k_{\lambda}$ are required for stability here as element shape regularity decreases, i.e.\ as \textit{df} increases, as demonstrated in Figure \ref{Cube_disp_error_SIPG_var_stab}.

A contour plot of an oblique cross-section of the cube shows by comparison with the exact solution that the accuracy of the $x$-displacement approximation of the SIPG method with cubic elements (\textit{df} = 0.0) is maintained for a mesh of distorted elements (with \textit{df} = 0.3) when $\nu = 0.49995$ (Figure \ref{Cube_disp_x}). The accuracy of the $y$-displacement approximation is similar, and these results extend to the other two IP methods.

\subsubsection{Postprocessed stress}


\begin{figure}[!ht]
\centering
\subfloat[SG $Q_1$]{\includegraphics[trim={3.5cm 9cm 4cm 9cm},clip,width=.40\columnwidth]{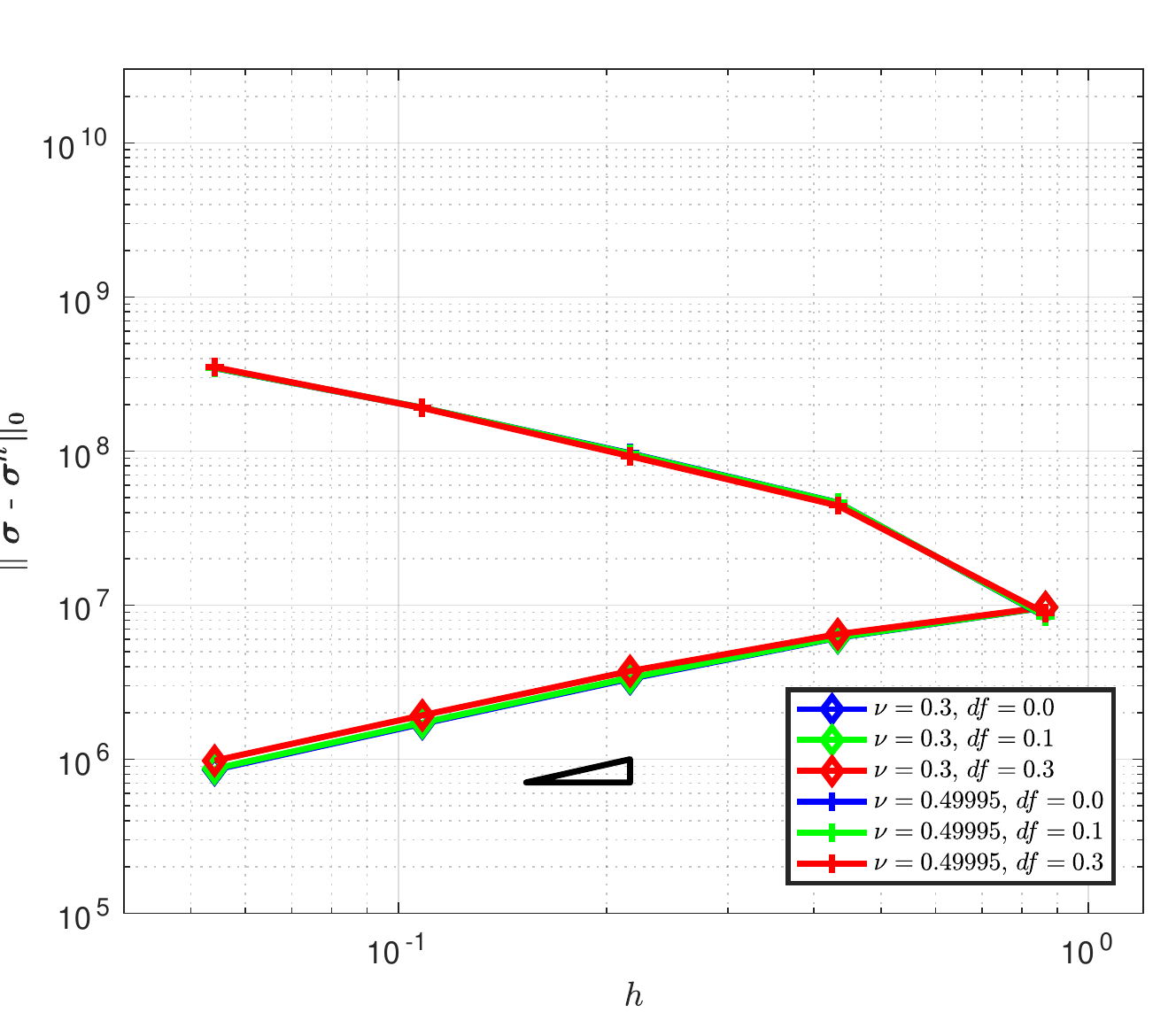}
\label{Cube_stress_error_SG}}
\subfloat[SG $Q_1$ with SRI]{\includegraphics[trim={3.5cm 9cm 4cm 9cm},clip,width=.40\columnwidth]{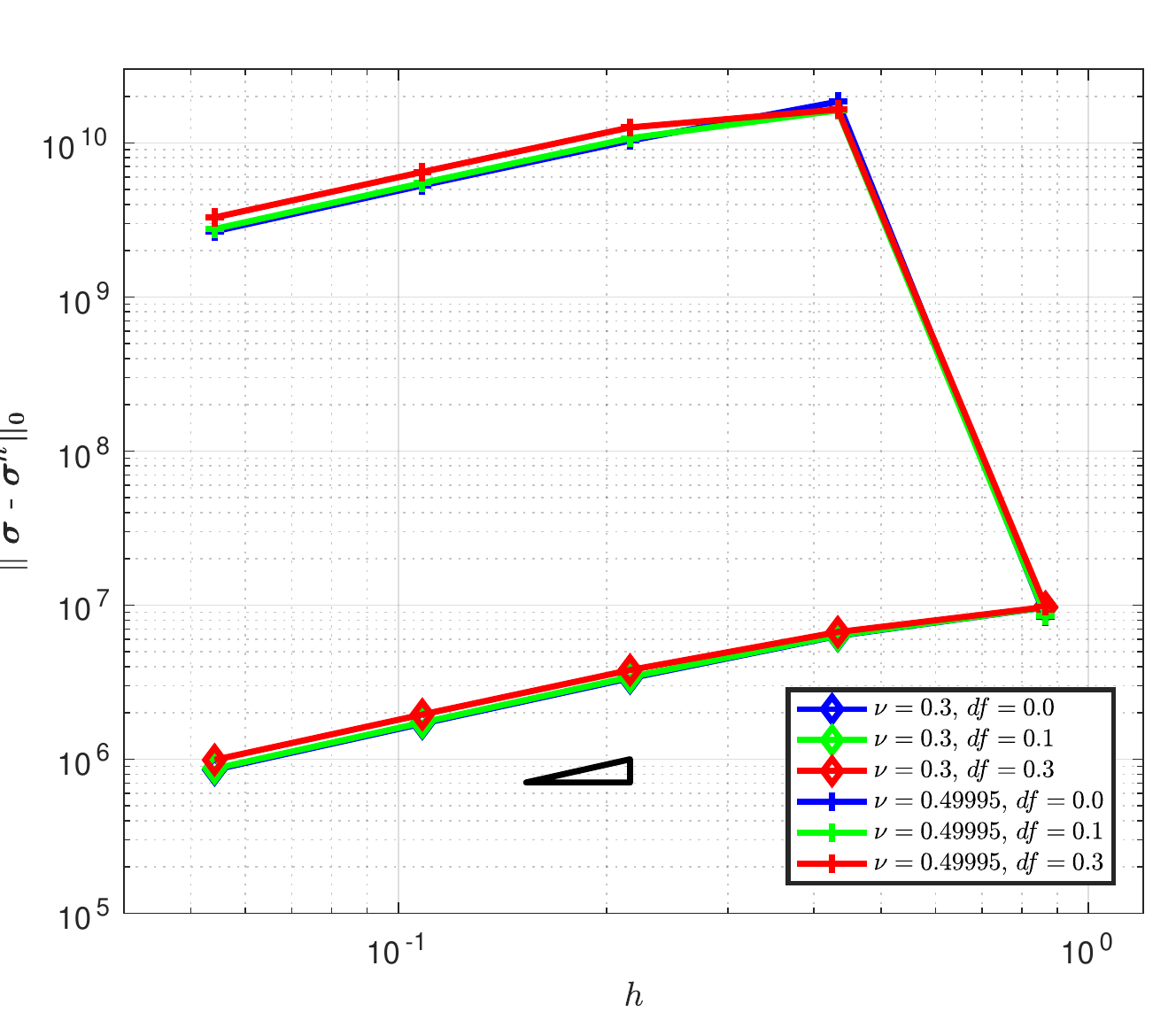}\label{Cube_stress_error_SG_sri}}\\
\subfloat[Original IP, $\nu = 0.49995$]{\includegraphics[trim={3.5cm 9cm 4cm 9cm},clip,width=.40\columnwidth]{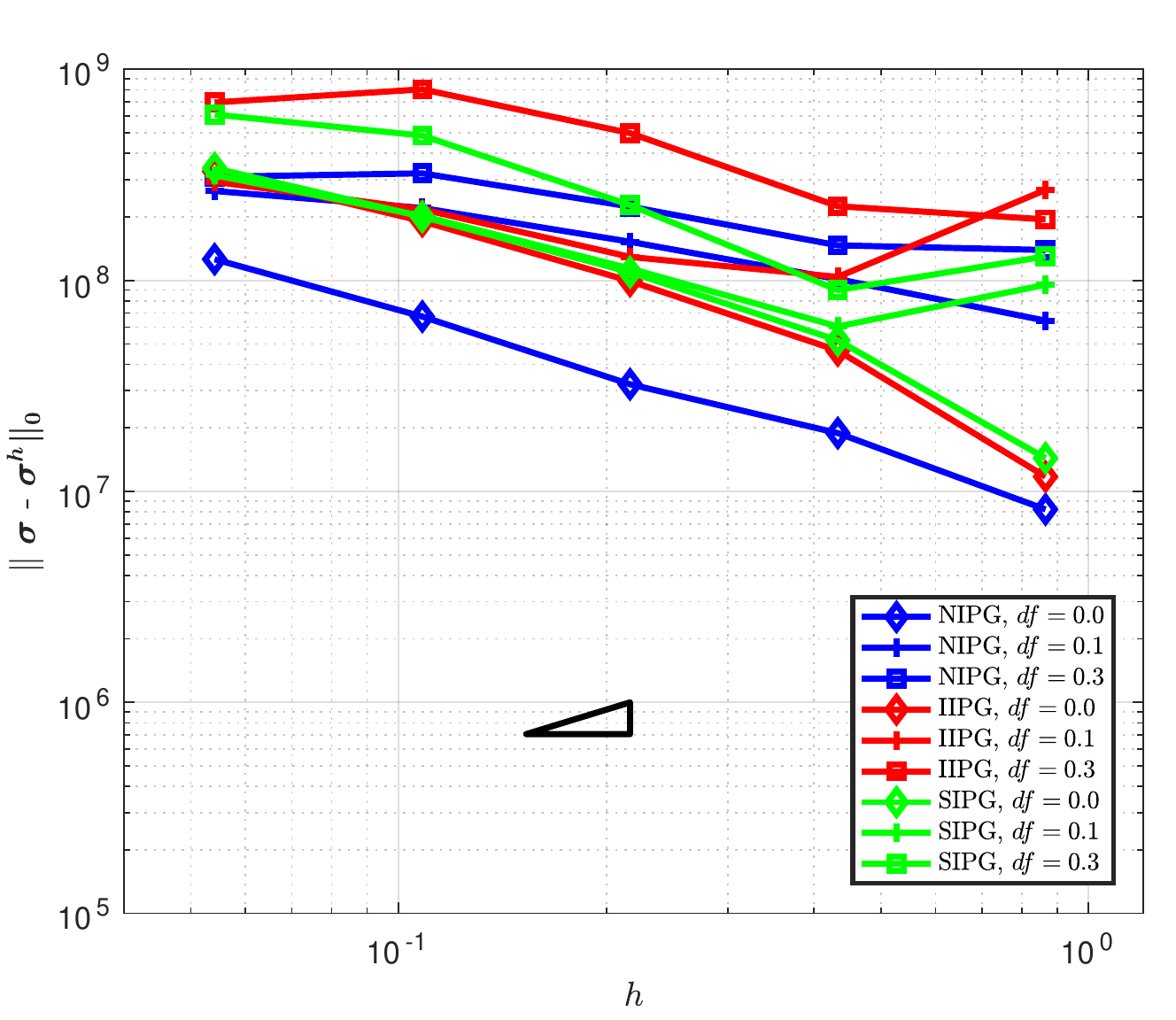}\label{Cube_stress_error_IPorig_nu2}}
\subfloat[New NIPG]{\includegraphics[trim={3.5cm 9cm 4cm 9cm},clip,width=.40\columnwidth]{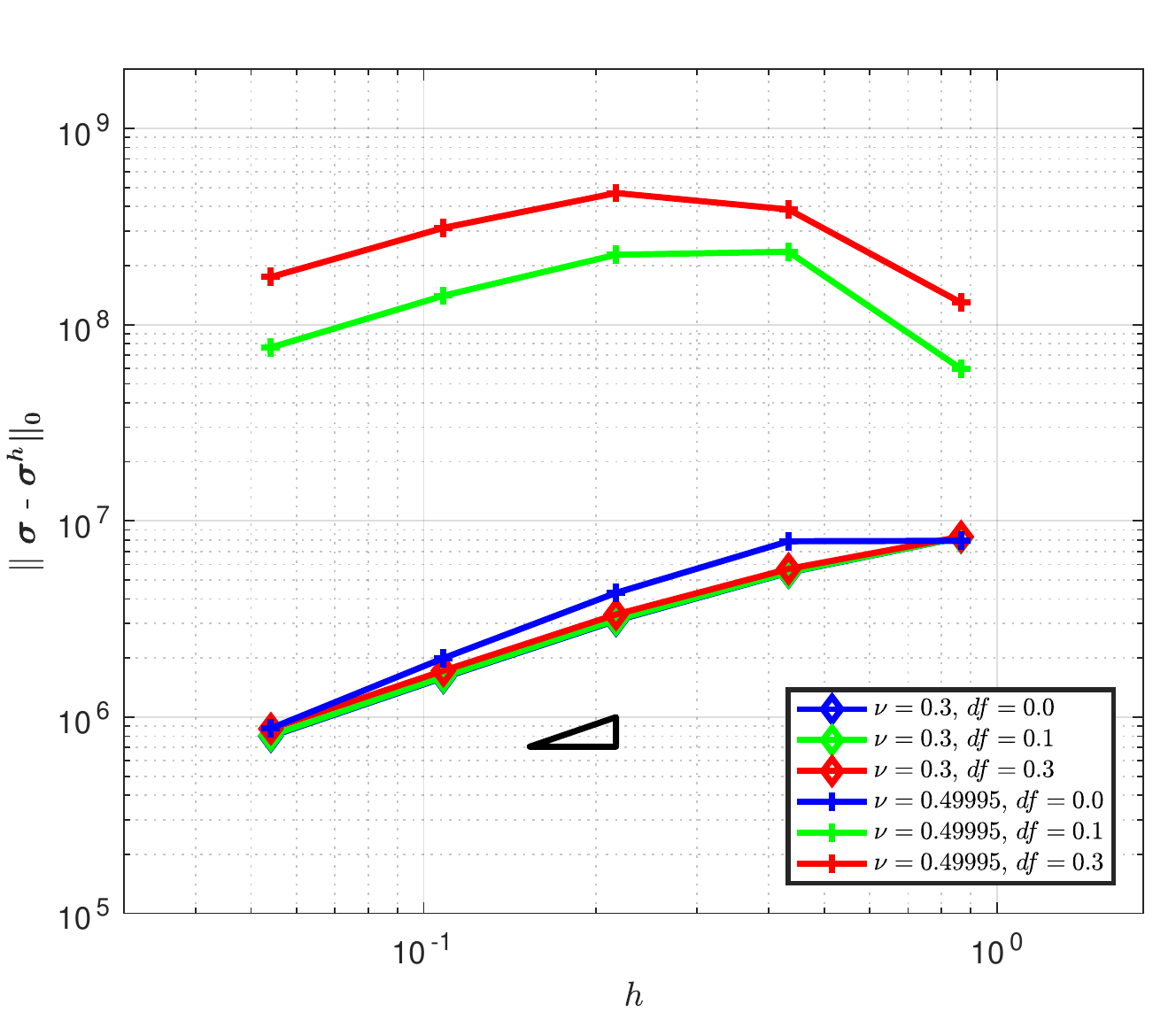}\label{Cube_stress_error_NIPG}}\\
\subfloat[New IIPG]{\includegraphics[trim={3.5cm 9cm 4cm 9cm},clip,width=.40\columnwidth]{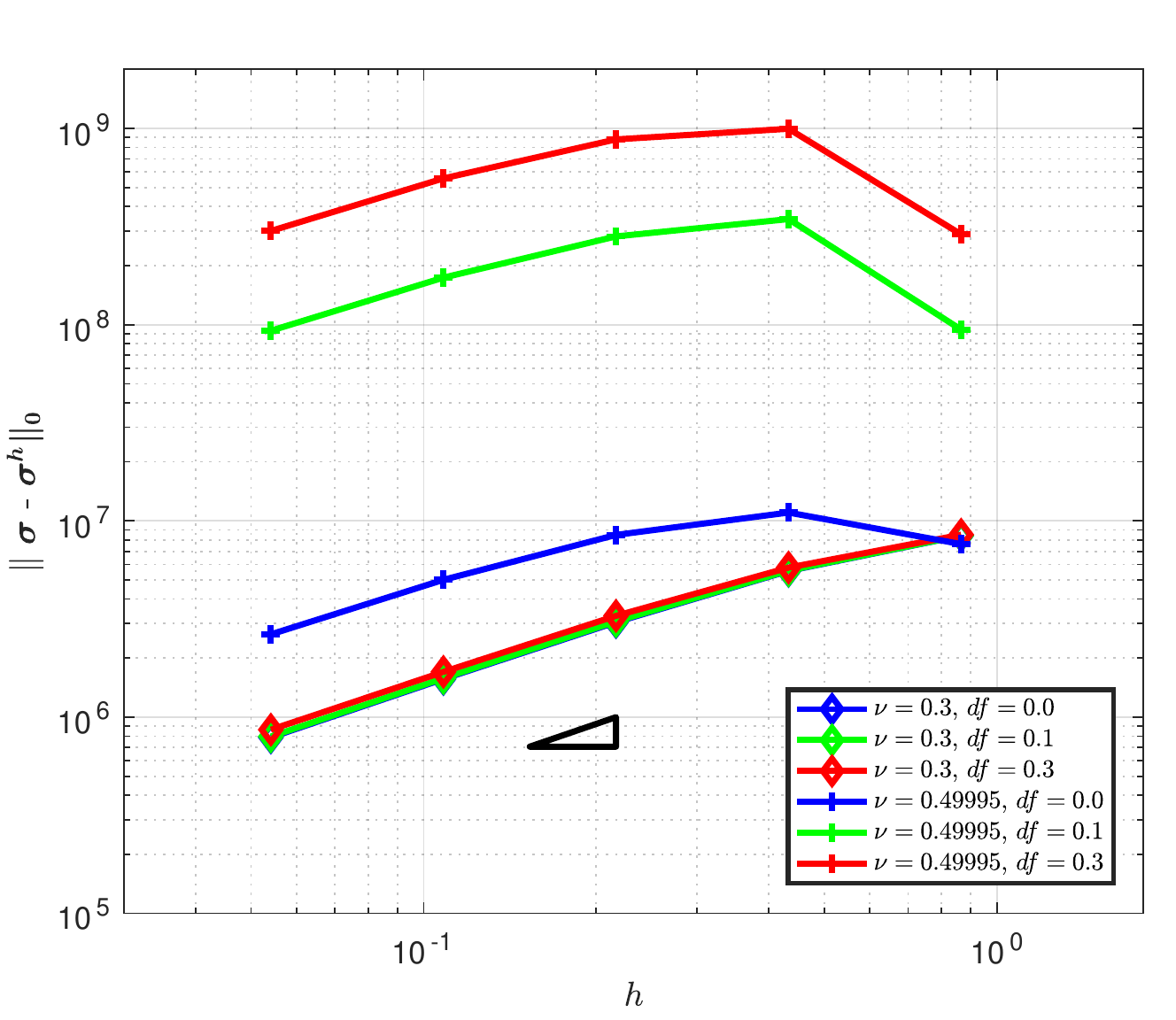}\label{Cube_stress_error_IIPG}}
\subfloat[New SIPG]{\includegraphics[trim={3.5cm 9cm 4cm 9cm},clip,width=.40\columnwidth]{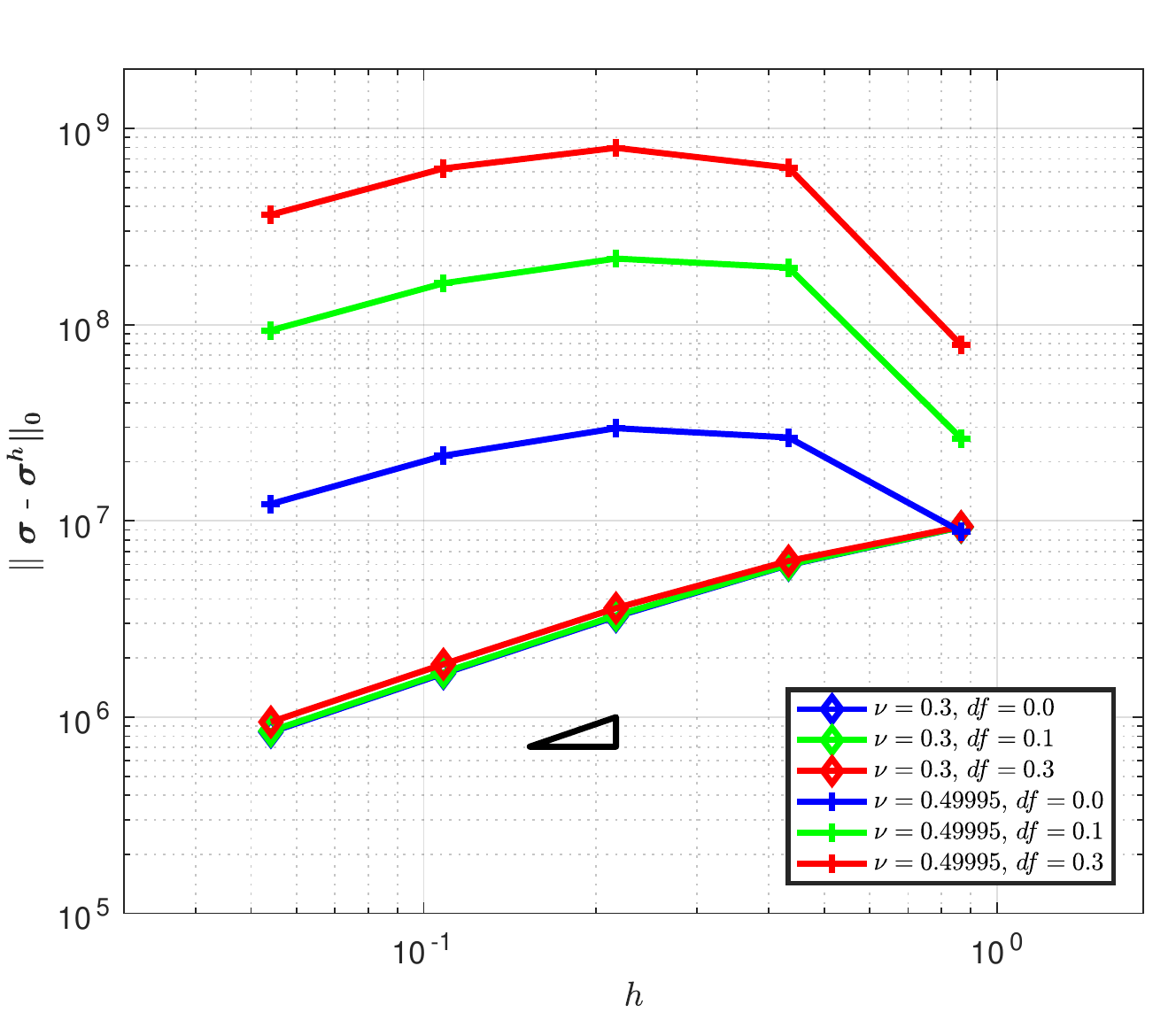}\label{Cube_stress_error_SIPG_50_250}}
\caption{Cube: Stress $L^2$ error convergence. The hypotenuse of the triangle has a slope of 1 in each case.}
\end{figure}


\begin{figure}[!ht]
\centering
\subfloat[$\sigma_{xx}$]{\includegraphics[width=.275\columnwidth]{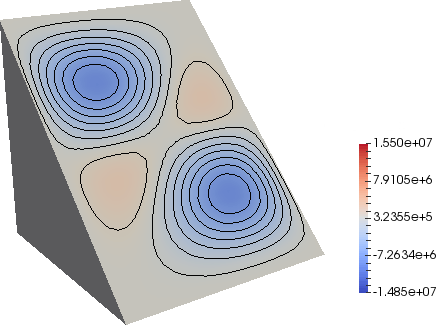}\label{Cube_exact_xx}} \hspace{10mm}
\subfloat[$\sigma_{xz}$]{\includegraphics[width=.275\columnwidth]{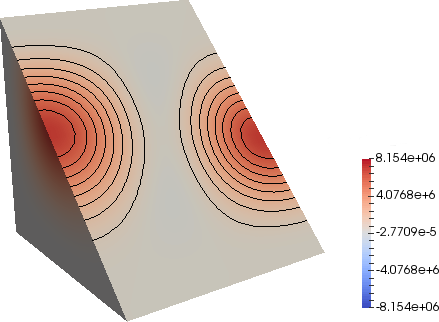}\label{Cube_exact_xz}}
\caption{Cube: Stress, exact solution, $\nu = 0.49995$}\label{Cube_stress_exact}
\end{figure}

\begin{figure}[!ht]
\centering
\subfloat[\textit{df} = 0.0]{\includegraphics[width=.20\columnwidth]{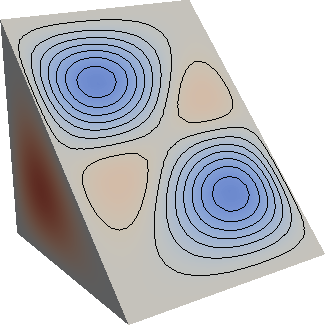}\label{Cube_SIPG_0_xx}} \hspace{10mm}
\subfloat[\textit{df} = 0.1, without contour lines]{\includegraphics[width=0.2\columnwidth]{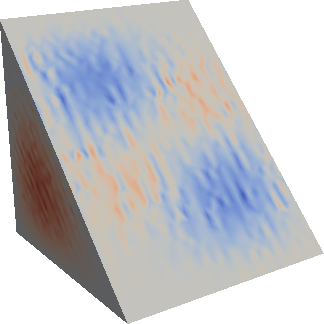}\label{Cube_SIPG_1_xx} }\hspace{10mm}
\subfloat[\textit{df} = 0.3, without contour lines]{\includegraphics[width=.275\columnwidth]{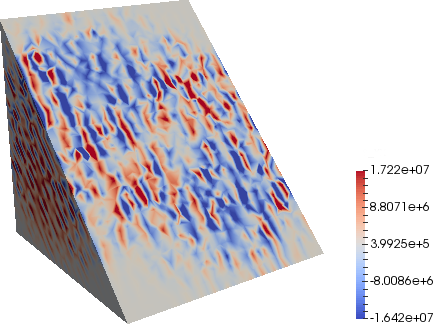}\label{Cube_SIPG_3_xx}}
\caption{Cube: $\sigma_{xx}$, SIPG, $\nu = 0.49995$}\label{Cube_stress_xx_SIPG}
\end{figure}

\begin{figure}[!ht]
\centering
\subfloat[\textit{df} = 0.0]{\includegraphics[width=.2\columnwidth]{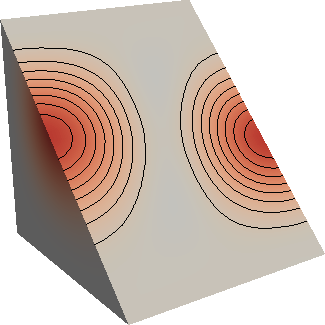}\label{Cube_SIPG_0_xz}} \hspace{10mm}
\subfloat[\textit{df} = 0.3]{\includegraphics[width=.27\columnwidth]{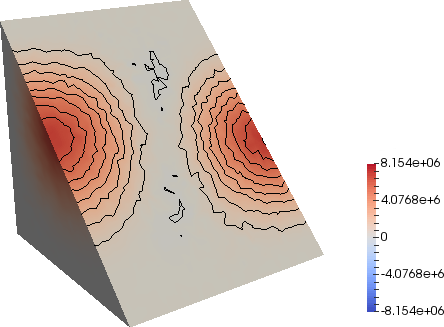}\label{Cube_SIPG_3_xz}}
\caption{Cube: $\sigma_{xz}$, SIPG, $\nu = 0.49995$}\label{Cube_stress_xz_SIPG}
\end{figure}

As in the first two examples, the $L^2$ error of the postprocessed stress of the SG method with $Q_1$ elements shows first-order convergence for all \textit{df} when $\nu = 0.3$, with a slight decline for \textit{df} = 0.3, and diverges when $\nu = 0.49995$ (Figure \ref{Cube_stress_error_SG}). When SRI is applied, first-order convergence is displayed for all \textit{df} for both values of $\nu$, with larger errors for the nearly incompressible case (Figure \ref{Cube_stress_error_SG_sri}). The original IP methods show diverging  stress $L^2$ errors for the nearly incompressible case (Figure \ref{Cube_stress_error_IPorig_nu2}). All three new IP methods tend towards first-order convergence for all \textit{df} when $\nu = 0.3$; when $\nu = 0.49995$ the rates are lower than first-order, declining with increasing \textit{df} and poorest for SIPG, but improve with refinement (Figures \ref{Cube_stress_error_NIPG}, \ref{Cube_stress_error_IIPG} and \ref{Cube_stress_error_SIPG_50_250}). The error magnitudes are larger for the nearly incompressible case than for the compressible case, with decreasing element regularity also resulting in greater errors. The IP errors are, however, smaller than those of the SG with SRI.

The qualitative behaviour of the SIPG method on mesh 5 for the individual components $\sigma_{xx}$ and $\sigma_{xz}$, for $\nu = 0.49995$, is shown in Figures \ref{Cube_stress_xx_SIPG} and \ref{Cube_stress_xz_SIPG}, with the exact solutions shown in Figure \ref{Cube_stress_exact} for comparison. There is a clear deterioration in accuracy in the approximate $\sigma_{xx}$ field as \textit{df} increases: with \textit{df} = 0.0, the result is visually accurate, with \textit{df} = 0.1 there is a loss of smoothness of the field, and with \textit{df} = 0.3 the field is significantly more patchy, although the general areas of high and low stresses are still reflected. The coarser the meshes, the less clearly these zones are displayed (images not shown here), indicating that refinement improves the quality of the postprocessed field in this component. The results for the direct stresses $\sigma_{yy}$ and $\sigma_{zz}$ are similar. The approximate shear stress component $\sigma_{xz}$ (Figure \ref{Cube_stress_xz_SIPG}) shows comparatively little deterioration in accuracy, though some loss of smoothness, as \textit{df} increases, as do the shear components $\sigma_{xy}$ and $\sigma_{yz}$. 

The other two IP methods similarly display increasing ``patchiness" in the direct stresses as \textit{df} increases, and similarly maintain good overall approximations, though with loss of smoothness as \textit{df} increases, in the shear components.

\section{Conclusion}\label{conc}

The first and primary aim of this paper was to establish that the uniform convergence, with respect to the compressibility parameter, of the new IP methods of \cite{Grieshaber2015} extends to the use of general quadrilateral/hexahedral elements. In all four model problems presented here, this uniform convergence is indeed seen. 
While the IP approximation errors are slightly larger in many cases for greater degrees of distortion, this is not a feature of near-incompressibility only: this behaviour is seen in the case $\nu = 0.3$ as well, and is also displayed by the SG method for $\nu = 0.3$.
Contour plots of the displacement approximations indicate that the use of general quadrilateral/hexahedral, instead of rectangular, elements does not visibly diminish the quality of the displacement results.

The second aim was to ascertain the extent to which a stress field postprocessed from the IP displacement approximation is accurate, particularly in the near-incompressible case and for non-rectangular elements. We see in the examples of the cantilever beam and the square plate that, for the stress approximation, the first-order convergence rate of the $L^2$ error achieved for the case of $\nu = 0.3$ using all three IP methods is also achieved for the case of $\nu = 0.49995$, irrespective of element shape, except in the case of the NIPG method for the square plate, where the error is nevertheless smaller than for the other IP methods. In the example of the L-shaped domain, the convergence rate is lower than first-order for $\nu = 0.3$ and this rate is exceeded for $\nu = 0.49995$. For the example of the cube, convergence is first-order when $\nu = 0.3$, and for $\nu = 0.49995$, convergence tends towards or exceeds this rate. These examples strongly suggest that the stress error of the new IP methods converges uniformly in the $L^2$-norm with respect to the compressibility parameter.

To assess the quality of the stress approximation on a given mesh (i.e.\ at a given refinement level), a continuous field was obtained for visualisation by a projection of the stress values calculated at quadrature points onto the mesh nodes, and compared component-wise to the analytical stress expression through contour plots. 
In general, with $\nu = 0.49995$ and general quadrilateral/hexahedral elements, 
IP results in the direct stress components were of lower quality than with rectangular elements, and more evidently so the lower the element shape regularity, while in the shear stress components the approximate fields lost only smoothness. With sufficient refinement, the quality was recovered, but the refinement levels necessary for good-quality stress approximations corresponded in some cases to very large systems of equations. The displacement approximations typically give visually accurate contour plots at several refinement levels lower than the postprocessed stresses do, for a given problem. 

This phenomenon is not, however, unique to the IP methods. For example, the SG method with $Q_2$ elements likewise produces low-accuracy stress fields in the near-incompressible case at a refinement level at which the displacements are extremely good, as the error plots show. Poor performance in the stress in the near-incompressible case is, moreover, not limited to the case of non-rectangular elements. Obtaining accurate results from postprocessing would thus seem to require alternative strategies, even where error convergence rates are very good.

A significant addition  and useful complement to the computational investigation presented here would be a theoretical displacement error analysis of the new IP formulation for the case of general quadrilateral elements. More generally, there is scope for an extension of this formulation and the corresponding analyses to nonlinear problems.
Regarding a numerical, implementational aspect of this work, direct solvers were used in the examples in this study, but are computationally inefficient for large systems. An investigation into the design of an appropriate preconditioner to be used with an iterative solver in the near-incompressible regime, with the aim of decreasing computing time and memory usage, would therefore be of value.

\section*{Acknowledgements}

The authors acknowledge with thanks the support for this work by the National Research Foundation, through the South
African Research Chair in Computational Mechanics. AM gratefully acknowledges the support provided by the EPSRC Strategic Support Package: Engineering of Active Materials by Multiscale/Multiphysics Computational Mechanics - EP/R008531/1.

\bibliographystyle{plain}
\bibliography{GrieshaberMcBrideReddy_Oct2019}

\end{document}